\documentclass[11pt]{article}

\usepackage[utf8]{inputenc}
\usepackage{amssymb}
\usepackage{mathrsfs}

\usepackage{times}
\usepackage{amsmath}
\usepackage{amsthm}
\usepackage{amsfonts} 
\usepackage{color}
\usepackage[colorlinks,linkcolor=blue,anchorcolor=blue,citecolor=blue]{hyperref}
\usepackage{bbm}
\usepackage{tikz}
\usetikzlibrary{arrows}

\usepackage[top=1in, bottom=1in, left=1.25in, right=1.25in]{geometry}

\usepackage{titlesec}
       \titleformat{\chapter}[display]
             {\normalfont\Large\bfseries}{\thechapter}{11pt}{\Large}
       \titleformat{\section}
             {\normalfont\large\bfseries}{\thesection}{11pt}{\large}
       \titlespacing*{\chapter}{0pt}{0pt}{15pt} %left, beforesep, aftersep, right
       \titlespacing*{\section}{0pt}{3.5ex plus 1ex minus .2ex}{2.3ex plus .2ex}

    \usepackage[titletoc]{appendix}

\usepackage[all,cmtip]{xy}

\newtheorem{thm}{Theorem}[section]
\newtheorem{prop}[thm]{Proposition}
\newtheorem{lemma}[thm]{Lemma}
\newtheorem{cor}[thm]{Corollary}
\newtheorem{definition}[thm]{Definition}%[section]

\newtheorem{remark}[thm]{Remark}

\newtheorem{prop-defn}[thm]{Proposition/Definition}

\newcommand{\abs}[1]{\lvert#1\rvert}

 \newcommand{\ba}{\begin{eqnarray}}
   \newcommand{\na}{\end{eqnarray}}
   \newcommand{\ban}{\begin{eqnarray*}}
   \newcommand{\nan}{\end{eqnarray*}}

\newcommand{\bC}{{\mathbb C}}

\newcommand{\bP}{{\mathbb P}}

\newcommand{\bR}{{\mathbb R}}
\newcommand{\bZ}{{\mathbb Z}}

\newcommand{\cF}{{\mathcal F}}

\newcommand{\cO}{{\mathcal O}}
\newcommand{\cS}{{\mathcal S}}
\newcommand{\cZ}{{\mathcal Z}}

\newcommand{\Arg}{\mathrm{Arg}}
\newcommand{\id}{\mathrm{id}}
\newcommand{\res}{\mathrm{Res}}

\newcommand{\Eff}{\mathrm{Eff}}
\newcommand{\Ch}{\mathrm{Ch}}
\newcommand{\End}{\mathrm{End}}

\newcommand{\hpc}{h}   %hyperplane class
\newcommand{\qdc}{\mathcal{Q}}  %quadric
\newcommand{\qtp}{q} % quantum parameter

  \newcommand{\<}{\langle}
  \renewcommand{\>}{\rangle}
\newcommand{\suml}{\sum\limits}
\newcommand{\prodl}{\prod\limits}

\newtheorem*{gammaII}{Gamma conjecture II}

\begin{document}{\allowdisplaybreaks[4]

\title{GAMMA CONJECTURE II FOR QUADRICS}

\author{XIAOWEN HU, AND HUA-ZHONG KE}
%\address{School of Mathematics, Sun Yat-sen University, Guangzhou 510275, P.R. China}

%\curraddr{Department of Mathematics }
%\email{huxw06@gmail.com}
%\thanks{ %The first author  is supported in part by %a RGC research grant from the Hong Kong Government.
%The second author  is  supported in part %by KRF-2007-341-C00006.
 %}

%\author{Hua-Zhong Ke}
%\address{School of Mathematics, Sun Yat-sen University, Guangzhou 510275, P.R. China}
%\curraddr{Department of Mathematics }
%\email{kehuazh@mail.sysu.edu.cn}
%\thanks{ %The first author  is supported in part by %a RGC research grant from the Hong Kong Government.
%The second author  is  supported in part %by KRF-2007-341-C00006.}

\date{}
%\author{Changzheng Li}
 %\address{School of Mathematics, Sun Yat-sen University, Guangzhou 510275, P.R. China}
%\email{lichangzh@mail.sysu.edu.cn}
\maketitle

%\thanks{2010 Mathematics Subject Classification. Primary 14N35. Secondary 14J45, 14J33.}
%\date{%January 1, 2001 and, in revised form, June 22, 2001.}


%\keywords{Gamma conjecture II. Fano manifolds. Quadratic hypersurfaces.}

%\subtitle{}

\begin{abstract}
%Gamma conjecture II for a Fano manifold $X$ was proposed by Galkin-Golyshev-Iritani, under the assumption of convergence and semisimplicity of quantum cohomology of $X$ and the existence of full exceptional colletions in $\mathcal D^b(X)$. This conjecture expects to describe the asymptotic behavior of flat sections of Dubrovin connection near irregular singularities in terms of a full exceptional collection in $\mathcal D^b(X)$ via the $\widehat{\Gamma}$-integral structure. The main results of this note are the followings: (i) the full quantum cohomology of a smooth quadric is convergent; (ii) Gamma conjecture II holds for smooth quadrics.
The Gamma conjecture II for the quantum cohomology of a Fano manifold $F$, proposed by Galkin, Golyshev and Iritani, describes the asymptotic behavior of the flat sections of the Dubrovin connection near the irregular singularities, in terms of a full exceptional collection, if there exists, of $\mathcal D^b(F)$ and the $\widehat{\Gamma}$-integral structure. In this paper, for the smooth quadric hypersurfaces  we prove the convergence of the full quantum cohomology  and the Gamma conjecture II. 

For the proof, %firstly, based on Galkin, Golyshev and Iritani's result that certain flat sections are Laplace tranforms of line integrals,} 
we first give a criterion on Gamma II  for Fano manifolds with semisimple quantum cohomology, by Dubrovin's theorem of analytic continuations of semisimple Frobenius manifolds.
Then we work out a closed formula of  the Chern characters of spinor bundles on quadrics. 
By the deformation-invariance of Gromov-Witten invariants we show that the full quantum cohomology can be reconstructed by its ambient part, and use this to obtain estimations. Finally we complete the proof of Gamma II for quadrics by  explicit asymptotic expansions of flat sections corresponding to Kapranov's exceptional collections and an application of our criterion.
\end{abstract}

\maketitle
  \tableofcontents

\section{Introduction}\label{sec:introduction}

A Fano manifold is a compact, complex manifold with ample anti-canonical line bundle. For the quantum cohomology of Fano manifolds, Galkin, Golyshev and Iritani proposed Conjecture $\cO$ and Gamma conjectures (I and II) \cite{GGI16}. Both Gamma conjecture I and II concern the asymptotic behavior of flat sections of Dubrovin connection near irregular singularities, while Conjecture $\cO$ is an underlying assumption for Gamma conjecture I. In this paper, we focus on Gamma conjecture II, and we refer readers to \cite{BFSS,CL17,GaIr19,GZ16,HKLY,Ke181,LMS17,SS17,Wit18} for recent progress of Conjecture $\cO$ and Gamma conjecture I.

Assuming the convergence of the quantum cohomology, the Dubrovin connection of a Fano manifold $F$ is a PDE system over $\bP^1\times H^{\textrm{even}}(F)$. We can roughly view it as a family of linear ODE systems on $\bP^1$ parameterized by $\mathbf t\in H^{\textrm{even} }(F)$. Let $z$ be an in-homogeneous coordinate on $\bP^1$. For each $\mathbf t$, the corresponding OED system has only two singularities: the regular singularity at $z=\infty$ and the irregular singularity at $z=0$. Let $\cS_\mathbf t$ be the space of global (multi-valued) solutions of the ODE system corresponding to $\mathbf t$. Inspired by mirror symmetry, Iritani used a canonical fundamental solution near $z=\infty$ to define the $K$-group framing on the space $\cS_\mathbf t$, which is a group-homomorphism from $K(F)$ to $\cS_\mathbf t$. The $\widehat\Gamma$-integral structure is the image of the $K$-group framing, which is a full-rank lattice on $\cS_\mathbf t$. Now we assume furthermore the semisimplicity of quantum cohomology. Then the Dubrovin connection has asymptotically exponential fundamental solutions (AEFS for short). Let the parameter $\mathbf t$ be a semisimple point of the quantum cohomology. For a suitably chosen $\phi\in\bR$, an AEFS corresponding to $(\mathbf t,\phi)$ is a fundamental solution in $\cS_\mathbf t$ consisting of flat sections characterized by their asymptotic behavior near $z=0$ along the direction $\arg z=\phi$. Loosely speaking, assuming the existence of  full exceptional collections in the bounded derived category $\mathcal D^b(F)$ of coherent sheaves on $F$, Gamma conjecture II states that an AEFS corresponding to $(\mathbf t,\phi)$ lies in Iritani's $\widehat\Gamma$-integral structure of $\cS_\mathbf t$, and moreover, it is given by a suitably chosen full exceptional collection of $\mathcal D^b(F)$. We refer readers to Section \ref{sec:preliminaries} for a precise statement of the conjecture.

Gamma conjecture II implies a part of the original Dubrovin's conjecture \cite{Dub98}, of which the refined formulation was proposed in \cite{Dub13,CDG18}.
 %We remark that Gamma conjecture II 
 It is closely related to mirror symmetry, and we refer readers to introductions of \cite{GGI16,GaIr19} for more discussions.

%A series of conjectures relating the quantum cohomology of a Fano manifold to its derived category was proposed by Dubrovin \cite{Dub98}, which were refined in \cite{Dub13, CDG18}. 
%The above mentioned Dubrovin's conjecture has various refinements in \cite{Dub13, HMT09, CDG18}. It was shown in \cite[Theorem 1.3]{CDG18} that Gamma conjecture II is equivalent to the part 3.b of their version of refined Dubrovin's conjectures \cite[Conjecture 1.1]{CDG18}. 

Presently the Gamma conjecture II is proved in only a few cases (Grassmannians \cite{CDG18,GGI16} and toric Fano manifolds \cite{FaZh19}). 
%Both cases contain projective spaces. Recall that a Fano manifold of dimension $N$ has at most Fano index $N+1$; it has Fano index $N+1$ iff it is a projective space, and it has Fano index $N$ iff it is a smooth quadric, i.e. a smooth quadric hypersurface in a projective space. 
The main result of this paper is the following.
\begin{thm}(= Theorem \ref{thm-GammaII-quadric})
 The Gamma conjecture II holds for the smooth quadric hypersurfaces in the projective spaces.
\end{thm}

In the rest of this paper, a smooth quadric hypersurface will be called a quadric for short. The strategy of the proof is to build an AEFS from a full exceptional collection via $\widehat\Gamma$-integral structure. We start with Kapranov's full exceptional collections for quadrics, and compute their Chern characters by using some knowledge of the spinor bundles. Then, we apply Givental's mirror theorem to compute the corresponding flat sections in the $\widehat\Gamma$-integral structure of $\cS_\mathbf 0$ near $z=\infty$, where $\mathbf 0$ is the origin of the cohomology space. It turns out that these flat sections can be expressed by Meijer's $G$-functions, and this enables us to find asymptotic expansions of these flat sections near $z=0$. However, Kapranov's full exceptional collection does not directly give an AEFS in $\cS_\mathbf 0$, since these asymptotic expansions are obtained along different directions. Similar difficulty also appears in Galkin, Golyshev and Iritani's proof of the Gamma conjecture II for projective spaces, where they solved this problem by showing that suitable mutations of  Beilinson's exceptional collections give AEFS in $\cS_{\mathbf0}$ (\cite[Section 5.3]{GGI16}). 
%They also pointed out that, after isomonodromic deformation, Beilinson's exceptional collections give AEFS, though not in $\cS_\mathbf0$ \cite[Remark 5.3.3]{GGI16}. 
Inspired by  \cite[Remark 5.3.3]{GGI16}, we establish a criterion for the Gamma conjecture II for general Fano manifolds (Theorem \ref{sufficientcondition}). Then  Gamma conjecture II for quadrics follows by an application of this criterion. 

Probably Theorem \ref{sufficientcondition} has independent interests. Roughly speaking, 
it states the following: assume that for a fundamental solution $\{y_i(\mathbf t,z)\}$ of the Dubrovin connection (instead of a single ODE system),  there exist $\mathbf{t}_0\in H^{\textrm{even}}(F)$ and real numbers $\phi_i$'s, such that, among other assumptions  not spelled out here, for each $i$, $y_i(z):=y_i(\mathbf{t}_0,z)\in\cS_{\mathbf t_0}$ has a suitable asymptotic behavior near $z=\infty$ along the direction $\arg z=\phi_i$; then there exist $\mathbf{t}_1$ and a real number $\phi$, such that $\{\bar y_i(z):=y_i(\mathbf{t}_1,z)\in\cS_{\mathbf t_1}\}$ is an AEFS corresponding to $(\mathbf{t}_1,\phi)$. The point is that the final $\phi$ is common for all $\bar y_i(z)$'s. We give a brief account of the proof of Theorem \ref{sufficientcondition}, and refer readers to Section \ref{sec:sufficientCondition-GammaII} for details. Firstly we can conclude from Lemma \ref{respectextension} that each $y_i$ is in an AEFS corresponding to $(\mathbf t_0,\phi_i)$. Results of \cite[Section 2.5]{GGI16} imply that a flat section in an AEFS in $\cS_{\mathbf t}$ is the Laplace transform of a flat section of the dual connection, and vice versa. The Laplace transform (see \eqref{laplacetransform}) is given by a line integral on $\bC$ along a half line starting at a singularity of the dual connection, and the singularities are analytic (multivalued) functions of the parameter $\mathbf t$. This enables us to pictorially view $y_i$ as given by a half line $L_i$ on $\bC$ starting at a singularity $u_i$ with phase $\phi_i$. The idea is to move the parameter $\mathbf t$ starting from $\mathbf t_0$, which corresponds to moving the singularities of the dual connection. At the same time, we need to move the half lines starting at these singularities, and vary their directions continuously if necessary to guarantee that the half lines do not intersect in the moving process, so that at the end we can make these half lines parallel. At first glance, the moving of the parameter $\mathbf t$ is only possible in a neighborhood of $\mathbf t_0$ due to the convergence issue. A result of Dubrovin \cite[Section 4]{Dub00} states that a semisimple Frobenius manifold can be meromorphically continued to the universal covering of the configuration space $\mathcal C$ of singularities of the dual connection. So we can move the parameter $\mathbf t$ from $\mathbf t_0\in H^{\textrm{even}}(F)$ to a suitably chosen $\mathbf t_1\in\mathcal C$ to prove Theorem \ref{sufficientcondition}.

An important assumption in Theorem \ref{sufficientcondition} is the convergence of quantum cohomology. The convergence of the generating functions of Gromov-Witten invariants is a fundamental problem and is unsolved in the general cases. We verify it for the quadrics (Theorem \ref{convergencofquadrics}). The convergence of the ambient part of the quantum cohomology of the quadrics was known \cite{Iri07, Zin14}. We prove the convergence for invariants with primitive classes. Note that on the cohomology of a smooth quadric there is a $\mathbb{Z}/2 \mathbb{Z}$ monodromy action, which is induced by the whole family of the smooth quadrics.
%that smooth quadrics have natural monodromy action on the cohomology, and 
We use the deformation invariance of Gromov-Witten invariants to simplify their WDVV equations, which helps us to reconstruct the full quantum cohomology from its ambient part. This idea was used by the first author to determined the full quantum cohomology of the smooth cubic hypersurfaces \cite{Hu15}, and also by the second author to verify Conjecture $\cO$ for the projective complete intersections \cite{Ke181}. 

While we were writing this paper, we noticed that Theorem \ref{convergencofquadrics} can be deduced from a recent result of Cotti \cite[Theorem 6.6]{Cot20}. Cotti's theorem assumes the semisimplicity of the small quantum cohomology, which is valid for smooth quadrics. However we expect that our method can be applied to show the convergence in other circumstances, e.g. for  the projective complete intersections (the quantum cohomology of a projective complete intersection, other than the quadrics, the complete intersection of two even dimensional quadrics or the cubic surfaces, is not semisimple).

The rest of the paper is organized as follows. In Section \ref{sec:preliminaries}, we review some background materials and give a precise statement of Gamma conjecture II. In Section \ref{sec:sufficientCondition-GammaII}, we prove the above-mentioned criterion for Gamma conjecture II (Theorem \ref{sufficientcondition}). In Section \ref{sec:chernclass-spinorbundle}, we compute the Chern characters of the spinor bundles on quadrics, using Riemann-Roch and an inductive feature of the spinor bundles. In Section \ref{sec:convergence}, we show that  the primary genus zero Gromov-Witten invariants of the smooth quadrics involving primitive classes are determined by the 3-points invariants and the ambient invariants. Then by some inductive estimations  we verify a conjecture of Zinger \cite[Conjecture 1]{Zin14} in the case of  quadrics for the genus zero Gromov-Witten invariants without descendents  and thus prove the convergence of the quantum cohomology of them. In Section \ref{sec:proof-GammaII}, with  the help of the result in  Section \ref{sec:chernclass-spinorbundle}  we find expressions of the flat sections of the Dubrovin connection, corresponding to Kapranov's exceptional collections, in terms of  Meijer's $G$-functions, then by the asymptotic expansions of the latters and theorem \ref{sufficientcondition}  we prove the Gamma conjecture II for quadrics.

\section{Preliminaries}\label{sec:preliminaries}

%In this section, we assume that $F$ is a Fano manifold with vanishing odd cohomology, and set $H^*(F)=H^*(F,\bC)$. We also assume that $\dim H^*(F)=s$.
For a manifold $F$ we denote $H^*(F)=H^*(F,\bC)$. 
In this section, for simplicity we assume that $F$ is a Fano manifold with $H^{\textrm{odd}}(F)=0$, and denote $s=\dim H^*(F)$.

\subsection{Quantum cohomology and Dubrovin connection}

Let $\Eff(F)\subset H_2(F,\bZ)$ be the set of effective curve classes in $F$. For $\mathbf d\in \Eff(F)$, let $\overline{M}_{0,n}(F,\mathbf d)$ be
the moduli space of connected $n$-pointed stable maps from connected nodal curves of arithmetic genus zero to $F$ with degree $\mathbf d$. Let $e_i$ be the evaluation map at the $i$-th marked point. The (genus-zero, descendent) Gromov-Witten invariants of $F$ are defined as
$$
\langle\prodl_{i=1}^n \tau_{k_i}(\gamma_i)\rangle^F_{\mathbf d}
:=\int_{[\overline{M}_{0,n}(F,\mathbf d)]^{vir}}\prod\limits_{i=1}^n\psi_i^{k_i}e_i^*\gamma_i,
$$
where $\gamma_i\in H^*(F)$, $k_i\in\bZ_{\geqslant 0}$, $\psi_i$ is the first Chern class of the cotangent line bundle associated to the $i$-th marked point, and $[\overline{M}_{0,n}(F,\mathbf d)]^{vir}$ is the virtual fundamental class. Informally speaking, if $\gamma_i$'s are algebraic, then the primary invariant $$\<\gamma_1,\dots,\gamma_n\>_{\mathbf d}^F:=\<\prodl_{i=1}^n \tau_{0}(\gamma_i)\>_{\mathbf d}^F$$ is the (virtual) number of degree-$\mathbf d$ rational curves in $F$ intersecting subvarieties Poincar\'e dual to $\gamma_1,\dots,\gamma_n$ in general position. We refer readers to \cite[Section 7.3, Section 10.1]{CoKa99} or \cite[Chapter V-VI]{Man99} for general properties of Gromov-Witten invariants.

Let $T_0:=\mathbbm1,T_1,\dots,T_{s-1}$ be a homogeneous basis of $H^*(F)$. For $\mathbf t=\suml_{i=0}^{s-1}t^iT_i\in H^*(F)$, we formally define the (genus-zero, primary) potential function of $F$ by
\ban
\cF_0^F(\mathbf t):=\suml_{n=0}^\infty\frac1{n!}\suml_{\mathbf d\in\Eff(F)}\<\tau_0(\mathbf t)^n\>^F_{\mathbf d}.
\nan
By  the point mapping axiom and the fundamental class axiom, we can formally write
\ban
\cF_0^F(\mathbf t)&=&\suml_{n_0+n_1+\cdots+n_{s-1}=3}\bigg(\int_F\prodl_{i=0}^{s-1}T_{i}^{n_{i}}\bigg)\prodl_{i=0}^{s-1}\frac{(t^{i})^{n_{i}}}{n_{i}!}\\
&&\qquad+\suml_{n_{1},\dots,n_{s-1}\geq0}\bigg(\suml_{\mathbf d\in\Eff(F)\setminus\{0\}}\<\prodl_{i=1}^{s-1}\tau_0(T_{i})^{n_{i}}\>_{\mathbf d}^F\bigg)\prodl_{i=1}^{s-1}\frac{(t^{i})^{n_{i}}}{n_{i}!}.
\nan
Since $F$ is Fano, it follows from the degree axiom and \cite[Corollary 1.19]{KoMo98} that the coefficient of $\prodl_{i=1}^{s-1}\frac{(t^{i})^{n_{i}}}{n_{i}!}$ is actually a finite sum. So we can see that
\ban
\cF_0^F(\mathbf t)\in\bC[[t^0,t^1,\dots,t^{s-1}]].
\nan
For each $\mathbf t\in H^*(F)$, the quantum product $\bullet_\mathbf t$ on $H^*(F)$ parameterized by $\mathbf t$ is formally given by
\ban
\<T_i\bullet_\mathbf t T_j,T_k\>^F&:=&\partial_{t^i}\partial_{t^j}\partial_{t^k}\cF_0^F(\mathbf t)\\
&=&\int_{F}T_i\cup T_j\cup T_k+\mbox{correction terms}.
\nan
Here $\<\cdot,\cdot\>^F$ is the Poincar\'e pairing on $H^*(F)$. Let $g_{pq}:=\<T_p,T_q\>^F$, and let $(g^{pq})$ be the inverse matrix of $(g_{pg})$. Then for any $i,j,k,l\in\{0,1,\dots,s-1\}$, $\cF_0^F(\mathbf t)$ satisfies the famous WDVV equation:
\begin{equation}
\suml_{p,q}g^{pq}\partial_{t^i}\partial_{t^j}\partial_{t^p}\cF_0^F(\mathbf t)\cdot\partial_{t^q}\partial_{t^k}\partial_{t^l}\cF_0^F(\mathbf t)=\suml_{p,q}g^{pq}\partial_{t^i}\partial_{t^k}\partial_{t^p}\cF_0^F(\mathbf t)\cdot\partial_{t^q}\partial_{t^j}\partial_{t^l}\cF_0^F(\mathbf t).\label{generalWDVV}
\end{equation}

Assume that the quantum product $\bullet_\mathbf t$ is analytic for $\mathbf t$ in a neighborhood $B$ of $\mathbf 0\in H^*(F)$, i.e. $\cF_0^F(\mathbf t)$ is an absolutely convergent power series of $t^0,t^1,\dots,t^{s-1}$ for $\mathbf t\in B$. Then the quantum cohomology $(H^*(F),\bullet_\mathbf t)$ is a commutative, associative $\bC$-algebra with the unit element $\mathbbm1\in H^*(F)$ for each $\mathbf t\in B$ (the associativity comes from WDVV). The \emph{Dubrovin connection} \cite{Dub96,Dub98,Dub00} of $F$ is a meromorphic flat connection $\nabla$ on the trivial $H^*(F)$-bundle over $B\times\bP^1$ given as follows: for a point $(\mathbf t,z)\in B\times\bP^1$, we have
\ban
\nabla_{\partial_{t^i}}&=&\partial_{t^i}+\frac1z(T_i\bullet_\mathbf t),\\
\nabla_{z\partial_z}&=&z\partial_z-\frac1z(E\bullet_\mathbf t)+\mu.
\nan
Here $E=c_1(F)+\sum_{i=0}^{s-1}(1-\frac12\deg T_i)t^iT_i$ is the Euler vector field, and $\mu\in\End(H^*(F))$ is the Hodge grading operator given by $\mu|_{H^{2p}(F)}:=(p-\frac{\dim F}2)\id_{H^{2p}(F)}$. For each $\mathbf t\in B$, denote by $\nabla_{\mathbf t}$ the restriction of $\nabla$ to $\{\mathbf t\}\times\bP^1$. Then $\nabla_{\mathbf t}$ is a meromorphic flat connection on the trivial $H^*(F)$-bundle over $\bP^1=\{\mathbf t\}\times\bP^1$, given by
\[
(\nabla_{\mathbf t})_{z\partial_z}=z\partial_z-\frac1z(E\bullet_\mathbf t)+\mu.
\]
Dubrovin proved that $\nabla$ is the isomonodromic deformation of $\nabla_\mathbf t$ over $B$ (\cite[Theorem 2.1, 4.4 and 4.6]{Dub00}).

The equation $\nabla f=0$ for a section $f$ is called the \emph{quantum differential equation} of $F$. The space of global solutions to the quantum differential equation is
\ban
\cS_B:=\{f\in\Gamma(B\times\widetilde{\bC^*},\underline{H^*(F)})):\widetilde\nabla f=0\},
\nan
where $\widetilde{\bC^*}$ is the universal covering of $\bC^*$, $\underline{H^*(F)}$ is the trivial $H^*(F)$-bundle over $B\times\widetilde{\bC^*}$, and $\widetilde\nabla$ is the pullback of $\nabla$. The space of global flat sections of $\nabla_\mathbf{t}$ is
\ban
\cS_\mathbf{t}:=\{f\in\Gamma(\{\mathbf t\}\times\widetilde{\bC^*},\underline{H^*(F)}):\widetilde{\nabla_\mathbf t}f=0\},
\nan
where $\underline{H^*(F)}$ is the trivial $H^*(F)$-bundle over $\{\mathbf t\}\times\widetilde{\bC^*}$, and  $\widetilde{\nabla_\mathbf t}$ is the pullback of $\nabla_{\mathbf t}$.

\subsection{Flat sections around \texorpdfstring{$z=\infty$}{z=infty}}

%\textcolor{red}{$\widehat\Gamma$-integral structure of $\cS_{B}$ (or $\cS_{\mathbf t}$) was defined by Iritani \cite{Iri09}, which can be given explicitly near $z=\infty$.}

For $\mathbf t\in H^*(F)$, write $\mathbf t=(\mathbf t^{(2)},\mathbf t')\in H^2(F)\oplus\bigoplus_{p\neq1}H^{2p}(F)$, and we define a formal linear endomorphism $\mathcal L(\mathbf t,z)\in\End(H^*(F))[[\frac1z]]$ parameterized by $\mathbf t$ as follows:
\ban
&&\<\mathcal L(\mathbf t,z)T_i,T_j\>^F\\
&:=&\<e^{-\frac{\mathbf t^{(2)}}z}T_i,T_j\>^F+\suml_{n,m=0}^\infty(\frac{-1}z)^{m+1}\frac1{n!}\suml_{\mathbf d\in\Eff(F)}\<\tau_m(e^{-\frac{\mathbf t^{(2)}}z}T_i)\tau_0(\mathbf t')^n\tau_0(T_j)\>_{\mathbf d}^Fe^{\int_{\mathbf d}\mathbf t^{(2)}}.
\nan

Suppose that the quantum product $\bullet_\mathbf t$ is analytic for $\mathbf t$ in a neighborhood $B$ of $\mathbf0\in H^*(F)$. For any $M\in\End(H^*(F))$, we write $z^M:=e^{M\log z}$. Let $\rho:=(c_1(F)\cup)\in\End(H^*(F))$. Then we can use $L(\mathbf t,z)z^{-\mu}z^\rho$ to identify $H^*(F)$ with the space $\cS_B$. More precisely, the \emph{cohomology framing}, defined by
\ban
\cZ^{coh}_B:H^*(F)&\to&\cS_B\\
\alpha&\mapsto&\mathcal L(\mathbf t,z)z^{-\mu}z^\rho\alpha,
\nan
is a linear isomorphism. For $\mathbf t\in B$, let $\cZ^{coh}_\mathbf t:H^*(F)\to\cS_\mathbf t$ be the restriction of $\cZ^{coh}_B$ to $\{\mathbf t\}\times\widetilde{\bC^*}$. Then $\cZ^{coh}_\mathbf t$ is also a linear isomorphism. As a consequence, the natural restriction $\cS_B\to\cS_\mathbf t$ is a linear isomorphism of $s$-dimensional vector spaces.

Let $K(F)$ be the Grothendieck group of topological complex vector bundles on $F$. The \emph{K-group framing} is a homomorphism of abelian groups defined by
\ban
\cZ^K_B:K(F)&\to&\cS_B\\
V&\mapsto&(2\pi)^{-\frac{\dim_\bC F}2}\cZ^{coh}_B\bigg(\widehat\Gamma_F\cup\Ch(V)\bigg).
\nan
Here $\widehat\Gamma_F\in H^*(F)$ is the \emph{Gamma class} of $F$, which is defined by $\widehat\Gamma_F:=\prodl_{i=1}^{\dim_\bC F}\Gamma(1+\delta_i)$, where $\delta_i$'s are Chern roots of the tangent bundle of $F$ and $\Gamma(x)$ is Euler's Gamma function, and $\Ch(V):=\suml_{p=0}^{\infty}(2\pi\mathbf{i})^pch_p(V)$ is the \emph{modified Chern character} of $V$. For $\mathbf t\in B$, let $\cZ^{K}_\mathbf t:K(F)\to\cS_\mathbf t$ be the restriction of $\cZ^{K}_B$ to $\{\mathbf t\}\times\widetilde{\bC^*}$. Iritani's \emph{$\widehat\Gamma$-integral structure} of $\cS_B$ (resp. $\cS_\mathbf t$) is the image of $\cZ_B$ (resp. $\cZ_\mathbf t$), which is a full rank lattice in $\cS_B$ (resp. $\cS_\mathbf t$).

\begin{remark}
	The definition of $\widehat\Gamma$-integral structure was introduced by Iritani \cite{Iri09}, and was inspired by mirror symmetry. For toric Fano manifolds, Iritani showed that the $\widehat\Gamma$-integral structures match the natural integral structures in their Landau-Ginzburg $B$-models.
\end{remark}

\if{
By abuse of notation, define bilinear non-symmetric pairings on $\cS_B$, $\cS_\mathbf t$ and $K(F)$ by:

\ban
[y_1,y_2):=\<y_1(\mathbf t,e^{-\pi\mathbf{i}}z),y_2(\mathbf t,z)\>^F,\quad\forall y_1,y_2\in\cS_B,
\nan
\ban
[y_1,y_2):=\<y_1(e^{-\pi\mathbf{i}}z),y_2(z)\>^F,\quad\forall y_1,y_2\in\cS_\mathbf t,
\nan
and
\ban
[V_1,V_2):=
\chi(V_1^\vee\otimes V_2),\quad\forall V_1,V_2\in K(F).
\nan
From \cite[formula (20)]{Iri09}, the pairings on  $\cS_B$ and $\cS_\mathbf t$ are well-defined. Moreover, from \cite[Proposition 2.10(iii)]{Iri09}, we have
\ba\label{paringsarematched}
[\cZ^K_B(V_1),\cZ^K_B(V_2))=[\cZ^K_\mathbf t(V_1),\cZ^K_\mathbf t(V_2))=[V_1,V_2),\quad\forall V_1,V_2\in K(F).
\na
}\fi

\subsection{Flat sections around \texorpdfstring{$z=0$}{z=0}}

Assume that the quantum cohomology $(H^*(F),\bullet_\mathbf t)$ is analytic and semisimple for $\mathbf t$ in a domain $B$ of $H^*(F)$ ($B$ does not necessarily contain $\mathbf{0}$). Recall that $(H^*(F),\bullet_\mathbf t)$ is semisimple if it is isomorphic to a direct sum $\underbrace{\bC\oplus\cdots\oplus\bC}_s$ as a $\bC$-algebra. Let $\psi_1(\mathbf t),\dots,\psi_s(\mathbf t)$ be the idempotent basis of $(H^*(F),\bullet_\mathbf t)$, i.e. $\psi_i(\mathbf t)\bullet_\mathbf t\psi_j(\mathbf t)=\delta_{ij}\psi_i(\mathbf t)$. Let $\Psi_i(\mathbf t):=\frac{\psi_i(\mathbf t)}{\sqrt{\<\psi_i(\mathbf t),\psi_i(\mathbf t)\>^F}}$, which is called a \emph{normalized} idempotent (defined up to sign for each $\mathbf t$). Each $\Psi_i(\mathbf t)$ is an eigenvector of $(E\bullet_\mathbf t)$, and we assume that $E\bullet_\mathbf t\Psi_i(\mathbf t)=u_i(\mathbf t)\Psi(\mathbf t)$. We say that a phase $\phi\in\bR$ is \emph{admissible} at $\mathbf t$ if $e^{\mathbf{i}\phi}$ is not parallel to any nonzero difference $u_i(\mathbf t)-u_j(\mathbf t)$.

\begin{definition}
	Let $B\subset H^*(F)$ be a domain in which $\bullet_\mathbf t$ is analytic and semisimple. We say that $B$ is \emph{properly-chosen} with respect to $\{(u_i,\Psi_i)\}_{1\leq i\leq s}$, if the maps $u_i:B\to\bC$ form analytic coordinates of $B$, and the maps $\Psi_i:B\to H^*(F)$ are single-valued and analytic. In this case, we say that $u_i$'s are canonical coordinates on $B$.
\end{definition}

From \cite[Theorem 3.1]{Dub00}, when the quantum cohomology $(H^*(F),\bullet_\mathbf t)$ is analytic and semisimple around $\mathbf t_0\in H^*(F)$, there exists a properly-chosen neighborhood of $\mathbf t_0$.

\begin{prop}\label{aefs}(\cite[Proposition 2.5.1]{GGI16})
Assume that $\bullet_\mathbf t$ is analytic and semisimple in a neighborhood $B$ of $\mathbf t_0\in H^*(F)$ properly-chosen with respect to $\{(u_i,\Psi_i)\}_{1\leq i\leq s}$, and $\phi\in\bR$ is an admissible phase at $\mathbf t_0$. Then there exists a unique basis $y_1,\dots,y_s$ of $\cS_{B}$, such that, there exists a neighborhood $B'$ of $\mathbf t_0$ in $B$ and $\varepsilon>0$ so that for each $\mathbf t\in B'$ and $i=1,\dots,s$, we have
\begin{equation}\label{eq-prop-aefs}
  y_i(\mathbf t,z)e^{\frac{u_i(\mathbf t)}z}\to\Psi_i(\mathbf t),\textrm{ as }z\to0\textrm{ in the sector }|\arg z-\phi|<\frac\pi2+\varepsilon.
\end{equation}
\end{prop}

We call the above basis of $\cS_B$ the {\em asymptotically exponential fundamental solution} (AEFS for short) to the quantum differential equation associated to the phase $\phi$ with respect to $\Psi_1,\dots,\Psi_s$ over $B'$. We will follow \cite[Section 2.5]{GGI16} to sketch the construction of the AEFS in Section \ref{sec:AEFSLT}.

\begin{definition}\label{nablaflatrespect}
	Let $B\subset H^*(F)$ be properly-chosen with respect to $\{(u_i,\Psi_i)\}_{1\leq i\leq s}$. For $y\in \cS_B$, we say that $y$ \emph{respects} $(u_i,\Psi_i)$ with phase $\phi\in\bR$ over $B$, if there exists $\varepsilon>0$ so that for each $\mathbf t\in B$, we have
	\ban
	y(\mathbf t,z)e^{\frac{u_i(\mathbf t)}z}\to\Psi_i(\mathbf t),
	\nan
	as $z\to0$ with $|\arg z-\phi|<\frac\pi2+\varepsilon$.
\end{definition}

Recall that $\nabla_{\mathbf t}$ is the restriction of $\nabla$ to $\{\mathbf t\}\times\bP^1$.

\begin{definition}\label{nablaflatrespectatonepoint}
	Let $y\in\cS_{\mathbf t}$ for a fixed $\mathbf t$. We say that $y$ \emph{respects} $(u_i(\mathbf t),\Psi_i(\mathbf t))$ with phase $\phi\in\bR$, if there exists $\varepsilon>0$ so that
	\ban
	y(z)e^{\frac{u_i(\mathbf t)}z}\to\Psi_i(\mathbf t),
	\nan
	as $z\to0$ with $|\arg z-\phi|<\frac\pi2+\varepsilon$.
\end{definition}

For any $y\in\cS_B$ and $\mathbf t\in B$, we let $y_\mathbf t\in\cS_\mathbf t$ be the restriction of $y$ to $\{\mathbf t\}\times\widetilde{\bC^*}$.

\begin{lemma}\label{respectextension}
	Under the assumption of Proposition \ref{aefs}, let $y\in\cS_B$. 
	\begin{enumerate}
		\item If $y_{\mathbf t_0}$ respects $(u_i(\mathbf t_0),\Psi_i(\mathbf t_0))$ with phase $\phi$, then $y=y_i$, where $y_i$ is in the resulted AEFS in Proposition \ref{aefs}.
		\item $y_{\mathbf t_0}$ respects $(u_i(\mathbf t_0),\Psi_i(\mathbf t_0))$ with phase $\phi$ if and only if $y$ respects $(u_i,\Psi_i)$ with phase $\phi$ over an open neighborhood of $\mathbf t_0$ (in the sense of Definition \ref{nablaflatrespect}).
	\end{enumerate}
\end{lemma}
\begin{proof}
	Let $y_1,\dots,y_s\in\cS_B$ be the AEFS in Proposition \ref{aefs}. In particular $y_1,\dots,y_s$ form a basis of $\cS_B$. Since the restriction map $\mathcal{S}_B\rightarrow \mathcal{S}_{\mathbf{t}_0}$ is an isomorphism, $y_{1,\mathbf t_0},\dots, y_{s,\mathbf t_0}$ also form a basis of $\cS_{\mathbf t_0}$. So there exist $k_1,\dots,k_s\in\bC$ such that $y=\suml_{j=1}^sk_jy_j$, and thus $y_{\mathbf t_0}=\suml_{j=1}^sk_jy_{j,\mathbf t_0}$. 
  %We will use $y_{\mathbf t_0}=\suml_{j=1}^sk_jy_{j,\mathbf t_0}$ to show $k_j=\delta_{ji}$, which implies $y=y_i$. This proves the first statement. The second statement is a consequence of the first.
  We are going to show 
\begin{equation}\label{eq-respectextension-1}
  k_j=\delta_{ij},\ \mbox{for}\ 1\leq j\leq s.
\end{equation}
If this is done, then $y=y_i$ and the first statement is proved. The second statement is a consequence of the first.

	We will prove (\ref{eq-respectextension-1}) in the following five steps. To ease the notations, we set $u_{jj'}:=\frac{u_j(\mathbf t_0)-u_{j'}(\mathbf t_0)}{e^{\mathbf{i}\phi}}$. Then by (\ref{eq-prop-aefs}) we have
	\begin{equation}\label{eq-exponential-1}
	e^{\frac{u_i(\mathbf t_0)}{z}}y_{\mathbf t_0}(z)=\suml_{j=1}^{s}k_j\exp\Big(\frac{u_{ij}}{|z|}e^{\mathbf{i}(\phi-\arg z)}\Big)\big(\Psi_j(\mathbf t_0)+g_j(z)\big),
	\end{equation}
	where $g_j(z)$ is a holomorphic $H^*(F)$-valued function, and $g_j(z)=o(1)$ as $z\to0$ in the sector $|\arg z-\phi|<\frac\pi2+\varepsilon$.
  We will use the following observation :
	\begin{itemize}
		\item Let $V$ be a normed vector space over $\bC$, and $v_1,\dots,v_k\in V$ are linearly independent. 
	   %If $c_1,\cdots,c_k\in\bC$ are not all zero, then for any $\theta_j\in\bR$, $\suml_{j=1}^kc_je^{\mathbf{i}\theta_j}v_j\neq0$. So 
     Then by the compactness of $(S^1)^k$, for any given $c_1,...,c_k\in\bC$ and $\theta_1,...,\theta_k\in \mathbb{R}$, if $c_i$ are not all zero,  we have
	    \begin{equation}
	    \min_{\theta_j\in\bR}\big|\suml_{j=1}^kc_je^{\mathbf{i}\theta_j}v_j\big|>0.\label{nonzerominimum}
	    \end{equation}
	\end{itemize}
	
	\begin{itemize}
	\item[Step 1:] Let $M_1:=\max\limits_{\{j|k_j\neq0\}}\Im u_{ij}$. We  show that if $\Im u_{ij}>0$, then $k_j=0$. Suppose this is not true.  Then $M_1>0$. Let $z\to0$ with $\arg z=\phi+\frac\pi2$, then (\ref{eq-exponential-1}) reads
		\ban
		|e^{\frac{u_i(\mathbf t_0)}{z}}y_{\mathbf t_0}(z)|&=&\Big|\suml_{j}k_j\exp\Big(\frac{\Im u_{ij}-\mathbf{i}\Re u_{ij}}{|z|}\Big)(\Psi_j(\mathbf t_0)+g_j(z))\Big|\\
		&\geq&\exp\Big(\frac{M_1}{\abs{z}}\Big)\Big\{\Big|\suml_{\{j|\Im u_{ij}=M_1\}}k_j\exp\Big(\frac{-\mathbf{i}\Re u_{ij}}{\abs{z}}\Big)(\Psi_j(\mathbf t_0)+g_j(z))\Big|\\
		&&\qquad-\suml_{\{j|\Im u_{ij}<M_1\}}\abs{k_j}\exp\Big(\frac{\Im u_{ij}-M_1}{\abs{z}}\Big)\abs{\Psi_j(\mathbf t_0)+g_j(z)}\Big\}.
		\nan
		Note that
		\[
		\suml_{\{j|\Im u_{ij}<M_1\}}\abs{k_j}\exp\Big(\frac{\Im u_{ij}-M_1}{\abs{z}}\Big)\abs{\Psi_j(\mathbf t_0)+g_j(z)}\to0,
		\]
		and 
		\[
		\Big|\suml_{\{j|\Im u_{ij}=M_1\}}k_j\exp\Big(\frac{-\mathbf{i}\Re u_{ij}}{\abs{z}}\Big)\big(\Psi_j(\mathbf t_0)+g_j(z)\big)\Big|=\Big|\suml_{\{j|\Im u_{ij}=M_1\}}k_j\exp\Big(\frac{-\mathbf{i}\Re u_{ij}}{\abs{z}}\Big)\Psi_j(\mathbf t_0)\Big|+o(1).
		\]
		Moreover, note that $\Psi_j(\mathbf t_0)$'s are linearly independent, and by our assumption, there exists $j$ such that $\Im u_{ij}=M_1$ and $k_j\neq0$. As a consequence, from \eqref{nonzerominimum}, we have
		\[
		C_1:=\min\Big\{\Big|\suml_{\{j|\Im u_{ij}=M_1\}}k_j\exp(\mathbf{i}\theta_j)\Psi_j(\mathbf t_0)\Big|:\theta_j\in\bR\Big\}>0,
		\]
		and then
		\[
		\liminf\Big|\suml_{\{j|\Im u_{ij}=M_1\}}k_j\exp\Big(\frac{-\mathbf{i}\Re u_{ij}}{\abs{z}}\Big)(\Psi_j(\mathbf t_0)+g_j(z))\Big|\geq C_1>0.
		\]
		This implies that $|e^{\frac{u_i(\mathbf t_0)}{z}}y_{\mathbf t_0}(z)|\to\infty$ as $z\to0$ with $\arg z=\phi+\frac\pi2$, contradicting the assumption on $y_{\mathbf t_0}(z)$.
	\item[Step 2:] Let  $M_2:=
    \max\limits_{\{j|k_j\neq0\}}(-\Im u_{ij})$. We show that if $\Im u_{ij}<0$, then $k_j=0$. Suppose this is not true.  Then $M_2>0$. Now let $z\to0$ with $\arg z=\phi-\frac\pi2$, and we have
	    \ban
		|e^{\frac{u_i(\mathbf t_0)}{z}}y_{\mathbf t_0}(z)|&=&\Big|\suml_{j}k_j\exp\Big(\frac{-\Im u_{ij}+\mathbf{i}\Re u_{ij}}{|z|}\Big)\big(\Psi_j(\mathbf t_0)+g_j(z)\big)\Big|\\
		&\geq&\exp\Big(\frac{M_2}{|z|}\Big)\Big\{\big|\suml_{\{j|-\Im u_{ij}=M_2\}}k_j\exp\Big(\frac{\mathbf{i}\Re u_{ij}}{|z|}\Big)\big(\Psi_j(\mathbf t_0)+g_j(z)\big)\big|\\
		&&\quad-\suml_{\{j|-\Im u_{ij}<M_2\}}\abs{k_j}\exp\Big(\frac{-\Im u_{ij}-M_2}{\abs{z}}\Big)\abs{\Psi_j(\mathbf t_0)+g_j(z)}\Big\}.
		\nan
		Note that
		\[
		\suml_{\{j|-\Im u_{ij}<M_2\}}\abs{k_j}\exp\Big(\frac{-\Im u_{ij}-M_2}{\abs{z}}\Big)\abs{\Psi_j(\mathbf t_0)+g_j(z)}\to0,
		\]
		and 
		\[
		\Big|\suml_{\{j|-\Im u_{ij}=M_2\}}k_j\exp\Big(\frac{\mathbf{i}\Re u_{ij}}{\abs{z}}\Big)\big(\Psi_j(\mathbf t_0)+g_j(z)\big)\Big|=\Big|\suml_{\{j|-\Im u_{ij}=M_2\}}k_j\exp\Big(\frac{\mathbf{i}\Re u_{ij}}{\abs{z}}\Big)\Psi_j(\mathbf t_0)\Big|+o(1).
		\]
		Moreover, note that $\Psi_j(\mathbf t_0)$'s are linearly independent, and by our assumption, there exists $j$ such that $-\Im u_{ij}=M_2$ and $k_j\neq0$. As a consequence, from \eqref{nonzerominimum}, we have
		\[
		C_2:=\min\Big\{\Big|\suml_{\{j|-\Im u_{ij}=M_2\}}k_j\exp(\mathbf{i}\theta_j)\Psi_j(\mathbf t_0)\Big|:\theta_j\in\bR\Big\}>0,
		\]
		and then
		\[
		\liminf\Big|\suml_{\{j|-\Im u_{ij}=M_2\}}k_j\exp\Big(\frac{\mathbf{i}\Re u_{ij}}{\abs{z}}\Big)(\Psi_j(\mathbf t_0)+g_j(z))\Big|\geq C_2>0.
		\]
	    This implies that $|e^{\frac{u_i(\mathbf t_0)}{z}}y_{\mathbf t_0}(z)|\to\infty$ as $z\to0$ with $\arg z=\phi-\frac\pi2$, contradicting the assumption on $y_{\mathbf t_0}(z)$.
	\item[Step 3:] Let $M_3:=
    \max\limits_{\{j|k_j\neq0\}}\Re u_{ij}$. From Steps 1 and 2, we know that if $\Im u_{ij}>0$, then $k_j=0$.
   Now we show that if $\Re u_{ij}>0$, then $k_j=0$. Suppose this is not true. Then  $M_3>0$. Let $z\to0$ with $\arg z=\phi$, and we have
		\ban
		|e^{\frac{u_i(\mathbf t_0)}{z}}y_{\mathbf t_0}(z)|&=&\Big|\suml_{j}k_j\exp\Big(\frac{\Re u_{ij}}{|z|}\Big)\big(\Psi_j(\mathbf t_0)+g_j(z)\big)\Big|\\
		&\geq&\exp\Big(\frac{M_3}{|z|}\Big)\Big\{\big|\suml_{\{j|\Re u_{ij}=M_3\}}k_j\big(\Psi_j(\mathbf t_0)+g_j(z)\big)\big|\\
		&&\qquad-\suml_{\{j|\Re u_{ij}<M_3\}}\abs{k_j}\exp\Big(\frac{\Re u_{ij}-M_3}{|z|}\Big)\abs{\Psi_j(\mathbf t_0)+g_j(z)}\Big\}.
		\nan
		Note that
		\[
		\suml_{\{j|\Re u_{ij}<M_3\}}\abs{k_j}\exp\Big(\frac{\Re u_{ij}-M_3}{|z|}\Big)\abs{\Psi_j(\mathbf t_0)+g_j(z)}\to0,
		\]
		and 
		\[
		\lim\Big|\suml_{\{j|\Re u_{ij}=M_3\}}k_j\big(\Psi_j(\mathbf t_0)+g_j(z)\big)\Big|
    =\Big|\suml_{\{j|\Re u_{ij}=M_3\}}k_j\Psi_j(\mathbf t_0)\Big|.
		\]
		Moreover, note that $\Psi_j(\mathbf t_0)$'s are linearly independent, and by our assumption, there exists $j$ such that $\Re u_{ij}=M_3$ and $k_j\neq0$. As a consequence, we have
		\[
		\Big|\suml_{\{j|\Re u_{ij}=M_3\}}k_j\Psi_j(\mathbf t_0)\Big|>0.
		\]
	 This implies that $|e^{\frac{u_i(\mathbf t_0)}{z}}y_{\mathbf t_0}(z)|\to\infty$ as $z\to0$ with $\arg z=\phi$, contradicting the assumption on $y_{\mathbf t_0}(z)$.
	\item[Step 4:] Let $M_4:=
    \max\limits_{\{j|k_j\neq0\}}(-\Re u_{ij})$.  We show that if $\Re u_{ij}<0$, then $k_j=0$. Suppose this is not true. We have $M_4>0$. Now fix $\varepsilon'\in(0,\varepsilon)\cap(0,\frac\pi2)$. Then let $z\to0$ with $\arg z=\phi-(\frac\pi2+\varepsilon')$, and we have
		\ban
		\big|e^{\frac{u_i(\mathbf t_0)}{z}}y_{\mathbf t_0}(z)\big|&=&\Big|\suml_{j}k_j\exp\Big(\frac{\Re u_{ij}}{|z|}(-\sin\varepsilon'+\mathbf{i}\cos\varepsilon')\Big)\big(\Psi_j(\mathbf t_0)+g_j(z)\big)\Big|\\
		&\geq&\exp\Big(\frac{M_4}{|z|}\sin\varepsilon'\Big)\Big\{\big|\suml_{\{j|-\Re u_{ij}=M_4\}}k_j\exp\Big(\mathbf{i}\frac{\Re u_{ij}}{|z|}\cos\varepsilon'\Big)\big(\Psi_j(\mathbf t_0)+g_j(z)\big)\big|\\
		&&\qquad-\suml_{\{j|-\Re u_{ij}<M_4\}}\abs{k_j}\exp\Big(\frac{-\Re u_{ij}-M_4}{|z|}\sin\varepsilon'\Big)\abs{\Psi_j(\mathbf t_0)+g_j(z)}\Big\}.
		\nan
		Note that
		\[
		\suml_{\{j|-\Re u_{ij}<M_4\}}\abs{k_j}\exp\Big(\frac{-\Re u_{ij}-M_4}{|z|}\sin\varepsilon'\Big)\abs{\Psi_j(\mathbf t_0)+g_j(z)}\to0,\ (\mbox{where we use that}\ 0<\varepsilon'<\frac{\pi}{2})
		\]
		and 
		\begin{eqnarray*}
		&&\Big|\suml_{\{j|-\Re u_{ij}=M_4\}}k_j\exp\Big(\mathbf{i}\frac{\Re u_{ij}}{|z|}\cos\varepsilon'\Big)\big(\Psi_j(\mathbf t_0)+g_j(z)\big)\Big|\\
    &=&\Big|\suml_{\{j|-\Re u_{ij}=M_4\}}k_j\exp\Big(\mathbf{i}\frac{\Re u_{ij}}{|z|}\cos\varepsilon'\Big)\Psi_j(\mathbf t_0)\Big|+o(1).
		\end{eqnarray*}
		Moreover, note that $\Psi_j(\mathbf t_0)$'s are linearly independent, and by our assumption, there exists $j$ such that $-\Re u_{ij}=M_4$ and $k_j\neq0$. As a consequence, from \eqref{nonzerominimum}, we have
		\[
		C_4:=\min\Big\{\Big|\suml_{\{j|-\Re u_{ij}=M_4\}}k_j\exp(\mathbf{i}\theta_j)\Psi_j(\mathbf t_0)\Big|:\theta_j\in\bR\Big\}>0\ (\mbox{we use again that}\ 0<\varepsilon'<\frac{\pi}{2}),
		\]
		and then
		\[
		\liminf\Big|\suml_{\{j|-\Re u_{ij}=M_4\}}k_j\exp\Big(\mathbf{i}\frac{\Re u_{ij}}{|z|}\cos\varepsilon'\Big)\big(\Psi_j(\mathbf t_0)+g_j(z)\big)\Big|\geq C_4>0.
		\]
	 This implies that $|e^{\frac{u_i(\mathbf t_0)}{z}}y_{\mathbf t_0}(z)|\to\infty$ as $z\to0$ with $\arg z=\phi-(\frac\pi2+\varepsilon')$, contradicting the assumption on $y_{\mathbf t_0}(z)$.
	\item[Step 5:] We come to the conclusion that 
		\ban
		e^{\frac{u_i(\mathbf t_0)}{z}}y_{\mathbf t_0}(z)=\suml_{\{j|u_{ij}=0\}}k_j\big(\Psi_j(\mathbf t_0)+g_j(z)\big).
		\nan
		Since 
		\ban
		\Psi_i(\mathbf t_0)=\lim\limits_{\substack{z\to0\\\arg z=\phi}}e^{\frac{u_i(\mathbf t_0)}{z}}y_{\mathbf t_0}(z)=\suml_{\{j|u_{ij}=0\}}k_j\Psi_j(\mathbf t_0),
		\nan
		it follows from the linear independence of $\Psi_j(\mathbf t_0)$'s that $k_j=\delta_{ij}$.
	\end{itemize}
\end{proof}

\subsection{Statement of Gamma conjecture II}
Roughly speaking, Gamma conjecture II expects that an AEFS is in the $\widehat\Gamma$-integral structure. The precise statement (\cite[Conjecture 4.6.1]{GGI16}, \cite[Conjecture 4.9]{GaIr19}) is as follows.

\begin{gammaII}
Assume that: (i) $\bullet_\mathbf t$ is analytic and semisimple in a neighborhood $B$ of $\mathbf t_0\in H^*(F)$ properly-chosen with respect to $\{(u_i,\Psi_i)\}_{1\leq i\leq s}$; (ii) the bounded derived category of coherent sheaves $\mathcal D^b(F)$ admits a full exceptional collection.
Let $\phi$ be an admissible phase at $\mathbf t_0$, and we numbers $u_i$'s such that
\ban
\Im(e^{-\mathbf{i}\phi}u_1(\mathbf t_0))\geq\Im(e^{-\mathbf{i}\phi}u_2(\mathbf t_0))\geq\dots\geq\Im(e^{-\mathbf{i}\phi}u_s(\mathbf t_0)).
\nan
Then there exists a full exceptional collection $(E_1,\dots,E_s)$ of $\mathcal D^b(F)$ such that $\cZ^K_{\mathbf t_0}(E_i)$ respects $(u_i(\mathbf t_0),\Psi_i(\mathbf t_0))$ for $1\leq i\leq s$.
\end{gammaII}

\begin{remark}
	\begin{itemize}
		\item[(1)] By the definition of admissible phase, $\Im(e^{-\mathbf{i}\phi}u_i(\mathbf t_0))=\Im(e^{-\mathbf{i}\phi}u_j(\mathbf t_0))$ if and only if $u_i(\mathbf t_0)=u_j(\mathbf t_0)$.
		\item[(2)] From Lemma \ref{respectextension}, the basis $(\cZ^K_B(E_1),\dots,\cZ^K_B(E_s))$ of $\cS_B$ is actually the AEFS associated to the phase $\phi$ with respect to $\Psi_1,\dots,\Psi_s$ around $\mathbf t_0$.
    \item[(3)] A conjecture of Dubrovin \cite{Dub98}, made more precise in \cite{Bay04,HMT09}, states that the semisimplicity of the quantum cohomology of $F$ is equivalent to the existence of full exceptional collections in $\mathcal{D}^b(F)$.
	\end{itemize}
\end{remark}

As pointed out in \cite[Remark 4.6.3]{GGI16} and \cite[Remark 4.13]{GaIr19}, the validity of Gamma conjecture II does not depend on the choice of $(\mathbf t_0,\phi)\in H^*(F)\times\bR$, as long as the quantum cohomology is analytic and semisimple around $\mathbf t_0$, and $\phi$ is admissible at $\mathbf t_0$. This is because as $(\mathbf t_0,\phi)$ varies, the AEFS changes by mutations, and we can consider the corresponding mutations on full exceptional collections. We refer readers to \cite[Section 4]{GGI16} for detailed discussions.

\section{A sufficient condition for Gamma conjecture II}\label{sec:sufficientCondition-GammaII}
This section is devoted to proving Theorem \ref{sufficientcondition}, a criterion for Gamma conjecture II, which will be used in \S 6. For $u\in \mathbb{C}$ and $\phi\in \mathbb{R}$, define
\begin{equation}
  L(u,\phi):=u+\bR_{\geq0}e^{\mathbf{i}\phi},
\end{equation}
i.e. $L(u,\phi)$ is the oriented half line in $\mathbb{C}$ from $u$ to $\infty$ with phase $\phi$. To state Theorem \ref{sufficientcondition}, we need some preparations.

\subsection{AEFS via Laplace transformation}\label{sec:AEFSLT}
In this subsection, we follow \cite[Section 2.5]{GGI16} to sketch the construction of the AEFS in Proposition \ref{aefs}. Recall that $B$ is properly-chosen with respect to $\{(u_i,\Psi_i)\}_{1\leq i\leq s}$.

Let $\widehat\nabla$ be the \emph{Laplace-dual connection} of $\nabla$ (see \cite[formula (2.5.2)]{GGI16} for the precise definition of $\widehat\nabla$), which is a meromorphic flat connection on the trivial $H^*(F)$-bundle over $B\times\bC_\lambda$, where $\bC_\lambda$ is a copy of $\bC$ with coordinate $\lambda$. 

Let $D_i$ be the smooth divisor in $B\times\bC_\lambda$ defined by $\lambda=u_i$, and let $D_\infty:=B\times\{\infty\}$ be the divisor at infinity. Then $\widehat\nabla$ has only logarithmic singularities at $D_1, \dots, D_s, D_\infty$. 

\begin{definition}
	Let $\hat y$ be a $\widehat\nabla$-flat section near $D_i$. We say that $\hat y$ \emph{respects} $(u_i,\Psi_i)$ over $B$, if we can analytically continue $\hat y$ so that $\hat y(\mathbf t, u_i(\mathbf t))=\Psi_i(\mathbf t)$ for all $\mathbf t\in B$.
\end{definition}

Though $\widehat\nabla$ is singular along the normal crossing divisor $=D_1\cup\cdots\cup D_s$ in $B\times\bC_\lambda$, it was shown in \cite[Lemma 2.5.3]{GGI16} that there exists a $\widehat\nabla$-flat section $\hat y_i$ respecting $(u_i,\Psi_i)$ around $\mathbf t_0$. Now let $\phi$ be an admissible phase at $\mathbf t_0$, and we define the \emph{Laplace transform of $\hat y_i$ associated to the phase $\phi$ with respect to $u_i$} by
\ba\label{laplacetransform}
\bar y_i(\mathbf t,z)=\frac1z\int_{\lambda\in L(u_i(\mathbf t),\phi)}\hat y_i(\mathbf t,\lambda)e^{-\frac\lambda z}d\lambda,\quad\mathbf t\textrm{ around }\mathbf t_0, |\arg z-\phi|<\frac\pi2.
\na
As in \cite[Proof of Proposition 2.5.1]{GGI16}, by slightly varying the slope of the integration contour in \eqref{laplacetransform}, we can analytically continue $\bar y_i$ to $|\arg z-\phi|<\frac\pi2+\varepsilon$ for some $\varepsilon>0$. Now one can show that $\bar y_1,\dots,\bar y_s$ form the required AEFS.

\subsection{Analytic continuation of Dubrovin connection}
\label{sec:continuationofDubrovinconnection}

Assume that $B$ is a  properly-chosen neighborhood of $\mathbf t_0\in H^*(F)$ with respect to $\{(u_i,\Psi_i)\}_{1\leq i\leq s}$, such that for each $\mathbf t\in B$, we have $u_i(\mathbf t)\neq u_j(\mathbf t)$ for $i\neq j$. Recall that both $t_i$'s and $u_i$'s are coordinates on $B$, and we can view $t_i$'s as functions of the chosen canonical coordinates.

Let
\ban
\Delta:=\{(u_1,\dots,u_s)\in\bC^s:u_i= u_j\textrm{ for some }i\neq j\},
\nan
and write $\mathbf{w}_0:=(u_1(\mathbf t_0),\dots,u_s(\mathbf t_0))\in\bC^s\setminus\Delta$. Let $\mathcal C$ be the universal covering of $\bC^s\setminus\Delta$ constructed as homotopy classes of paths starting from $\mathbf{w}_0$. By abuse of notation we still denote by $\mathbf{w}_0$  the point in  $\mathcal C$ corresponding to the trivial loop at $\mathbf{w}_0$, and  use the chosen canonical coordinates $u_i$'s to identify the open neighborhood $B$ of $\mathbf t_0$ with a neighborhood of $\mathbf{w}_0$ in $\mathcal C$. Then via the diagram
\[
\xymatrix@R-1.5pc{
  B  \ar@{^{(}->}[r] & \mathbb{C}^s\setminus \Delta & \mathcal{C} \ar[l] \\
  t_0 \ar@{|->}[r] \ar@{}[u]|*[@]{\in} & \mathbf{w}_0 \ar@{}[u]|*[@]{\in} & \mathbf{w}_0 \ar@{|->}[l] \ar@{}[u]|*[@]{\in}
}
\]
we regard $u_i$'s as functions of $\mathbf{w}\in \mathcal{C}$.

From \cite[Theorem 4.5, Theorem 4.6]{Dub00}, $\nabla_{\mathbf t_0}$ admits a unique isomonodromic deformation over $\mathcal C$. More precisely, there is a unique meromorphic flat connection $\nabla$ on the trivial $H^*(F)$-bundle over $\mathcal C\times\bP^1$ of the form
\ban
\nabla_{\partial_{u_i}}&=&\partial_{u_i}+\frac1z\mathcal U_i,\\
\nabla_{z\partial_z}&=&z\partial_z-\frac1z\mathcal U+\mu,
\nan
where $\mathcal U_i$'s and $\mathcal U$ are $\End(H^*(F))$-valued meromorphic functions on $\mathcal C$, such that $\nabla$ restricts to the Dubrovin connection on $B\times\bP^1$ (via the above mentioned identification $B\hookrightarrow\mathcal C$). This isomonodromic deformation defines a semisimple Frobenius manifold structure on a dense open subset $\mathcal B$ of $\mathcal C$, such that the complement $\mathcal C\setminus\mathcal B$ is a divisor, and we have $B\hookrightarrow\mathcal B$ via the above mentioned identification. In particular, $\mathcal B$ is connected.

Now we can continue the functions $u_i$'s and $\Phi_i$'s on $B$ to analytic functions on $\mathcal B$. Moreover, as functions of chosen canonical coordinates on $B$, from \cite[(4.54)]{Dub00}, the functions $t_i$'s can be analytically continued to $\mathcal B$. So, based on $\mathcal Z^{coh}_B$ and $\mathcal Z^K_B$, we can use this analytic continuation to define cohomology framing and $K$-group framing on $\mathcal B$. Moreover, the Laplace-dual connection $\widehat\nabla$ can also be meromorphically continued to $\mathcal B\times\bC_\lambda$, and we can use Laplace transformation of $\widehat\nabla$-flat sections to study asymptotic behavior of $\nabla$-flat sections as in Section \ref{sec:AEFSLT}, which gives AEFS over an open neighborhood of a point in $\mathcal B$.

So we can consider Gamma conjecture II for $F$ over $\mathcal B$, that is, for $(\mathbf{w},\phi)\in\mathcal B\times\bR$ with $\phi$ admissible at $\mathbf{w}$, matching the corresponding AEFS with flat sections from a full exceptional collection via $\mathcal Z^K_\mathbf{w}$. Now it follows from \cite[Section 4]{GGI16} that the validity of Gamma conjecture II does not depend on the choice of $(\mathbf{w},\phi)$.

\subsection{A criterion for Gamma conjecture II}

In this subsection, we use notations from Section \ref{sec:continuationofDubrovinconnection}. The main result of this subsection is Theorem \ref{sufficientcondition}.

For $\mathbf{w}\in\mathcal B$, we regard $u_1(\mathbf{w}), \dots, u_s(\mathbf{w})$ as pairwise distinct points in $\bC_\lambda$. Let
\ban
L(u_i(\mathbf{w}),\phi)^\circ:=L(u_i(\mathbf{w}),\phi)\setminus\{u_i(\mathbf{w})\}
\nan
be an open half line in $\bC_\lambda$, and we define
\ban
A_i(\mathbf{w}):=\{\phi\in\bR:L(u_i(\mathbf{w}),\phi)^\circ\textrm{ does not contain any of }u_1(\mathbf{w}),\dots,u_s(\mathbf{w})\}.
\nan
Then $A_i(\mathbf{w})$ is an open subset in $\bR$ such that $\bR\setminus A_i(\mathbf{w})$ is discrete, and the set of admissible phases at $\mathbf{w}$ is
\ban
A_1(\mathbf{w})\cap\cdots\cap A_s(\mathbf{w}).
\nan

Let $\widehat\nabla_{\mathbf{w}}$ be the restriction of $\widehat\nabla$ on $\{\mathbf{w}\}\times\bC_\lambda$. 

\begin{definition}\label{LT}
	Let $\hat y$ be a $\widehat\nabla_{\mathbf{w}}$-flat section for a fixed $\mathbf{w}$. If $\hat y$ is holomorphic near $\lambda=u_i(\mathbf{w})$, and $\hat y(u_i(\mathbf{w}))=\Psi_i(\mathbf{w})$, then we say that $\hat y$ \emph{respects} $(u_i(\mathbf{w}),\Psi_i(\mathbf{w}))$, and we define the \emph{Laplace transform of $\hat y$ associated to phase $\phi\in A_i(\mathbf{w})$ with respect to $u_i(\mathbf{w})$} by
	\ban
	y(z)=\frac1z\int_{\lambda\in L(u_i(\mathbf{w}),\phi)}\hat y(\lambda)e^{-\frac\lambda z}d\lambda,\quad|\arg z-\phi|<\frac\pi2.
	\nan
\end{definition}

One can use arguments in \cite[Proof of Proposition 2.5.1]{GGI16} to prove that if a $\widehat\nabla_{\mathbf{w}}$-flat section $\hat y$ respects $(u_i(\mathbf{w}),\Psi_i(\mathbf{w}))$, then its Laplace transform $y$ is actually $\nabla_\mathbf{w}$-flat and respects $(u_i(\mathbf{w}),\Psi_i(\mathbf{w}))$ with phase $\phi$ (in the sense of Definition \ref{nablaflatrespectatonepoint}). 

The following Lemma \ref{varycontuor} is well-known to experts, and we state it here for convenience of readers.
\begin{lemma}\label{varycontuor}
	Suppose that two phases $\phi$ and $\phi'$ are in the same connected component of $A_i(\mathbf{w})$, and $\hat y$ is a $\widehat\nabla_\mathbf{w}$-flat section respecting $(u_i(\mathbf{w}),\Psi_i(\mathbf{w}))$. Let $y^\phi,y^{\phi'}$ be Laplace transforms of $y$ with respect to $u_i(\mathbf{w})$ associated the the phases $\phi,\phi'$ respectively. Then $y^\phi=y^{\phi'}$.
\end{lemma}
\begin{proof}
	Note that $\widehat\nabla_\mathbf{w}$ is regular singular at $\lambda=\infty$, which implies that $\hat y(\lambda)$ grows at most polynomially as $\lambda\to\infty$, i.e. $|\hat y(\lambda)|$ is bounded by a polynomial of $|\lambda|$ as $\lambda\to\infty$. One can use this observation to compare the two contour integration defining $y^\phi$ and $y^{\phi'}$. We leave the details to interested readers.
\end{proof}

\begin{lemma}\label{existenceofintrep}
	Let $y$ be the $\nabla_{\mathbf{w}}$-flat section respecting $(u_i(\mathbf{w}),\Psi_i(\mathbf{w}))$ with phase $\phi\in A_i(\mathbf{w})$. Then there exists a $\widehat\nabla_{\mathbf{w}}$-flat section $\hat y$ respecting $(u_i(\mathbf{w}),\Psi_i(\mathbf{w}))$, such that $y$ is the Laplace transform of $\hat y$ associated to the phase $\phi$ with respect to $u_i(\mathbf{w})$.
\end{lemma}
\begin{proof}
	From Definition \ref{nablaflatrespect}, we can choose $\phi'\in A(\mathbf{w})$ near $\phi$, such that the two phases $\phi$ and $\phi'$ are in the same connected component of $A_i(\mathbf{w})$, and that $y$ also respects $(u_i(\mathbf{w}),\Psi_i(\mathbf{w}))$ with phase $\phi'$. Let $y_1,\dots,y_s$ be the AEFS associated to the phase $\phi'$ with respect to $\Psi_1,\dots,\Psi_s$. Then from Lemma \ref{respectextension}, we have $y=y_{i,\mathbf t_0}$. Note that we can use Laplace transforms of $\widehat\nabla$-flat sections to construct an AEFS. So from the uniqueness of AEFS, we see that $y$ is the Laplace transform of a $\widehat\nabla_\mathbf{w}$-flat section $\hat y$ associated to the phase $\phi'$ with respect to $u_i(\mathbf{w})$. Now the conclusion follows from Lemma \ref{varycontuor}.
\end{proof}

\begin{lemma}\label{varycontuor2}
	Suppose that two phases $\phi$ and $\phi'$ are in the same connected component of $A_i(\mathbf{w})$, and $y\in\cS_\mathbf{w}$ respects $(u_i(\mathbf{w}),\Psi_i(\mathbf{w}))$ with phase $\phi$. Then $y$ also respects $(u_i(\mathbf{w}),\Psi_i(\mathbf{w}))$ with phase $\phi'$.
\end{lemma}
\begin{proof}
    From Lemma \ref{existenceofintrep}, $y$ is the Laplace transform of a $\widehat\nabla_\mathbf{w}$-flat section $\hat y$ associated to the phase $\phi$ with respect to $u_i(\mathbf{w})$. Now the required result follows from Lemma \ref{varycontuor}.
\end{proof}

\begin{prop}\label{moveu}
	Let $y\in\mathcal S_{\mathcal B}$ be such that $y_{\mathbf{w}_0}$ respects $(u_i(\mathbf{w}_0),\Phi_i(\mathbf{w}_0))=(u_i(\mathbf t_0),\Phi_i(\mathbf t_0))$ with phase $\phi_0\in A_i(\mathbf{w}_0)$. Given a path $\{\mathbf{w}_t\}_{0\leq t\leq 1}$ in $\mathcal B$, assume that there is a continuous map $\phi:[0,1]\to\bR$ starting from $\phi_0$ such that $\phi(t)\in A_i(\mathbf{w}_t)$ for each $t\in[0,1]$. Then $y_{\mathbf{w}_1}$ respects $(u_i(\mathbf{w}_1),\Phi_i(\mathbf{w}_1))$ with phase $\phi_1$.
\end{prop}
\begin{proof}
	From \cite[Lemma 2.5.3]{GGI16}, for each $t\in[0,1]$, there exists a $\widehat\nabla$-flat section $\hat y_t$ respecting $(u_i,\Psi_i)$ around $\mathbf{w}_t$, and we let $\bar y_t\in\cS_\mathcal B$ be the Laplace transform of $\hat y_t$ associated to the phase $\phi(t)$ with respect to $u_i$. Then from Lemma \ref{respectextension}, for $t'\in[0,1]$ around $t$, $\bar y_{t,\mathbf{w}_{t'}}\in\cS_{\mathbf{w}_{t'}}$ respects $(u_i(\mathbf{w}_{t'}),\Psi_i(\mathbf{w}_{t'}))$ with phase $\phi(t)$. Let $(a_t,b_t)$ be the connected component of $A_i(\mathbf{w}_t)$ containing $\phi(t)$. Then there exists a connected open neighborhood $N_t$ of $t$ in $[0,1]$ such that $\phi(t')\in(a_t,b_t)$ for $t'\in N_t$. So from Lemma \ref{varycontuor2}, $\bar y_{t,\mathbf{w}_{t'}}\in\cS_{\mathbf{w}_{t'}}$ also respects $(u_i(\mathbf{w}_{t'}),\Psi_i(\mathbf{w}_{t'}))$ with phase $\phi(t')$. From Lemma \ref{respectextension}, we have $\bar y_t=\bar y_{t'}$ for $t'\in N_t$. Since $[0,1]$ is compact and connected, it follows that $y=\bar y_t$ for all $t\in[0,1]$. This finishes the proof of the proposition.
\end{proof}

\begin{cor}\label{moveuu}
	For $1\leq i\leq s$, let $y_i\in\cS_\mathcal B$ be such that $y_{i,\mathbf{w}_0}\in\cS_{\mathbf{w}_0}$ respects $(u_i(\mathbf{w}_0),\Phi_i(\mathbf{w}_0))$ with phase $\phi_{i,0}\in A_i(\mathbf{w}_0)$. Given a path $\{\mathbf{w}_t\}_{0\leq t\leq 1}$ in $\mathcal B$, assume that there are continuous maps $\phi_i:[0,1]\to\bR$ starting from $\phi_{i,0}$ such that $\phi_i(t)\in A_i(\mathbf{w}_t)$ for each $t\in[0,1]$. Then for each $i$, $y_{i,\mathbf{w}_1}$ respects $(u_i(\mathbf{w}_1),\Phi_i(\mathbf{w}_1))$ with phase $\phi_i(1)$.
\end{cor}
\begin{proof}
	This is a consequence of Proposition \ref{moveu}.
\end{proof}

\begin{thm}\label{sufficientcondition}
For a Fano manifold $F$,
	let $\mathcal E=(E_1,\dots,E_s)$ be a full exceptional collection of $\mathcal D^b(F)$, and set $y_i=\mathcal Z^K_\mathcal B(E_i)\in\cS_\mathcal B$, where $\mathcal{B}$ and $\mathcal{Z}_{\mathcal{B}}^K(E_i)$ are defined as in \S 3.2. Assume that the followings hold:
	\begin{enumerate}
		\item $y_i$ respects $(u_i,\Psi_i)$ with phase $\phi_i$ in an open neighborhood of  $\mathbf{w}_0$, with
		\ban
		\phi_1>\dots>\phi_s,\quad\phi_1-\phi_s<2\pi.
		\nan
		\item The half lines $L(u_i(\mathbf{w}_0),\phi_i)$ are pairwise disjoint.
	\end{enumerate}
	Then  Gamma conjecture II holds for $F$. 
\end{thm} 
\begin{proof} Without loss of generality, we assume that
	\ban
	2\pi>\phi_1>\dots>\phi_s>0,
	\nan
	and for any nonzero complex number $\lambda$, we choose its principal argument $\Arg(\lambda)\in[0,2\pi)$. Here we give a pictorial proof, with the illustrations in the next page. A more formal and more detailed proof will be given after the pictorial proof.

The basic idea is to move all the half lines $L(u_i(\mathbf{w}_0),\phi_i)$ to $L(u_i(\mathbf{w}^*),\phi^*_i)$ in $\mathbb{C}$, such that $\phi^*_1=\dots=\phi^*_{s}$, and $u_i(\mathbf{w}^*)$ are located in the fourth quadrant, and such that
\begin{subequations}\label{eq-requirementforpictorialproof}
\begin{eqnarray}
&2\pi>\mathrm{Arg}\big(u_1(\mathbf{w}^*)\big)>\dots>\mathrm{Arg}\big(u_s(\mathbf{w}^*)\big)>\frac{3\pi}{2},&\\
&0<|u_1(\mathbf{w}^*)|<\dots<|u_s(\mathbf{w}^*)|,&
\end{eqnarray}
\end{subequations}
and thus Gamma conjecture II holds. But we require that the movement of the half lines satisfies:
\begin{enumerate}
  	\item[(i)] the configuration of the starting points of the half lines lies in $\mathcal{C}\backslash \mathcal{B}$;
  	\item[(ii)] the half lines are disjoint in $\mathbb{C}$.
  \end{enumerate}
The way to find a movement satisfying both conditions (i) and (ii) are illustrated in the following, where $s$ equals $5$ as an example. In the first step we take an open neighborhood $V_i$ of each half line $L_i=L(u_i(\mathbf{w}_0),\phi_i)$, such that $V_i$ are pairwise disjoint for $1\leq i\leq s$, and move the half line $L_i$ towards $\infty$ along the direction $\phi_i$, until its starting point lay on a circle centered at the origin. Recall that the configuration of starting points of these half lines is a point in $\mathcal C$. In the process if the configuration meets the divisor $\mathcal{B}$, we slightly move these half lines in $V_i$ such that the resulted configuration avoids $\mathcal{B}$; this is possible because $\mathcal{B}$ has real codimension 2. After this step we obtain $L_i'=L(u_i(\mathbf{w}'),\phi'_i)$ for $1\leq i\leq s$.

In the second step, we take an open neighborhood $U$ of the circle, and rotate successively
$L'_1,\dots,L'_s$ counterclockwise, such that all 
the starting points of $L'_i$
lie in the fourth quadrant part of the thickened circle, and the direction of the half lines $L'_i$
is the direction of the radius from the center to its starting point.
In the process, for the same reason as in the first step, if the configuration of starting points of these half lines meets the divisor $\mathcal{B}$, we slightly move these half lines with their starting points remained in $U$, such that the resulted configuration avoids $\mathcal{B}$. After this step we obtain $L''_i=L(u_i(\mathbf{w}''),\phi''_i)$ for $1\leq i\leq s$.

Write $u_j=x_j+\mathbf{i}y_j$ for $1\leq j \leq s$, where $x_j,y_j\in \mathbb{R}$.
In the third step, we move each $L''_i$ along its direction, such that all the starting points of the resulted half lines, for $1\leq i\leq s$,
have the same $x$-coordinates.
After this process we also manage to make the $y$-coordinates of the starting points of $L''_i$ satisfy $0>y_1>\dots>y_s$.
After this step we obtain $L'''_i=L(u_i(\mathbf{w}'''),\phi'''_i)$ for $1\leq i\leq s$. The resulted positions of $u_i(\mathbf{w}''')$, satisfying the above two conditions, enable us to perform the last step.

In the last step we rotate the half lines $L'''_1,\dots,L'''_s$ counterclockwise successively such that their directions become horizontal. Then finally we obtain the half lines $L(u_i(\mathbf{w}^*),\phi^*_i)$ satisfying \eqref{eq-requirementforpictorialproof}.

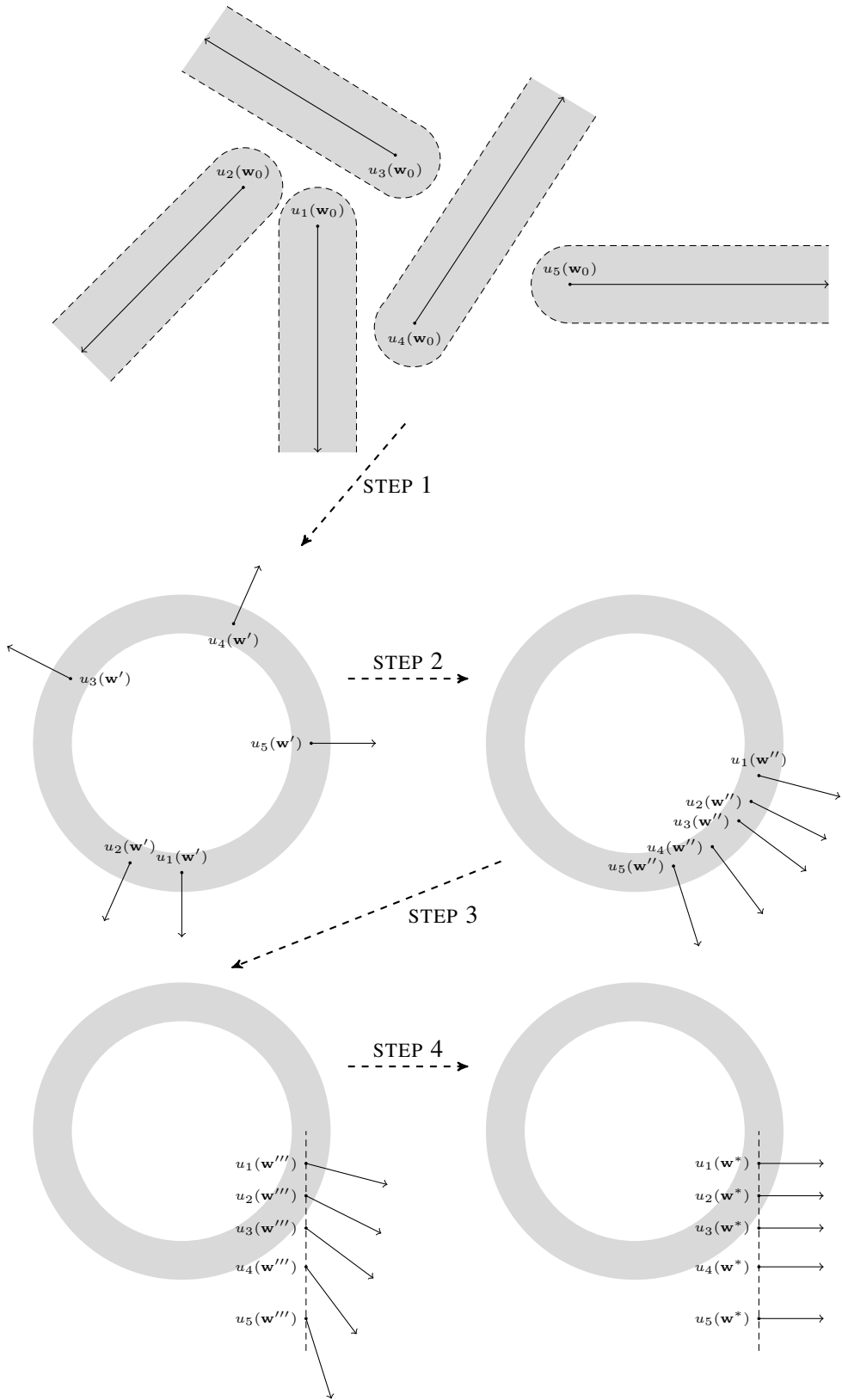
\begin{tikzpicture}
 [pile/.style={thick, ->, >=stealth', shorten <=2pt, shorten
    >=2pt}]
  \coordinate (u_1) at (0.6,0);
  \coordinate (u_2) at (-0.55,0.6);
  \coordinate (u_3) at (1.8,1.1);
  \coordinate (u_4) at (2.1,-1.5);
  \coordinate (u_5) at (4.5,-0.9);

  \coordinate (c_1) at (-1.5,-8);
  \coordinate (c_2) at (5.5,-8);
  \coordinate (c_3) at (-1.5,-14);
  \coordinate (c_4) at (5.5,-14);

    \draw[fill = gray!30, densely dashed] (0,-3.5)-- (0,0) arc [start angle = 180, end angle = 0, x radius = 6mm, y radius = 6mm]  -- (1.2,-3.5);

   \draw[arrows=->] (0.6,0) -- (0.6,-3.5);

   \filldraw (u_1) circle (0.5pt) node[above] {\tiny$u_1(\mathbf{w}_0)$};

   \draw[fill = gray!30, densely dashed] (-3.5,-1.5)-- (-1,1) arc [start angle = 140, end angle = -40, x radius = 6mm, y radius = 6mm]  -- (-2.6,-2.4);

   \draw[arrows=->] (-0.55,0.6) -- (-3.05,-1.95);

   \filldraw (u_2) circle (0.5pt) node[above] {\tiny$u_2(\mathbf{w}_0)$};

   \draw[fill = gray!30, densely dashed] (-0.8,3.4)-- (2.1,1.6) arc [start angle = 70, end angle = -110, x radius = 6mm, y radius = 6mm]  -- (-1.5,2.4);

    \draw[arrows=->] (1.8,1.1) -- (-1.15,2.9);

    \filldraw (u_3) circle (0.5pt) node[below] {\tiny$u_3(\mathbf{w}_0)$};

    \draw[fill = gray!30, densely dashed] (4.9,1.7)-- (2.5,-2) arc [start angle = -45, end angle = -225, x radius = 6mm, y radius = 6mm]  -- (3.9,2.3);

    \draw[arrows=->] (2.1,-1.5) -- (4.4,2);

    \filldraw (u_4) circle (0.5pt) node[below] {\tiny$u_4(\mathbf{w}_0)$};

    \draw[fill = gray!30, densely dashed] (8.5,-1.5)-- (4.5,-1.5) arc [start angle = -90, end angle = -270, x radius = 6mm, y radius = 6mm]  -- (8.5,-0.3);

    \draw[arrows=->] (4.5,-0.9) -- (8.5,-0.9);

    \filldraw (u_5) circle (0.5pt) node[above] {\tiny$u_5(\mathbf{w}_0)$};

  \fill[gray!30] (c_1) circle (23mm);
    \fill[white] (c_1) circle (17mm);
      
      \draw[arrows=->] (-1.5,-10) -- (-1.5,-11);
      \filldraw (-1.5,-10) circle (0.5pt) node[above] {\tiny$u_1(\mathbf{w}')$};

      \draw[arrows=->] (-2.3,-9.85) -- (-2.7,-10.75);
      \filldraw (-2.3,-9.85) circle (0.5pt) node[above] {\tiny$u_2(\mathbf{w}')$};

      \draw[arrows=->] (-3.22,-7) -- (-4.2,-6.5);
      \filldraw (-3.22,-7) circle (0.5pt) node[right] {\tiny$u_3(\mathbf{w}')$};

      \draw[arrows=->] (-0.7,-6.15) -- (-0.3,-5.25);
      \filldraw (-0.7,-6.15)circle (0.5pt) node[below] {\tiny$u_4(\mathbf{w}')$};

      \draw[arrows=->] (0.5,-8) -- (1.5,-8);
      \filldraw (0.5,-8)circle (0.5pt) node[left] {\tiny$u_5(\mathbf{w}')$};

      \fill[gray!30] (c_2) circle (23mm);
    \fill[white] (c_2) circle (17mm);

      \draw[arrows=->] (7.42,-8.5) -- (8.67804, -8.82761);
            \filldraw (7.42,-8.5) circle (0.5pt) node[above] {\tiny$u_1(\mathbf{w}'')$};

      \draw[arrows=->] (7.30,-8.9) -- (8.46276, -9.48138);
            \filldraw (7.30,-8.9) circle (0.5pt) node[left] {\tiny$u_2(\mathbf{w}'')$};

      \draw[arrows=->] (7.11,-9.2) -- (8.15233, -9.97689);
          \filldraw (7.11,-9.2) circle (0.5pt) node[left] {\tiny$u_3(\mathbf{w}'')$};

      \draw[arrows=->] (6.7,-9.6) -- (7.48, -10.64);
            \filldraw (6.7,-9.6) circle (0.5pt) node[left] {\tiny$u_4(\mathbf{w}'')$};

      \draw[arrows=->] (6.1,-9.9) -- (6.49147, -11.1397);
            \filldraw (6.1,-9.9) circle (0.5pt) node[left] {\tiny$u_5(\mathbf{w}'')$};

  \fill[gray!30] (c_3) circle (23mm);
    \fill[white] (c_3) circle (17mm);

      \draw[densely dashed] (0.42,-17.4) -- (0.42,-14);

      \draw[arrows=->] (0.42,-14.5) -- (1.67804, -14.82761);
            \filldraw (0.42,-14.5) circle (0.5pt) node[left] {\tiny$u_1(\mathbf{w}''')$};

      \draw[arrows=->] (0.42,-15) -- (1.58276, -15.58138);
            \filldraw (0.42,-15) circle (0.5pt) node[left] {\tiny$u_2(\mathbf{w}''')$};

      \draw[arrows=->] (0.42,-15.5) -- (1.46233, -16.27689);
          \filldraw (0.42,-15.5) circle (0.5pt) node[left] {\tiny$u_3(\mathbf{w}''')$};

      \draw[arrows=->] (0.42,-16.1) -- (1.2, -17.14);
            \filldraw (0.42,-16.1) circle (0.5pt) node[left] {\tiny$u_4(\mathbf{w}''')$};

      \draw[arrows=->] (0.42,-16.9) -- (0.81147, -18.1397);
            \filldraw (0.42,-16.9) circle (0.5pt) node[left] {\tiny$u_5(\mathbf{w}''')$};

    \fill[gray!30] (c_4) circle (23mm);
    \fill[white] (c_4) circle (17mm);

    \draw[densely dashed] (7.42,-17.4) -- (7.42,-14);

    \draw[arrows=->] (7.42,-14.5) -- (8.42,-14.5);
      \filldraw (7.42,-14.5) circle (0.5pt) node[left] {\tiny$u_1(\mathbf{w}^*)$};

      \draw[arrows=->] (7.42,-15) -- (8.42,-15);
      \filldraw (7.42,-15) circle (0.5pt) node[left] {\tiny$u_2(\mathbf{w}^*)$};

      \draw[arrows=->] (7.42,-15.5) -- (8.42,-15.5);
        \filldraw (7.42,-15.5) circle (0.5pt) node[left] {\tiny$u_3(\mathbf{w}^*)$};

      \draw[arrows=->] (7.42,-16.1) -- (8.42,-16.1);
        \filldraw (7.42,-16.1) circle (0.5pt) node[left] {\tiny$u_4(\mathbf{w}^*)$};

      \draw[arrows=->] (7.42,-16.9) -- (8.42,-16.9);
        \filldraw (7.42,-16.9) circle (0.5pt) node[left] {\tiny$u_5(\mathbf{w}^*)$};

    \draw[dashed, pile] (2,-3) -- (0.3,-5) node[midway,right] {\textsc{step 1}};
    \draw[dashed, pile] (1,-7) -- (3,-7) node[midway,above] {\textsc{step 2}};
    \draw[dashed, pile] (3.5,-9.8) -- (-0.8,-11.5) node[midway,right] {\hspace{0.5cm}\textsc{step 3}};
    \draw[dashed, pile] (1,-13) -- (3,-13) node[midway,above] {\textsc{step 4}};

\end{tikzpicture}
%\caption{} \label{fig:1}
\end{proof}

\newpage

\begin{proof}[More detailed proof of theorem \ref{sufficientcondition}]
	Without loss of generality, we assume 
	\ban
	2\pi>\phi_1>\dots>\phi_s>0,
	\nan
	and for any nonzero complex number $\lambda$, we choose its principal argument $\Arg(\lambda)\in[0,2\pi)$.
	
	We start by defining certain open subsets and sectors in $\mathbb{C}_{\lambda}$.
	\begin{itemize}
		\item[(i)] For $i\neq j$, let $d_{ij}$ be the distance between $L(u_i(\mathbf{w}_0),\phi_i)$ and $L(u_j(\mathbf{w}_0),\phi_j)$, which is nonzero by our assumption. Let
		\ban
		d_0:=\frac13\min\{d_{ij}:i,j=1,\dots,s,i\neq j\}.
		\nan
		Then the open subsets $W_i$ in $\bC_\lambda$, defined by
	    \ban
	    W_i:=\{\lambda_i\in\bC_\lambda:d(\lambda_i,L(u_i(\mathbf{w}_0),\phi_i))<d_0\},
	    \nan
	    are pairwise disjoint.
	    \item[(ii)] 
	    
	    We find sectors $S_i$, $S'_i$ and $S''_i$ in $\bC_\lambda$ satisfying:
	    \begin{enumerate}
	    	\item[(a)] $W_i\setminus S_i$ is a bounded subset of $\bC_\lambda$;
	     	\item[(b)] $S_i\neq\emptyset$ for $1\leq i\leq s$, $S_1,\dots,S_s$ are pairwise disjoint, all lying in the cut plane $\{\mathrm{Arg}(\lambda)\neq 0\}$, and $S_i$ lies in the counter-clockwise side of $S_{i+1}$;
	     	\item[(c)] $S''_i\neq\emptyset$, $S''_1,\dots,S''_s$ are pairwise disjoint, all lying in the fourth quadrant $\{\frac{3\pi}{2}<\mathrm{Arg}(\lambda)<2\pi\}$;
	     	\item[(d)]  For $1\leq i\leq s$, $S_i\cap S''_i=\emptyset$, and $S''_i$ lies in the counter-clockwise  side of $S_i$; 
	     	\item[(e)] $S'_i$ is the smallest open sector lying  in the cut plane $\{\mathrm{Arg}(\lambda)\neq 0\}$ that contains both $S_i$ and $S'_i$.
	     \end{enumerate}  
	     For example, putting
	    \[
	    \varepsilon'_0:=\frac13\min\big(\{2\pi-\phi_1\}\cup\{\phi_{i}-\phi_{i+1}:1\leq i\leq s-1\}\cup\{\phi_s\}\big),
	    \]
	     and take
	    \ban
	    \varepsilon_0<\min\{\varepsilon'_0,\frac{2\pi-\phi_1}{3s},\frac{\pi}{2(3s-1)}\},
	    \nan
	    one can  take the sectors $S_i$, $S'_i$ and $S''_i$ in $\bC_\lambda$ defined by
	    \begin{equation}\label{eq-def-Si-0}
	    		    S_i:=\{\lambda\neq0:|\Arg(\lambda)-\phi_i|<\varepsilon_0\},
	    \end{equation}
	    \begin{equation}\label{eq-def-Si-1}
	    	S'_i:=\{\lambda\neq0:\phi_i-\varepsilon_0<\Arg(\lambda)<2\pi-(3i-3)\varepsilon_0\},
	    \end{equation}
	    \begin{equation}\label{eq-def-Si-2}
	    	S''_i:=\{\lambda\neq0:|\Arg(\lambda)-\big(2\pi-(3i-2)\varepsilon_0\big)|<\varepsilon_0\}.
	    \end{equation}
    \end{itemize}
   By the choice of $\varepsilon_0$, these sectors meet all the above requirements (a) to (e). In fact the form of the definition (\ref{eq-def-Si-0}) implies (a); $\varepsilon<\varepsilon'_0$ implies (b); $\varepsilon_0< \frac{\pi}{2(3s-1)}$ implies (c); $\varepsilon_0<\frac{2\pi-\phi_1}{3s}$ implies (d).
   
Now let $W^{(0)}:=W_1\times\cdots\times W_s$. Then $W^{(0)}$ is a simply-connected open subset of $\bC^s\setminus\Delta$, and we identify $W^{(0)}$ with a neighborhood of $\mathbf{w}_0$ in $\mathcal C$. Recall that $W^{(0)}\setminus\mathcal B$ is a divisor. So we can find a path $\{\mathbf{w}_t\}_{0\leq t\leq1}$ in $W^{(0)}\cap\mathcal B$ such that $u_i(\mathbf{w}_1)\in W_i\cap S_i$ for each $i$. For each $t\in[0,1]$, since $u_i(\mathbf{w}_t)\in W_i$, it follows that $\phi_i\in A_i(\mathbf{w}_t)$. So from Corollary \ref{moveuu}, $y_i$ respects $(u_i,\Psi_i)$ with phase $\phi_i$ around $\mathbf{w}_1$. Since $u_i(\mathbf{w}_1)\in S_i$, it follows that $\phi_i$ and $\Arg(u_i(\mathbf{w}_1))$ are in the same connected component of $A_i(\mathbf{w}_1)$. So $y_i$ also respects $(u_i,\Psi_i)$ with phase $\Arg(u_i(\mathbf{w}_1))$ around $\mathbf{w}_1$.

Let $W^{(1)}:=S'_1\times S_2\times\cdots\times S_s$. Then $W^{(1)}$ is a simply-connected open subset of $\bC^s\setminus\Delta$, and we can identify $W^{(1)}$ with a neighborhood of $\mathbf{w}_1$ in $\mathcal C$. Recall that $W^{(1)}\setminus\mathcal B$ is a divisor. So we can find a path $\{\mathbf{w}_t\}_{1\leq t\leq1}$ in $W^{(1)}\cap\mathcal B$ such that $u_1(\mathbf{w}_2)\in S''_1$. Note that for each $t\in[1,2]$, we have $\Arg(u_i(\mathbf{w}_t))\in A_i(\mathbf{w}_t)$. So from Corollary \ref{moveuu}, $y_i$ respects $(u_i,\Psi_i)$ with phase $\Arg(u_i(\mathbf{w}_2))$ around $\mathbf{w}_2$.

Similarly, for $j=2,3,\dots,s-1$, assume that we have moved $\mathbf{w}_{j-1}$ to $\mathbf{w}_j$ via a path in $\mathcal B$ such that $u_i(\mathbf{w}_j)\in S''_i(1\leq i\leq j)$, $u_i(\mathbf{w}_j)\in S_i(j+1\leq i\leq s)$, and $y_i$ respects $(u_i,\Psi_i)$ with phase $\Arg(u_i(\mathbf{w}_j))$ around $\mathbf{w}_j$. Let $W^{(j)}:=S''_1\times\cdots S''_{j-1}\times S'_j\times S_{j+1}\times\cdots\times S_s$. Then $W^{(j)}$ is a simply-connected open subset of $\bC^s\setminus\Delta$, and we can identify $W^{(j)}$ with a neighborhood of $\mathbf{w}_j$ in $\mathcal C$. Recall that $W^{(j)}\setminus\mathcal B$ is a divisor. So we can find a path $\{\mathbf{w}_t\}_{j\leq t\leq j+1}$ in $W^{(j)}\cap\mathcal B$ such that $u_j(\mathbf{w}_{j+1})\in S''_j$. Note that for each $t\in[j,j+1]$, we have $\Arg(u_i(\mathbf{w}_t))\in A_i(\mathbf{w}_t)$. So from Corollary \ref{moveuu}, $y_i$ respects $(u_i,\Psi_i)$ with phase $\Arg(u_i(\mathbf{w}_{j+1}))$ around $\mathbf{w}_{j+1}$.

So we obtain $u_{s+1}\in W^{(s+1)}:=S''_1\times\cdots\times S''_s$ such that $y_i$ respects $(u_i,\Psi_i)$ with phase $\Arg(u_i(\mathbf{w}_{s+1}))$ around $\mathbf{w}_{s+1}$. Note that $W^{(j)}$ is a simply-connected open subset of $\bC^s\setminus\Delta$, and we can identify $W^{(s+1)}$ with a neighborhood of $\mathbf{w}_{s+1}$ in $\mathcal C$. Recall that $W^{(s+1)}\setminus\mathcal B$ is a divisor. So we can find a path $\{\mathbf{w}_t\}_{s+1\leq t\leq s+2}$ in $W^{(s+1)}\cap\mathcal B$ such that for $i=1,\dots,s-1$, we have $|u_{i+1}(\mathbf{w}_{s+2})|>|u_{i}(\mathbf{w}_{s+2})|$.  Note that for each $t\in[s+1,s+2]$, we have $\Arg(u_i(\mathbf{w}_t))\in A_i(\mathbf{w}_t)$. So from Corollary \ref{moveuu}, $y_i$ respects $(u_i,\Psi_i)$ with phase $\Arg(u_i(\mathbf{w}_{s+2}))$ around $\mathbf{w}_{s+2}$.

Now we have
\ban
&2\pi>\Arg(u_1(\mathbf{w}_{s+2}))>\dots>\Arg(u_s(\mathbf{w}_{s+2}))>\frac{3\pi}{2},&\\&
0<|u_1(\mathbf{w}_{s+2})|<\dots<|u_s(\mathbf{w}_{s+2})|.&
\nan
So we can check that $2\pi$ and $\Arg(u_i(\mathbf{w}_{s+2}))$ are in the same connected component of $A_i(\mathbf{w}_{s+2})$. So from Corollary \ref{moveuu}, $y_i$ respects $(u_i,\Psi_i)$ with phase $2\pi$ around $\mathbf{w}_{s+2}$. Then $w^*:=w_{s+2}$ satisfies (\ref{eq-requirementforpictorialproof}).
 This implies that $F$ satisfies Gamma conjecture II. 
\end{proof}

\if{
\subsection{Mutation}

Let $Z$ be a vector space equipped with a bilinear form $[\cdot,\cdot)$. For any $u$, $v\in Z$, the right muation of $u$ with respect to $v$ is defined to be
\ban
R_v(u):=u-[u,v)v,
\nan
and the left mutation of $u$ with respect to $v$ is defined to be
\begin{equation*}
	L_v(u):=u-[v,u)v.
\end{equation*}
Given an ordered basis $\mathcal{B}=(u_1,\dots,u_s)$ of $Z$, for $1\leq i\leq s-1$, its $i$-th right mutation is the ordered basis $R_i\mathcal{B}$ defined by
\begin{equation*}
	R_i\mathcal{B}=(u_1,\dots,u_{i-1},u_{i+1},R_{u_{i+1}}(u_i),u_{i+2},\dots,u_s),
\end{equation*}
and for $2\leq i\leq s$, its $i$-th left mutation is the ordered basis $L_i\mathcal{B}$ defined by
\begin{equation*}
	L_i\mathcal{B}=(u_1,\dots,u_{i-2},L_{u_{i-1}}(u_i),u_{i-1},u_{i+1},\dots,u_s).
\end{equation*}

\begin{lemma}
     Assume that $\bullet_\mathbf t$ is analytic and semisimple in a neighborhood $B$ of $\mathbf t_0\in H^*(F)$. Then for $V_1,V_2\in K(F)$, we have
     \ban
     \cZ^K_{B}(R_{V_1}(V_2))=R_{\cZ^K_{B}(V_1)}(\cZ^K_{B}(V_2)),\quad\cZ^K_{\mathbf t_0}(R_{V_1}(V_2))=R_{\cZ^K_{\mathbf t_0}(V_1)}(\cZ^K_{\mathbf t_0}(V_2)),\\
     \cZ^K_{B}(L_{V_1}(V_2))=L_{\cZ^K_{B}(V_1)}(\cZ^K_{B}(V_2)),\quad\cZ^K_{\mathbf t_0}(L_{V_1}(V_2))=L_{\cZ^K_{\mathbf t_0}(V_1)}(\cZ^K_{\mathbf t_0}(V_2)).
     \nan
\end{lemma}
\begin{proof}
	This follows from \eqref{paringsarematched}.
\end{proof}

\subsection{Mutation of full exceptional collection}

\begin{lemma}\label{frommutationofbundlestoFEC}	
Let $\mathcal E$ be an FEC of $\mathcal D^b(F)$, and $\mathcal B$ the corresponding ordered basis in $K(F)$. Then $R_i\mathcal{B}$ (resp. $L_i\mathcal{B}$) is the correspoinding ordered basis of $R_i\mathcal{E}$ (resp. $L_i\mathcal{E}$).
\end{lemma}
}\fi

\section{Chern characters of spinor bundles}\label{sec:chernclass-spinorbundle}
In this section, we compute the Chern characters of the spinor bundles over quadric hypersurfaces. We denote by $\qdc_N$ a smooth quadric hypersurface of dimension $N$ in $\mathbb{P}^{N+1}$. 
First we recall some facts on the integral cohomology ring of quadric hypersurfaces \cite[p. 315]{Fri83}. Denote by $\hpc^i$ be the pullback of the class of a codimension-$i$ linear subspace of the ambient projective space. Then for odd dimensional quadrics,
\begin{eqnarray*}
 H^*(\qdc_{2k+1},\mathbb{Z})=\mathbb{Z}1+\mathbb{Z}\hpc+\cdots+\mathbb{Z}\hpc^{k}+\mathbb{Z}\frac{1}{2}\hpc^{k+1}+\cdots+\mathbb{Z}\frac{1}{2}\hpc^{2k+1}.
\end{eqnarray*}
 For $\qdc_{4k+2}$, there  are two primitive classes $e_{2k+1}, e'_{2k+1}\in H^{2k+1}(\qdc_{4k+2},\mathbb{Z})$ such that
\begin{eqnarray*}
H^*(\qdc_{4k+2},\mathbb{Z})=\mathbb{Z}1+\mathbb{Z}\hpc+\cdots+\mathbb{Z}\hpc^{2k}
+\mathbb{Z}e_{2k+1}+\mathbb{Z}e_{2k+1}'+\mathbb{Z}\frac{1}{2}\hpc^{2k+2}+\cdots+\mathbb{Z}\frac{1}{2}\hpc^{4k+2},
\end{eqnarray*}
and
\begin{equation}\label{eq-cohring-4k+2}
\begin{split}
e_{2k+1}+e_{2k+1}'=\hpc^{2k+1}, &\quad e_{2k+1}\cdot e_{2k+1}=e_{2k+1}'\cdot e_{2k+1}'=0, \\
e_{2k+1}\cdot e_{2k+1}'=\frac{1}{2}\hpc^{4k+2}, &\quad \hpc\cdot e_{2k+1}=\hpc\cdot e_{2k+1}'=\frac{1}{2}\hpc^{2k+2}.
\end{split}
\end{equation}
For $\qdc_{4k}$, there  are two primitive classes $e_{2k}, e'_{2k}\in H^{2k}(\qdc_{4k},\mathbb{Z})$ such that 
\begin{eqnarray*}
H^*(\qdc_{4k},\mathbb{Z})=\mathbb{Z}1+\mathbb{Z}\hpc+\cdots+\mathbb{Z}\hpc^{2k-1}+\mathbb{Z}e_{2k}+\mathbb{Z}e_{2k}'+\mathbb{Z}\frac{1}{2}\hpc^{2k+1}+\cdots+\mathbb{Z}\frac{1}{2}\hpc^{4k},
\end{eqnarray*}
and
\begin{equation}\label{eq-cohring-4k}
\begin{split}
e_{2k}+e_{2k}'=\hpc^{2k}, &\quad e_{2k}\cdot e_{2k}=e_{2k}'\cdot e_{2k}'=\frac{1}{2}\hpc^{4k}, \\
e_{2k}\cdot e_{2k}'=0, &\quad \hpc\cdot e_{2k}=\hpc\cdot e_{2k}'=\frac{1}{2}\hpc^{2k+1}.
\end{split}
\end{equation}

Next we collect some facts  about the spinor bundles  on quadric hypersurfaces from  \cite[theorem 1.4, theorem 2.8, theorem 2.3]{Ott88}.
Denote the spinor bundle on $\qdc_{2k+1}$ by $S_{2k+1}$, and the two spinor bundles on $\qdc_{2k}$ by $S'_{2k}$ and $S''_{2k}$, $k\geq 1$. We have 
\begin{subequations}\label{eq-spinorbundles}
\begin{eqnarray}
& S_1=\mathcal{O}_{\mathbb{P}^1}(-1),\ S_2'=\mathcal{O}_{\mathbb{P}^1\times\mathbb{P}^1}(0,-1), S_2''=\mathcal{O}_{\mathbb{P}^1\times\mathbb{P}^1}(-1,0),\label{eq-sp-initial}\\
&\mathrm{rk} (S_{2k+1})=2^k,\  \mathrm{rk} (S'_{2k})=\mathrm{rk} (S''_{2k})=2^{k-1},\label{eq-sp-rank}\\
&S_{2k+1}|_{\qdc_{2k}}=S'_{2k}\oplus S''_{2k},\ S_{2k}'|_{\qdc_{2k-1}}=S_{2k}''|_{\qdc_{2k-1}}=S_{2k-1},\label{eq-sp-restrictions}\\
&S_{2k+1}^{\vee}\cong S_{2k+1}(1),\label{eq-sp-dual-1}\\
&S'^{\vee}_{4k}\cong S_{4k}'(1),\ S''^{\vee}_{4k}\cong S_{4k}''(1),\label{eq-sp-dual-2}\\
&S'^{\vee}_{4k+2}\cong S_{4k+2}''(1),\ S''^{\vee}_{4k+2}\cong S_{4k}'(1),\label{eq-sp-dual-3}\\
&
\begin{cases}
H^i(\qdc_{2k+1}, S_{2k+1})=0,& \mbox{for}\ 0\leq i\leq 2k+1,\\
H^i(\qdc_{2k}, S'_{2k})=0=H^i(\qdc_{2k}, S''_{2k}),& \mbox{for}\ 0\leq i\leq 2k.
\end{cases}\label{eq-sb-vanishing-1}
\end{eqnarray}
\end{subequations}
The vanishing (\ref{eq-sb-vanishing-1}) for $\mathcal{Q}_n$ in the range $0\leq i< n$ is \cite[theorem 2.3]{Ott88}, and for $i=n$ it follows from Serre duality and (\ref{eq-sp-dual-1}) and (\ref{eq-sp-dual-2}).

Note that for an even dimensional quadric $\qdc_{2k}$, the presented ring structure (\ref{eq-cohring-4k+2}) or (\ref{eq-cohring-4k})  cannot completely determine $e_k$ and $e'_k$. Namely, we have a choice to name which one of the two  $e_k$ or $e'_k$. This freedom will be deprived by the following theorem and eventually by (\ref{eq-sp-initial}), the restrictions (\ref{eq-sp-restrictions}), and a choice of isomorphism $\qdc_2\cong \mathbb{P}^1\times \mathbb{P}^1$.

\begin{thm}\label{computationofch}
\begin{enumerate}
\item For $k\geq 0$,
\begin{eqnarray*}
\mathrm{ch}(S_{2k+1})=\frac{2^{k+1}}{1+e^\hpc}.
\end{eqnarray*}
\item For $k\geq 1$, exchanging $e_{2k}$ and $e_{2k}'$ if necessary, we have
\begin{eqnarray*}
\mathrm{ch}(S_{2k}')=\frac{2^{k}}{1+e^\hpc}+ \frac{e_{k}-e_{k}'}{2}, \\
\mathrm{ch}(S_{2k}'')=\frac{2^{k}}{1+e^\hpc}- \frac{e_{k}-e_{k}'}{2},
\end{eqnarray*}
and
\begin{eqnarray*}
e_{k}-e_{k}'=(-1)^{k-1}\frac{c_{k}(S_{2k}')-c_{k}(S_{2k}'')}{(k-1)!}.
\end{eqnarray*}
\end{enumerate}
\end{thm}
\begin{proof}
By the isomorphisms (\ref{eq-sp-restrictions}), there exist rational numbers $a_i$, $i=0,1,2,\dots$, such that
\begin{multline*}
\mathrm{ch} (S_{2k}'\oplus S_{2k}'')=2^k a_0+2^{k}a_{1}\hpc+2^{k-1}a_{2}\hpc^2+2^{k-1}a_{3}\hpc^3\\
+\cdots+2a_{2k-2}\hpc^{2k-2}+2a_{2k-1}\hpc^{2k-1}+a_{2k}\hpc^{2k},\ \mbox{for any}\ k\geq 1,
\end{multline*}
and
\begin{multline*}
\mathrm{ch}(S_{2k+1})=2^k a_0+2^{k}a_{1}\hpc+2^{k-1}a_{2}\hpc^2+2^{k-1}a_{3}\hpc^3\\
+\cdots+2a_{2k-2}\hpc^{2k-2}+2a_{2k-1}\hpc^{2k-1}+a_{2k}\hpc^{2k}+a_{2k+1}\hpc^{2k+1},\ \mbox{for any}\ k\geq 0.
\end{multline*}
Taking the Chern characters of both sides of (\ref{eq-sp-dual-1}) we obtain
\begin{eqnarray*}
&&2^k a_0-2^{k}a_{1}\hpc+2^{k-1}a_{2}\hpc^2-2^{k-1}a_{3}\hpc^3\\
&&+\cdots+2a_{2k-2}\hpc^{2k-2}-2a_{2k-1}\hpc^{2k-1}+a_{2k}\hpc^{2k}-a_{2k+1}\hpc^{2k+1}\\
&=& e^{\hpc}(2^k a_0+2^{k}a_{1}\hpc+2^{k-1}a_{2}\hpc^2+2^{k-1}a_{3}\hpc^3\\
&&+\cdots+2a_{2k-2}\hpc^{2k-2}+2a_{2k-1}\hpc^{2k-1}+a_{2k}\hpc^{2k}+a_{2k+1}\hpc^{2k+1}),
\end{eqnarray*}
thus
\begin{eqnarray*}
a_{2k+1}=-\frac{2^{k-1} a_0}{(2k+1)!}-\frac{2^{k-1} a_1}{(2k)!}-\frac{2^{k-2} a_2}{(2k-1)!}-\frac{2^{k-2} a_3}{(2k-2)!}-\cdots-a_{2k}.
\end{eqnarray*}
By (\ref{eq-sb-vanishing-1}) and Riemann-Roch, we have
\begin{eqnarray*}
\int_{\qdc_{2k}}\mathrm{Td}(\qdc_{2k})\mathrm{ch} (S_{2k}'\oplus S_{2k}'')=0, & \int_{\qdc_{2k+1}}\mathrm{Td}(\qdc_{2k+1})\mathrm{ch}(S_{2k+1})=0.
\end{eqnarray*}

Since
\begin{eqnarray*}
\mathrm{Td}(\qdc_n)=\frac{\hpc^{n+2}}{(1-e^{-\hpc})^{n+2}}\Big/\frac{2\hpc}{1-e^{-2\hpc}}
=\frac{1+e^{-\hpc}}{2}\cdot \Big(\frac{\hpc}{1-e^{-\hpc}}\Big)^{n+1},
\end{eqnarray*}
we have, for $k\geq 1$,
\begin{eqnarray}\label{eq-coefficient-1}
\textrm{Coeff}_{\hpc^{2k}}\Bigg(\frac{1+e^{-\hpc}}{2}\cdot \Big(\frac{\hpc}{1-e^{-\hpc}}\Big)^{2k+1}\cdot(2^k a_0+2^{k}a_{1}\hpc+2^{k-1}a_{2}\hpc^2+2^{k-1}a_{3}\hpc^3\notag\\
+\cdots+2a_{2k-2}\hpc^{2k-2}+2a_{2k-1}\hpc^{2k-1}+a_{2k}\hpc^{2k})\Bigg)=0,
\end{eqnarray}
\begin{eqnarray}\label{eq-coefficient-2}
\textrm{Coeff}_{\hpc^{2k+1}}\Bigg(\frac{1+e^{-\hpc}}{2}\cdot \Big(\frac{\hpc}{1-e^{-\hpc}}\Big)^{2k+2}\cdot(2^k a_0+2^{k}a_{1}\hpc+2^{k-1}a_{2}\hpc^2+2^{k-1}a_{3}\hpc^3\notag\\
+\cdots+2a_{2k-2}\hpc^{2k-2}+2a_{2k-1}\hpc^{2k-1}+a_{2k}\hpc^{2k}+a_{2k+1}\hpc^{2k+1})\Bigg)=0,
\end{eqnarray}
where $\textrm{Coeff}_{t^i}g(t)$ denotes the coefficient of $t^i$ in a formal series $g(t)$ of $t$. 
One directly checks
\begin{eqnarray}\label{eq-coefficient-3}
\textrm{Coeff}_{h^0}\Bigg(\frac{1+e^{-\hpc}}{2}\cdot \Big(\frac{\hpc}{1-e^{-\hpc}}\Big)\cdot a_0\Bigg)=1,
\end{eqnarray}
and
\begin{eqnarray}\label{eq-coefficient-4}
\textrm{Coeff}_\hpc\Bigg(\frac{1+e^{-\hpc}}{2}\cdot \Big(\frac{\hpc}{1-e^{-\hpc}}\Big)^2\cdot (a_0+a_1\hpc)\Bigg)=0.
\end{eqnarray}
Let
\begin{eqnarray*}
f(x)=\frac{1+e^{-x}}{2}\cdot (a_0+a_{1}x+\frac{a_{2}}{2}x^2+\frac{a_{3}}{2}x^3+\cdots).
\end{eqnarray*}
Then (\ref{eq-coefficient-1})-(\ref{eq-coefficient-4}) amount to
\begin{eqnarray*}
\mathrm{Res}_{x=0}\frac{f(x)}{(1-e^{-x})^{k}}=
\begin{cases}
1, & k=1,\\
0, & k\geq 2.
\end{cases}
\end{eqnarray*}
Let $y=1-e^{-x}$, we compute by the residue theorem
\begin{eqnarray*}
&&\mathrm{Res}_{x=0}\frac{f(x)}{(1-e^{-x})^{k}}=\frac{1}{2\pi\sqrt{-1}}\oint_{|x|=\epsilon}\frac{f(x)dx}{(1-e^{-x})^{k}}\\
&=& \frac{1}{2\pi\sqrt{-1}}\oint_{|y|=\epsilon}\frac{f(x)e^xdy}{y^{k}}=\mathrm{Res}_{y=0}\frac{f(x)e^x}{y^{k}},
\end{eqnarray*}
Therefore
\begin{eqnarray*}
f(x)e^x=1,
\end{eqnarray*}
i.e.,
\begin{eqnarray*}
a_0+a_{1}x+\frac{a_{2}}{2}x^2+\frac{a_{3}}{2}x^3+\cdots=\frac{2}{1+e^x}.
\end{eqnarray*}
So
\begin{eqnarray*}
\mathrm{ch}(S_{2k+1})=\frac{2^{k+1}}{1+e^\hpc},
\end{eqnarray*}
and
\begin{eqnarray*}
\mathrm{ch} (S_{2k}')=\frac{2^k}{1+e^\hpc}+\alpha_{2k}', & \mathrm{ch} (S_{2k}'')=\frac{2^k}{1+e^\hpc}+\alpha_{2k}'',
\end{eqnarray*}
where $\alpha_{2k}', \alpha_{2k}''\in H_{\mathrm{prim}}^{2k}(\qdc_{2k},\mathbb{\qdc})$. Note that
\begin{eqnarray*}
c_{k}(S_{2k}')=(-1)^{k-1}(k-1)!\alpha_{2k}'+\mathrm{ambient\ class},\\
c_{k}(S_{2k}'')=(-1)^{k-1}(k-1)!\alpha_{2k}''+\mathrm{ambient\ class}.
\end{eqnarray*}
By Riemann-Roch, we have
\begin{eqnarray*}
&& \sum_{i=0}^{4k}(-1)^i\dim \mathrm{Ext}^{i}(S_{4k}',S_{4k}')\\
&=& \int_{\qdc_{4k}}\mathrm{Td}(\qdc_{4k})\mathrm{ch}(S_{4k}'\otimes S_{4k}'(1))\\
&=&2\textrm{Coeff}_{\hpc^{4k}}\Bigg(\frac{1+e^{-\hpc}}{2}\cdot \Big(\frac{\hpc}{1-e^{-\hpc}}\Big)^{4k+1}\cdot  \frac{2^{4k}}{(1+e^\hpc)^2}\cdot e^{\hpc}\Bigg)\\
&&+\int_{\qdc_{4k}}\alpha_{4k}'\cup \alpha_{4k}'\\
&=&2^{4k}\textrm{Coeff}_{\hpc^{4k}}\Bigg(\Big(\frac{\hpc}{1-e^{-\hpc}}\Big)^{4k+1}\cdot \frac{1}{1+e^\hpc}\Bigg)+\int_{\qdc_{4k}}\alpha_{4k}'\cup \alpha_{4k}'\\
&=&2^{4k} \mathrm{Res}_{y=0}\frac{1}{(2-y)y^{4k+1}}+\int_{\qdc_{4k}}\alpha_{4k}'\cup \alpha_{4k}'\\
&=&\frac{1}{2}+\int_{\qdc_{4k}}\alpha_{4k}'\cup \alpha_{4k}',
\end{eqnarray*}
and
\begin{eqnarray*}
&& \sum_{i=0}^{4k}(-1)^i\dim \mathrm{Ext}^{i}(S_{4k}',S_{4k}'')=\sum_{i=0}^{4k}(-1)^i\dim \mathrm{Ext}^{i}(S_{4k}'',S_{4k}')\\
&=& \int_{\qdc_{4k}}\mathrm{Td}(\qdc_{4k})\mathrm{ch}(S_{4k}'\otimes S_{4k}''(1))\\
&=&2\textrm{Coeff}_{\hpc^{4k}}\Bigg(\frac{1+e^{-\hpc}}{2}\cdot \Big(\frac{\hpc}{1-e^{-\hpc}}\Big)^{4k+1}\cdot  \frac{2^{4k}}{(1+e^\hpc)^2}\cdot e^{\hpc}\Bigg)\\
&&+\int_{\qdc_{4k}}\alpha_{4k}'\cup \alpha_{4k}''\\
&=&\frac{1}{2}+\int_{\qdc_{4k}}\alpha_{4k}'\cup \alpha_{4k}''.
\end{eqnarray*}
But by \cite{Kap88}, $S'_{4k}$ and $S''_{4k}$ form an exceptional pair, and thus
\[
\sum_{i=0}^{4k}(-1)^i\dim \mathrm{Ext}^{i}(S_{4k}',S_{4k}')=1,\quad
\sum_{i=0}^{4k}(-1)^i\dim \mathrm{Ext}^{i}(S_{4k}',S_{4k}'')=0.
\]
Then from the cohomology ring structure (\ref{eq-cohring-4k}) of $\qdc_{4k}$, we obtain
\begin{eqnarray*}
\mathrm{ch}(S'_{4k})=\frac{2^{2k}}{1+e^\hpc}\pm \frac{e_{2k}-e_{2k}'}{2}, & \mathrm{ch}(S''_{4k})=\frac{2^{2k}}{1+e^\hpc}\mp \frac{e_{2k}-e_{2k}'}{2}.
\end{eqnarray*}
For $\qdc_{4k+2}$, again by Riemann-Roch, we have
\begin{eqnarray*}
&& \sum_{i=0}^{4k+2}(-1)^i\dim \mathrm{Ext}^{i}(S_{4k+2}',S_{4k+2}')= \sum_{i=0}^{4k+2}(-1)^i \dim \mathrm{Ext}^{i}(S_{4k+2}'',S_{4k+2}'')\\
&=& \int_{\qdc_{4k}}\mathrm{Td}(\qdc_{4k+2})\mathrm{ch}(S_{4k+2}''\otimes S_{4k+2}'(1))\\
&=&2\textrm{Coeff}_{\hpc^{4k+2}}\Bigg(\frac{1+e^{-\hpc}}{2}\cdot \Big(\frac{\hpc}{1-e^{-\hpc}}\Big)^{4k+3}\cdot  \frac{2^{4k+2}}{(1+e^\hpc)^2}\cdot e^{\hpc}\Bigg)\\
&&+\int_{\qdc_{4k+2}}\alpha_{4k+2}'\cup \alpha_{4k+2}''\\
&=&\frac{1}{2}+\int_{\qdc_{4k+2}}\alpha_{4k+2}'\cup \alpha_{4k+2}'',
\end{eqnarray*}
\begin{eqnarray*}
&& \sum_{i=0}^{4k+2}(-1)^i\dim \mathrm{Ext}^{i}(S_{4k+2}',S_{4k+2}'')\\
&=& \int_{\qdc_{4k}}\mathrm{Td}(\qdc_{4k+2})\mathrm{ch}(S_{4k+2}''\otimes S_{4k+2}''(1))\\
&=&\frac{1}{2}+\int_{\qdc_{4k+2}}\alpha_{4k+2}''\cup \alpha_{4k+2}'',
\end{eqnarray*}
and 
\begin{eqnarray*}
&& \sum_{i=0}^{4k+2}(-1)^i\dim \mathrm{Ext}^{i}(S_{4k+2}'',S_{4k+2}')\\
&=& \int_{\qdc_{4k}}\mathrm{Td}(\qdc_{4k+2})\mathrm{ch}(S_{4k+2}'\otimes S_{4k+2}'(1))\\
&=&\frac{1}{2}+\int_{\qdc_{4k+2}}\alpha_{4k+2}'\cup \alpha_{4k+2}'.
\end{eqnarray*}
So by (\ref{eq-cohring-4k+2}) we get
\begin{eqnarray*}
\mathrm{ch}(S_{4k+2}')=\frac{2^{2k+1}}{1+e^\hpc}\pm \frac{e_{2k+1}-e_{2k+1}'}{2}, &
\mathrm{ch}(S_{4k+2}'')=\frac{2^{2k+1}}{1+e^\hpc}\mp \frac{e_{2k+1}-e_{2k+1}'}{2}.
\end{eqnarray*}
\end{proof}

\section{Convergence of quantum cohomology for quadrics}\label{sec:convergence}

Let $\qdc$ be a smooth quadric hypersurface in a projective space. The main result of this section is the following theorem.
\begin{thm}\label{convergencofquadrics}
The quantum cohomology of $\qdc$ is analytic in a neighborhood of $\mathbf0\in H^*(\qdc)$.
\end{thm}

This theorem follows from Proposition \ref{prop-oddconvergent} and \ref{prop-evenconvergent}. The growth of invariants with only ambient classes is known from \cite[Theorem 1]{Zin14}. So the odd-dimensional case is relatively easy. For the even-dimensional case, we use deformation invariance of GW and WDVV to 
deduce some estimates on the growth of invariants with primitive classes (Proposition \ref{thm-estimate-primitive}).

\subsection{Odd-dimensional case}

\begin{prop}\label{prop-oddconvergent}
	If $\dim_\bC \qdc$ is odd, then the quantum cohomology of $\qdc$ is analytic in a neighborhood of $\mathbf0\in H^*(\qdc)$.
\end{prop}
\begin{proof}
	Suppose that $\dim_\bC \qdc=2N+1$. Then $\hpc^0,\hpc^1,\dots,\hpc^{2N+1}$ is a homogeneous basis of $H^*(\qdc)$. Recall that $\Eff(\qdc)$ can be naturally identified with $\bZ_{\geq0}$. From \cite[Theorem 1]{Zin14}, there is a positive number $C$ such that
	\ban
	|\<\prodl_{i=1}^{2N+1}\tau_0(\hpc^i)^{n_i}\>_d^\qdc|\leq(\suml_{i=1}^{2N+1}n_i)!\cdot C^{\suml_{i=1}^{2N+1}n_i+d},\quad\forall n_1,\dots,n_{2N+1}\geq0.
	\nan
	From the degree axiom, 
	\ban
	\<\prodl_{i=1}^{2N+1}\tau_0(\hpc^i)^{n_i}\>_d^\qdc\neq0\Rightarrow d=\frac1{2N+1}\suml_{i=1}^{2N+1}(i-1)n_i+\frac3{2N+1}-1.
	\nan
	Therefore we have
	\ban
	|\<\prodl_{i=1}^{2N+1}\tau_0(\hpc^i)^{n_i}\>_d^\qdc\prodl_{i=1}^{2N+1}\frac{(t^i)^{n_i}}{n_i!}|\leq{\binom{n_1+\cdots+n_{2N+1}}{n_1,\dots,n_{2N+1}}}\cdot\prodl_{i=1}^{2N+1}\Big(C^{1+\frac{i-1}{2N+1}}|t^i|\Big)^{n_i}\cdot C^{\frac{3}{2N+1}-1}.
	\nan
	So for $\mathbf t=\suml_{i=0}^{2N+1}t^i\hpc^i\in H^*(\qdc)$ satisfying $\suml_{i=1}^{2N+1}C^{1+\frac{i-1}{2N+1}}|t^i|<1$, we have
	\ban
	&&\suml_{n_{1},\dots,n_{2N+1}\geq0}\bigg|\suml_{d>0}\<\prodl_{i=1}^{2N+1}\tau_0(\hpc^{i})^{n_{i}}\>_{d}^F\bigg|\prodl_{i=1}^{2N+1}\frac{|t^{i}|^{n_{i}}}{n_{i}!}\\
	&\leq&\suml_{n_{1},\dots,n_{2N+1}\geq0}{\binom{n_1+\cdots+n_{2N+1}}{n_1,\dots,n_{2N+1}}}\cdot\prodl_{i=1}^{2N+1}\Big(C^{1+\frac{i-1}{2N+1}}|t^i|\Big)^{n_i}\cdot C^{\frac3{2N+1}-1}\\
	&=&\frac{C^{\frac3{2N+1}-1}}{1-\suml_{i=1}^{2N+1}C^{1+\frac{i-1}{2N+1}}|t^i|}.
	\nan
	This implies that $\cF_0^\qdc(\mathbf t)$ is a absolutely convergent power series of $t^0,t^1,\dots,t^{2N+1}$ for $\mathbf t$ sufficiently close to $\mathbf 0\in H^*(\qdc)$.
\end{proof}

\subsection{Even-dimensional case}
In this subsection, let $\qdc=\qdc_{2N}$ be a smooth quadric hypersurface in $\mathbb{P}^{2N+1}$ for a fixed $N\geq2$. Then $\Eff(\qdc)$ is naturally identified with $\bZ_{\geq0}$. We denote by $\wp\in H^{2N}(\qdc)$ a primitive class such that $\wp\cup\wp=\hpc^{2N}$. Then $\mathbf{B}
:=\{\mathbbm1,\hpc,\hpc^2,\dots,\hpc^{2N},\wp\}$ is a basis of $H^*(\qdc)$. For $\gamma_1,\dots,\gamma_{n}\in \mathbf{B}$, if the invariant $\langle \gamma_1,\dots,\gamma_n\rangle_{d}^\qdc$ is non-zero, then the degree axiom implies that
\begin{equation}\label{eq-dim-constraint}
	\mathrm{deg}_{\mathbb{C}}(\gamma_1)+\cdots+\mathrm{deg}_{\mathbb{C}}(\gamma_n)=2N-3+n+2Nd.
\end{equation}
Here $\mathrm{deg}_\bC(\gamma)=p$ for $\gamma\in H^{2p}(\qdc)$, and we have identified the effective curve class with a non-negative integer $d$.
For a nonzero invariant $\langle \gamma_1,\dots,\gamma_n\rangle_{d}^\qdc$ with $\gamma_i\in\mathbf{B}$, the class $d$ is determined by $\gamma_1,\dots,\gamma_n$. So we will often use the abbreviated notation 
\[
\langle \gamma_1,\dots,\gamma_n\rangle:=\langle \gamma_1,\dots,\gamma_n\rangle_{d}:=\langle \gamma_1,\dots,\gamma_n\rangle_{d}^\qdc,\mbox{ for }\gamma_i\in\mathbf{B}.
\]

\begin{lemma}
	An invariant of the form
	\[
	\langle \underbrace{\wp,\dots,\wp}_m, \hpc^{k_1},\dots,\hpc^{k_n}\rangle
	\]
	is zero unless $m$ is even. 
\end{lemma}
\begin{proof}
Note that there is a monodromy action on $H^*(\qdc)$ which preserves $\hpc^k$ and transforms $\wp$ to $-\wp$ (see e.g. \cite[Proposition 5.2]{Del73}). Then the result follows from the deformation invariance of Gromov-Witten invariants \cite[Theorem 4.2']{LiTi98}.
\end{proof}

 %Now by the dimension constraint (\ref{eq-dim-constraint}), the monodromy reason, and the fundamental class axiom we have the followng key observation:
 
\begin{lemma}\label{lem-key-1}
Suppose $n\geq 4$. Given $\gamma_1,\dots,\gamma_{n-1}\in \mathbf{B}$, there exists at most one class $\gamma_n$ in $\mathbf{B}$ such that 
\[
\langle \gamma_1,\dots,\gamma_n\rangle\neq 0.
\]
\end{lemma}
\begin{proof} Suppose $\langle \gamma_1,\dots,\gamma_n\rangle\neq 0$. If $\wp$ appears in $\gamma_1,\dots,\gamma_{n-1}$ odd times, then $\gamma_n$ must be $\wp$. If $\wp$ appears in $\gamma_1,\dots,\gamma_{n-1}$ even times, then $\gamma_n\in \{\mathbbm1,\hpc,\hpc^2,\dots,\hpc^{2N}\}$, and there is exactly one such class matching the dimension constraint \eqref{eq-dim-constraint}.
\end{proof}

\begin{lemma}
\begin{equation}\label{eq-vanish-1}
\langle\underbrace{\wp,\dots,\wp}_{n\ \wp's}\rangle_{d}=0,
\end{equation}
\begin{equation}\label{eq-vanish-2}
\langle \underbrace{\wp,\dots,\wp}_{m\ \wp's},\underbrace{\hpc,\dots,\hpc}_{n\ h's}\rangle_{d}=0.
\end{equation}
\end{lemma}
\begin{proof}
The first equality holds for even $n$ by the dimension constraint \eqref{eq-dim-constraint}, and for odd $n$ by the monodromy reason. If $d>0$ the second equality follows  from the first equality by the divisor axiom. If $d=0$ the second equality follows  from $\wp\cup h=0$ when $m\geq 1$, and from the dimension reason when $m=0$.
\end{proof}

Now we recall Zinger's estimate on invariants with only ambient insertions. Define
\[
D(k,n)=\frac{k-2N+3-n}{2N}.
\]
We have
\begin{equation}\label{eq-D}
	D(k_1,n_1)+D(k_2,n_2)=D(k_1+k_2-2N,n_1+n_2-3).
\end{equation}

\begin{prop}\label{thm-zinger}
There exists $C>0$, such that for all $n\geq 1$ and any $i_1,\dots,i_n\geq 0$,
\begin{equation}\label{eq-estimate-0}
|\langle \hpc^{i_1},\dots,\hpc^{i_n}\rangle|\leq n!\cdot C^{n+D(\sum_{j=1}^n i_j,n)}.
\end{equation}
\end{prop}
\begin{proof}
It suffices to consider non-zero invariants. For a non-zero invariant $\<\hpc^{i_1},\dots,\hpc^{i_n}\>_d$, we use the dimension constraint \eqref{eq-dim-constraint} to get $d=D(\sum_{j=1}^n i_j,n)$. Then the inequality follows from \cite[Theorem 1]{Zin14}.
\end{proof}

\begin{lemma}\label{nonzerothreepointinvariants}
	The non-zero 3-point invariants with insertions in $\mathbf{B}$ are \begin{equation}\label{eq-3pointsinvariants}
	\begin{cases}
	\langle \wp,\wp, \mathbbm1\rangle=2,\ \langle \wp,\wp, \hpc^{2N}\rangle=-4,\\
	\langle \hpc^{i},\hpc^{j},\hpc^{k}\rangle=2,\ \mbox{if}\ i+j+k=2N,\\
	\langle \hpc^{i},\hpc^{2N-i},\hpc^{2N}\rangle=4,\ \mbox{if}\ 1\leq i\leq 2N-1,\\
	\langle \hpc^{i},\hpc^{j},\hpc^{k}\rangle=8,\ \mbox{if}\ i+j+k=4N,\ \mbox{and}\ i,j,k\neq 2N,\\
	\langle \hpc^{2N},\hpc^{2N},\hpc^{2N}\rangle=8.
	\end{cases}
	\end{equation}
\end{lemma}
\begin{proof}
	Let $\circ:=\bullet_{\mathbf t=\mathbf0}$ be the quantum product on $H^*(\qdc)$ specialized at $\mathbf t=\mathbf0$. Then $3$-point invariants are coefficients of the product $\circ$. From \cite[Proposition 1, (1.6), (2.2), (2.3)]{Bea95}, the specialized quantum product $\circ$ is given by
	\ban
	\hpc\circ \hpc^i&=&\left\{\begin{array}{cc}
		\hpc^{i+1},&0\leq i\leq 2N-2,\\
		\hpc^{2N}+2\hpc^0,&i=2N-1,\\
		2\hpc,&i=2N,
	\end{array}\right.\\
	\hpc\circ\wp&=&0,\\
	\wp\circ\wp&=&\hpc^{2N}-2\hpc^0.
	\nan
	We can use these identities to obtain all products. For example, we have
	\ban
	\hpc^{2N}\circ \hpc^{2N}&=&(\hpc\circ \hpc^{2N-1}-2\hpc^0)\circ \hpc^{2N}\\
	&=&\hpc^{2N-1}\circ(\hpc\circ \hpc^{2N})-2\hpc^{2N}\\
	&=&\hpc^{2N-1}\circ2\hpc-2\hpc^{2N}\\
	&=&2(\hpc^{2N}+2\hpc^0)-2\hpc^{2N}\\
	&=&4\hpc^0.
	\nan
	For one more example, for $1\leq i\leq 2N-1$, we have
	\[\hpc^i\circ\wp=\hpc^{i-1}\circ \hpc\circ\wp=0,\]
	and for $i=2N$, we have
	\ban
	\hpc^{2N}\circ\wp&=&(\hpc\circ \hpc^{2N-1}-2\hpc^0)\circ\wp=\hpc^{2N-1}\circ(\hpc\circ\wp)-2\wp=-2\wp.
	\nan
	We leave the rest cases to interested readers.
\end{proof}

%\subsubsection{\textcolor{red}{An estimate on invariants with primitive insertions}}

Denote $\xi_i=\hpc^i$ for $0\leq i\leq 2N$, and $\xi_{2N+1}=\wp$. Denote the dual coordinates by $t^0,\dots,t^{2N},t^{2N+1}$. 
First let us see how to compute all the genus-zero primary invariants inductively from invariants with only ambient insertions, and the  WDVV equation \eqref{generalWDVV}
\begin{equation}\label{eq-WDVV1}
\sum_{\alpha=0}^{2N+1}\sum_{\beta=0}^{2N+1}\partial_{t^{2N+1}}\partial_{t^{2N+1}}\partial_{t^\alpha}\mathcal F_{0}^\qdc
g^{\alpha \beta}\partial_{t^\beta}\partial_{t^i}\partial_{t^j}\mathcal F_{0}^\qdc=\sum_{\alpha=0}^{2N+1}\sum_{\beta=0}^{2N+1}\partial_{t^{2N+1}}\partial_{t^i}\partial_{t^\alpha}\mathcal F_{0}^\qdc g^{\alpha \beta}\partial_{t^\beta}\partial_{t^{2N+1}}\partial_{t^j}\mathcal F_{0}^\qdc,
\end{equation}
where $0\leq i,j\leq 2N$, and 
\begin{equation}\label{eq-gcoefficients}
	\begin{cases}
	g^{\alpha \beta}=\frac{1}{2}\delta_{\alpha+\beta,2N},& \mbox{for}\ 0\leq \alpha,\beta\leq 2N,\\
	g^{2N+1,\alpha}=\frac{1}{2}\delta_{\alpha,2N+1},& \mbox{for}\ 0\leq \alpha\leq 2N+1.
	\end{cases}
\end{equation}
By (\ref{eq-vanish-1}) and (\ref{eq-vanish-2}), it suffices  to compute invariants of the form
\begin{equation}\label{eq-invariants-with-primitive}
	\langle \underbrace{\wp,\dots,\wp}_{2m}, \hpc^{k_1},\dots,\hpc^{k_n},\hpc^{l}\rangle_{d},
\end{equation}
which satisfies (\ref{eq-dim-constraint}), where $2m+n+1\geq 4$, $m\geq 1$ and $l\geq 2$. For a subset $S\subset \{1,2,\dots,2m-2+n\}$, we introduce a convenient notation
\[
\xi_S=\ \mbox{the $|S|$-tuple}\ \underbrace{\wp,\dots,\wp}_{|S\cap \{1,\dots,2m-2\}|}, \hpc^{k_{i_1}},\dots,\hpc^{k_{i_\sigma}},
\]
and 
\[
\mathrm{deg}_\bC(\xi_S)=\mathrm{deg}_\bC(\wp)\cdot |S\cap \{1,\dots,2m-2\}|+k_{i_1}+\dots+k_{i_\sigma},
\]
where $\{i_1,\dots,i_\sigma\}=S\cap\{2m-2+1,\dots,2m-2+n\}$.

\begin{lemma}
Suppose $2\leq l\leq 2N$, and take $1\leq i,j\leq 2N-1$ such that $i+j=l$. Then 
\begin{multline}\label{eq-WDVV3}
\langle \underbrace{\wp,\dots,\wp}_{2m}, \hpc^{k_1},\dots,\hpc^{k_n},\hpc^l\rangle \\
=- \sum_{S\subsetneqq \{1,\dots,2m-2+n\}}\sum_{\alpha=0}^{2N+1}\sum_{\beta=0}^{2N+1}\langle \wp,\wp,\xi_S,\xi_\alpha\rangle g^{\alpha \beta}\langle \xi_\beta, \xi_{\{1,\dots,2m-2+n\}-S},\hpc^i,\hpc^j\rangle\\
+ \sum_{\begin{subarray}{c}S\subset \{1,\dots,2m-2+n\}\\ \emptyset \neq S\neq\{1,\dots,2m-2+n\} \end{subarray}}\sum_{\alpha=0}^{2N+1}\sum_{\beta=0}^{2N+1}
\langle \wp,\xi_S, \hpc^i,\xi_\alpha\rangle g^{\alpha \beta}\langle \xi_{\beta}, \xi_{\{1,\dots,2m-2+n\}-S},\wp,\hpc^j\rangle.
\end{multline}
\end{lemma}
\begin{proof}
 Taking the coefficients of $(t^{2N+1})^{2m-2}t^{k_1}\dots t^{k_n}$ on both sides of (\ref{eq-WDVV1}), we obtain
\begin{subequations}\label{eq-WDVV2}
\begin{eqnarray}
&& \sum_{\alpha=0}^{2N+1}\sum_{\beta=0}^{2N+1}\langle \underbrace{\wp,\dots,\wp}_{2m}, \hpc^{k_1},\dots,\hpc^{k_n},\xi_{\alpha}\rangle g^{\alpha \beta}\langle \xi_{\beta},\hpc^i,\hpc^j\rangle\label{subeq-WDVV2-1}\\
&&+ \sum_{S\subsetneqq \{1,\dots,2m-2+n\}}\sum_{\alpha=0}^{2N+1}\sum_{\beta=0}^{2N+1}
\langle \wp,\wp,\xi_S,\xi_\alpha\rangle g^{\alpha \beta}\langle \xi_\beta, \xi_{\{1,\dots,2m-2+n\}-S},\hpc^i,\hpc^j\rangle\notag\\
&=&  \sum_{\alpha=0}^{2N+1}\sum_{\beta=0}^{2N+1} \langle \underbrace{\wp,\dots,\wp}_{2m-1},\hpc^{k_1},\dots,\hpc^{k_n}, \hpc^{i},\xi_{\alpha}\rangle g^{\alpha \beta}\langle \xi_{\beta},\wp,\hpc^j\rangle\label{subeq-WDVV2-2}\\
&&+ \sum_{\alpha=0}^{2N+1}\sum_{\beta=0}^{2N+1} \langle \wp, \hpc^{i},\xi_{\alpha}\rangle g^{\alpha \beta}\langle \xi_{\beta},\underbrace{\wp,\dots,\wp}_{2m-1},\hpc^{k_1},\dots,\hpc^{k_n},,\hpc^j\rangle\label{subeq-WDVV2-3}\\
&&+ \sum_{\begin{subarray}{c}S\subset \{1,\dots,2m-2+n\}\\ \emptyset \neq S\neq\{1,\dots,2m-2+n\} \end{subarray}}\sum_{\alpha=0}^{2N+1}\sum_{\beta=0}^{2N+1}
\langle \wp,\xi_S, \hpc^i,\xi_\alpha\rangle g^{\alpha \beta}\langle \xi_{\beta}, \xi_{\{1,\dots,m-2+n\}-S},\wp,\hpc^j\rangle.\notag
\end{eqnarray}
\end{subequations}
Since $2\leq l\leq 2N$, we can choose $1\leq i,j\leq 2N-1$ such that $i+j=l$. Then by (\ref{eq-3pointsinvariants}), 
\[
\langle \xi_{\beta},\hpc^i,\hpc^j\rangle=\begin{cases} 
2 \delta_{\beta,2N-l}, & \mbox{if}\ 2\leq l\leq 2N-1,\\
 2 \delta_{\beta,0}+4 \delta_{\beta,2N},    & \mbox{if}\ l=2N.
 \end{cases}
\]
and
\[
\langle \wp,\hpc^i,\xi_\alpha\rangle=0,\ \forall\ 0\leq \alpha\leq 2N+1.
\]
So (\ref{subeq-WDVV2-2}) and (\ref{subeq-WDVV2-3}) vanish, and by  lemma \ref{lem-key-1} and (\ref{eq-gcoefficients}), the sum (\ref{subeq-WDVV2-1}) equals
\[
\langle \underbrace{\wp,\dots,\wp}_{2m}, \hpc^{k_1},\dots,\hpc^{k_n},\hpc^l\rangle.
\]
So we obtain (\ref{eq-WDVV3}).
\end{proof}

The invariants appearing on the right handside of (\ref{eq-WDVV3}), either have less $\wp$, or have the same number of $\wp$ but have less ambient insertions than that of the left handside. So we can inductively compute all the invariants of the form (\ref{eq-invariants-with-primitive}). Now we are ready to make the following estimate.
\begin{thm}\label{thm-estimate-primitive}
There exists $C>0$ such that for any $n\geq 3$, and for $\gamma_i\in\mathtt{B}$, $1\leq i\leq n$, 
\begin{equation}\label{eq-estimate-1}
|\langle \gamma_1,\dots,\gamma_n\rangle|\leq n! C^{n+D(\sum_{i=1}^n \mathrm{deg}_\bC(\gamma_i),n)}.
\end{equation}
\end{thm}
\begin{proof}
Suppose there are exactly $2m$ $\wp$'s in $\gamma_1,\dots,\gamma_n$, and without loss of generality we assume $n \geq 4$.
We prove (\ref{eq-estimate-1}) by induction on  $m\geq0$, and on the length $n$.  When $m=0$, we have Zinger's estimate (\ref{eq-estimate-0})
\begin{equation*}
|\langle \gamma_1,\dots,\gamma_n\rangle|\leq n! C^{n+D(\sum_{i=1}^n \mathrm{deg}_\bC(\gamma_i),n)},
\end{equation*}
for some $C>0$. Replacing $C$ by $C'$ such that 
\[
\left(\frac{C'}{C}\right)^n\geq n(n-1)(n-2),\ \forall n\geq 3,
\]
we have
\begin{equation*}
|\langle \gamma_1,\dots,\gamma_n\rangle|\leq (n-3)! C'^{n+D(\sum_{i=1}^n \mathrm{deg}_\bC(\gamma_i),n)},
\end{equation*}
Replacing $C'$ by a larger one again, denoted still by $C$, we have 
\begin{equation*}
|\langle \gamma_1,\dots,\gamma_n\rangle|\leq (n-3)! C^{n-3+D(\sum_{i=1}^n \mathrm{deg}_\bC(\gamma_i),n)}.
\end{equation*}
Suppose we are given an invariant of the form (\ref{eq-invariants-with-primitive}); recall that in (\ref{eq-invariants-with-primitive}) we have assumed $l\geq 2$.
Now we use (\ref{eq-WDVV3}) to estimate (\ref{eq-invariants-with-primitive}), inductively from  the invariants with less $\wp$'s or with the same number of $\wp$ but less ambient insertions  than that of (\ref{eq-invariants-with-primitive}).  Therefore
\begin{subequations}\label{eq-estimate-2}
\begin{eqnarray}
&&|\langle \underbrace{\wp,\dots,\wp}_{2m}, \hpc^{k_1},\dots,\hpc^{k_n},\hpc^l\rangle| \notag\\
&\leq & \sum_{S\subsetneqq \{1,\dots,2m-2+n\}}\sum_{\alpha=0}^{2N+1}\sum_{\beta=0}^{2N+1}|\langle \wp,\wp,\xi_S,\xi_\alpha\rangle g^{\alpha \beta}\langle \xi_\beta, \xi_{\{1,\dots,2m-2+n\}-S},\hpc^i,\hpc^j\rangle|\notag\\
&&+ \sum_{\begin{subarray}{c}S\subset \{1,\dots,2m-2+n\}\\ \emptyset \neq S\neq\{1,\dots,2m-2+n\} \end{subarray}}\sum_{\alpha=0}^{2N+1}\sum_{\beta=0}^{2N+1}
|\langle \wp,\xi_S, \hpc^i,\xi_\alpha\rangle g^{\alpha \beta}\langle \xi_{\beta}, \xi_{\{1,\dots,2m-2+n\}-S},\wp,\hpc^j\rangle|\notag\\
&= & |\langle \xi_{2N}, \xi_{\{1,\dots,2m-2+n\}},\hpc^i,\hpc^j\rangle|\label{eq-estimate-2-0}\\
&&+\sum_{\begin{subarray}{c}S\subset \{1,\dots,2m-2+n\}\\ \emptyset \neq S\neq\{1,\dots,2m-2+n\} \end{subarray}}\sum_{\alpha=0}^{2N+1}\sum_{\beta=0}^{2N+1}|\langle \wp,\wp,\xi_S,\xi_\alpha\rangle g^{\alpha \beta}\langle \xi_\beta, \xi_{\{1,\dots,2m-2+n\}-S},\hpc^i,\hpc^j\rangle|\label{eq-estimate-2-1}\\
&&+ \sum_{\begin{subarray}{c}S\subset \{1,\dots,2m-2+n\}\\ \emptyset \neq S\neq\{1,\dots,2m-2+n\} \end{subarray}}\sum_{\alpha=0}^{2N+1}\sum_{\beta=0}^{2N+1}
|\langle \wp,\xi_S, \hpc^i,\xi_\alpha\rangle g^{\alpha \beta}\langle \xi_{\beta}, \xi_{\{1,\dots,2m-2+n\}-S},\wp,\hpc^j\rangle|,\label{eq-estimate-2-2}
\end{eqnarray}
\end{subequations}
where we have used (\ref{eq-3pointsinvariants}) to obtain (\ref{eq-estimate-2-0}). Moreover, as to the sums (\ref{eq-estimate-2-1}) we note that by lemma \ref{lem-key-1}, for every $S\subset\{1,\dots,2m-2+n\}$, there is at most one pair $(\alpha,\beta)\in \{0,1,\dots,2N+1\}^2$ such that
\[
\langle \wp,\wp,\xi_S,\xi_\alpha\rangle g^{\alpha \beta}\langle \xi_\beta, \xi_{\{1,\dots,2m-2+n\}-S},\hpc^i,\hpc^j\rangle\neq 0,
\]
and for such a pair $(\alpha,\beta)$, due to (\ref{eq-gcoefficients}) we have
\begin{eqnarray}\label{eq-deg-sum-1}
&&\big(2\mathrm{deg}_\bC(\wp)+\mathrm{deg}_\bC(\xi_S)+\mathrm{deg}_\bC(\alpha)
\big)\notag\\
&&+\big(\mathrm{deg}_\bC(\beta)+\mathrm{deg}_\bC(\xi_{\{1,\dots,2m-2+n\}-S})+\mathrm{deg}_\bC(\hpc^i)+\mathrm{deg}_\bC(\hpc^j)
\big) \notag\\
&=& 2m N+k_1+\dots +k_n+l+2N.
\end{eqnarray}
The same remark holds for the  sums (\ref{eq-estimate-2-2}).

\textbf{Claim} 1: 
\begin{eqnarray}\label{eq-claim1}
&& |\wp,\wp,\hpc^{k_1},\dots,\hpc^{k_n},\hpc^l\rangle|\notag\\
&\leq& \frac{1}{4n-2}\binom{2n}{n} n! (2C)^{n+D(2N+\sum_{a=1}^n k_a+l,n+3)}.
\end{eqnarray}

Here we set the factors
\[
a_n:=\frac{1}{4n-2}\binom{2n}{n},\ n\geq 1,
\]
so that the generating function 
\[
f(X)=\sum_{n=1}^{+\infty}a_n X^n=\frac{1-\sqrt{1-4X}}{2}
\]
satisfies
\[
f(X)^2-f(X)+X=0,
\]
and thus 
\begin{equation}\label{eq-an-iteration}
	a_n=\sum_{i=1}^{n-1}a_i a_{n-i},\ n\geq 2.
\end{equation}
We are going to prove Claim 1 by induction on $n$. For $n=1$,
\begin{eqnarray*}
&& |\langle \wp,\wp, \hpc^{k_1},\hpc^l\rangle|= |\langle \hpc^{2N}, \hpc^{k_1},\hpc^i,\hpc^j\rangle|\\
&\leq & C^{1+D(2N+k_1+i+j,4)}=C^{1+D(2N+k_1+l,4)}< (2C)^{1+D(2N+k_1+l,4)}.
\end{eqnarray*}
Suppose that Claim 1 holds for $n\leq r-1$.
\begin{eqnarray*}
&&|\langle \wp,\wp, \hpc^{k_1},\dots,\hpc^{k_r},\hpc^l\rangle| \\
&\leq & |\langle \hpc^{2N}, \hpc^{k_1},\dots,\hpc^{k_r},\hpc^i,\hpc^j\rangle|\\
&&+\sum_{\begin{subarray}{c}S\subset \{1,\dots,r\}\\ \emptyset \neq S\neq\{1,\dots,r\} \end{subarray}}\sum_{\alpha=0}^{2N+1}\sum_{\beta=0}^{2N+1}|\langle \wp,\wp,\xi_S,\xi_\alpha\rangle g^{\alpha \beta}\langle \xi_\beta, \xi_{\{1,\dots,r\}-S},\hpc^i,\hpc^j\rangle|\\
&&+ \sum_{\begin{subarray}{c}S\subset \{1,\dots,r\}\\ \emptyset \neq S\neq\{1,\dots,r\} \end{subarray}}\sum_{\alpha=0}^{2N+1}\sum_{\beta=0}^{2N+1}
|\langle \wp,\xi_S, \hpc^i,\xi_\alpha\rangle g^{\alpha \beta}\langle \xi_{\beta}, \xi_{\{1,\dots,r\}-S},\wp,\hpc^j\rangle|\\
&=& |\langle \hpc^{2N}, \hpc^{k_1},\dots,\hpc^{k_r},\hpc^i,\hpc^j\rangle|\\
&& + \sum_{s=1}^r \sum_{\alpha=0}^{2N+1}\sum_{\beta=0}^{2N+1}|\langle \wp,\wp,\xi_{\{s\}},\xi_\alpha\rangle g^{\alpha \beta}\langle \xi_\beta, \xi_{\{1,\dots,r\}-\{s\}},\hpc^i,\hpc^j\rangle|\\
&&+\sum_{\begin{subarray}{c}S\subset \{1,\dots,r\}\\ \emptyset \neq S\neq\{1,\dots,r\}\\ |S|\geq 2 \end{subarray}}\sum_{\alpha=0}^{2N+1}\sum_{\beta=0}^{2N+1}|\langle \wp,\wp,\xi_S,\xi_\alpha\rangle g^{\alpha \beta}\langle \xi_\beta, \xi_{\{1,\dots,r\}-S},\hpc^i,\hpc^j\rangle|\\
&&+ \sum_{\begin{subarray}{c}S\subset \{1,\dots,r\}\\ \emptyset \neq S\neq\{1,\dots,r\} \end{subarray}}\sum_{\alpha=0}^{2N+1}\sum_{\beta=0}^{2N+1}
|\langle \wp,\xi_S, \hpc^i,\xi_\alpha\rangle g^{\alpha \beta}\langle \xi_{\beta}, \xi_{\{1,\dots,r\}-S},\wp,\hpc^j\rangle|\\
& \leq & r! C^{r+D(2N+\sum_{a=1}^r k_a+l,r+3)}
+r\cdot \frac{1}{2}\cdot(r-1)!C^{r+D(2N+\sum_{a=1}^r k_a+l,r+3)}\\
&&+\frac{1}{2}\sum_{k=2}^{r-1}\binom{r}{k} a_k\cdot k!\cdot a_{r-k}\cdot(r-k)! (2C)^{r+D(2N+\sum_{a=1}^r k_a+l,r+3)}\\
&&+\frac{1}{2}\sum_{k=1}^{r-1}\binom{r}{k} a_k\cdot k!\cdot a_{r-k}\cdot(r-k)! (2C)^{r+D(2N+\sum_{a=1}^r k_a+l,r+3)}\\
&\leq & \sum_{k=1}^{r-1}\binom{r}{k} a_k\cdot k!\cdot a_{r-k}\cdot(r-k)! (2C)^{r+D(2N+\sum_{a=1}^r k_a+l,r+3)}\\
&=& r! (2C)^{r+D(2N+\sum_{a=1}^r k_a+l,r+3)}\cdot \sum_{k=1}^{r-1}a_k a_{r-k}.
\end{eqnarray*}
In the above deduction we used (\ref{eq-claim1}) for smaller $n$, with $l$ possibly being $1$. This is legal because by (\ref{eq-vanish-2}) we can take another $k_i$ as $l$ if $l=1$. 
Now combining with (\ref{eq-an-iteration}) we obtain
\begin{eqnarray*}
|\langle \wp,\wp, \hpc^{k_1},\dots,\hpc^{k_r},\hpc^l\rangle|\leq a_r \cdot r! (2C)^{r+D(2N+\sum_{a=1}^r k_a+l,r+3)},
\end{eqnarray*}
which is Claim 1. Now we go to the general cases  $m\geq 1$. 

\textbf{Claim} 2: 
\begin{eqnarray}\label{eq-claim-2}
&&|\langle \underbrace{\wp,\dots,\wp}_{2m}, \hpc^{k_1},\dots,\hpc^{k_n},\hpc^l\rangle|\notag\\
&\leq& \frac{1}{4(2m-2+n)-2}\binom{2(2m-2+n)}{2m-2+n}(2m-2+n)!\\
&&\cdot 9^{m-1}(2C)^{2m-2+n+D(2mN+\sum_{a=1}^n k_a+l,2m+n+1)}.\notag
\end{eqnarray}
We are going to prove (\ref{eq-claim-2}) by induction on $m$. The case $m=1$ is Claim 1. Suppose that the Claim 2 holds for invariants with less $\wp$ or with  the same number of $\wp$ but smaller $n$, than that of $\langle \underbrace{\wp,\dots,\wp}_{2m}, \hpc^{k_1},\dots,\hpc^{k_n},\hpc^l\rangle$. Consider (\ref{eq-estimate-2-1}) and (\ref{eq-estimate-2-2}).

For $\emptyset \neq S\neq\{1,\dots,2m-2+n\}$, when $|S\cap \{1,\dots,2m-2\}|$ is even, say $2c$, the pair $(\alpha,\beta)\in \{0,1,\dots,2N+1\}^2$ such that
\[
\langle \wp,\wp,\xi_S,\xi_\alpha\rangle g^{\alpha \beta}\langle \xi_\beta, \xi_{\{1,\dots,2m-2+n\}-S},\hpc^i,\hpc^j\rangle\neq 0,
\]
lies in $\{0,1,\dots,2N\}^2$, so we have, when $c<m-1$, 
\begin{eqnarray*}
&&\sum_{\alpha=0}^{2N+1}\sum_{\beta=0}^{2N+1}|\langle \wp,\wp,\xi_S,\xi_\alpha\rangle g^{\alpha \beta}\langle \xi_\beta, \xi_{\{1,\dots,2m-2+n\}-S},\hpc^i,\hpc^j\rangle|\\
&\leq & a_{|S|}\cdot |S|! \cdot 9^{(c+1)-1}\cdot  \frac{1}{2}\cdot a_{2m-2+n-|S|}\cdot(2m-2+n-|S|)!\cdot 9^{(m-1-c)-1}\\
&& \cdot (2C)^{2m-2+n+D(2mN+\sum_{a=1}^n k_a+l,2m+n+1)},
\end{eqnarray*}
and when $c=m-1$,
\begin{eqnarray*}
&&\sum_{\alpha=0}^{2N+1}\sum_{\beta=0}^{2N+1}|\langle \wp,\wp,\xi_S,\xi_\alpha\rangle g^{\alpha \beta}\langle \xi_\beta, \xi_{\{1,\dots,2m-2+n\}-S},\hpc^i,\hpc^j\rangle|\\
&\leq & a_{|S|}\cdot |S|! \cdot 9^{m-1}\cdot  \frac{1}{2}\cdot a_{2m-2+n-|S|}\cdot(2m-2+n-|S|)!\\
&& \cdot (2C)^{2m-2+n+D(2mN+\sum_{a=1}^n k_a+l,2m+n+1)}.
\end{eqnarray*}
On the other hand when $|S\cap \{1,\dots,2m-2\}|$ is odd, say $2c+1$,  the pair $(\alpha,\beta)\in \{0,1,\dots,2N+1\}^2$ such that
\[
\langle \wp,\wp,\xi_S,\xi_\alpha\rangle g^{\alpha \beta}\langle \xi_\beta, \xi_{\{1,\dots,2m-2+n\}-S},\hpc^i,\hpc^j\rangle\neq 0,
\]
is $(2N+1,2N+1)$, so  we have
\begin{eqnarray*}
&&\sum_{\alpha=0}^{2N+1}\sum_{\beta=0}^{2N+1}|\langle \wp,\wp,\xi_S,\xi_\alpha\rangle g^{\alpha \beta}\langle \xi_\beta, \xi_{\{1,\dots,2m-2+n\}-S},\hpc^i,\hpc^j\rangle|\\
&\leq & a_{|S|}\cdot |S|! \cdot 9^{(c+2)-1}\cdot \frac{1}{2}\cdot a_{2m-2+n-|S|}\cdot(2m-2+n-|S|)!\cdot 9^{(m-1-c)-1}\\
&& \cdot (2C)^{2m-2+n+D(2mN+\sum_{a=1}^n k_a+l,2m+n+1)}.
\end{eqnarray*}
In any case, we have
\begin{eqnarray*}
&&\sum_{\alpha=0}^{2N+1}\sum_{\beta=0}^{2N+1}|\langle \wp,\wp,\xi_S,\xi_\alpha\rangle g^{\alpha \beta}\langle \xi_\beta, \xi_{\{1,\dots,2m-2+n\}-S},\hpc^i,\hpc^j\rangle|\\
&\leq &  \frac{1}{2}  a_{|S|}\cdot |S|! \cdot a_{2m-2+n-|S|}\cdot(2m-2+n-|S|)!\cdot 9^{m-1}\\
&& \cdot (2C)^{2m-2+n+D(2mN+\sum_{a=1}^n k_a+l,2m+n+1)}.
\end{eqnarray*}
The same argument shows, for $\emptyset \neq S\neq\{1,\dots,2m-2+n\}$,
\begin{eqnarray*}
&&\sum_{\alpha=0}^{2N+1}\sum_{\beta=0}^{2N+1}
|\langle \wp,\xi_S, \hpc^i,\xi_\alpha\rangle g^{\alpha \beta}\langle \xi_{\beta}, \xi_{\{1,\dots,2m-2+n\}-S},\wp,\hpc^j\rangle|\\
&\leq &  \frac{1}{2}  a_{|S|}\cdot |S|! \cdot a_{2m-2+n-|S|}\cdot(2m-2+n-|S|)!\cdot 9^{m-1}\\
&& \cdot (2C)^{2m-2+n+D(2mN+\sum_{a=1}^n k_a+l,2m+n+1)}.
\end{eqnarray*}
We need to treat (\ref{eq-estimate-2-1}) in the case $|S|=1$ individually. In this case, 
when $|S\cap\{1,\dots,2m-2\}|=1$, due to (\ref{eq-vanish-1}) we have 
\[
\langle \wp,\wp,\xi_S,\xi_\alpha\rangle=\langle \wp,\wp,\wp,\xi_\alpha\rangle=0.
\]
So for $|S|=1$, only the $S$ such that $|S|\cap \{1,\dots,2m-2\}=\emptyset$ contributes to (\ref{eq-estimate-2-1}). Thus if $|S|=1$,
\begin{eqnarray*}
&&\sum_{\alpha=0}^{2N+1}\sum_{\beta=0}^{2N+1}|\langle \wp,\wp,\xi_S,\xi_\alpha\rangle g^{\alpha \beta}\langle \xi_\beta, \xi_{\{1,\dots,2m-2+n\}-S},\hpc^i,\hpc^j\rangle|\\
&\leq & \frac{1}{2} a_{a}\cdot  \cdot a_{2m-3+n}\cdot(2m-3+n)!\cdot 9^{m-2}\\
&& \cdot (2C)^{2m-2+n+D(2mN+\sum_{a=1}^n k_a+l,2m+n+1)}.
\end{eqnarray*}
Hence
\begin{eqnarray*}
&&|\langle \underbrace{\wp,\dots,\wp}_{2m}, \hpc^{k_1},\dots,\hpc^{k_n},\hpc^l\rangle| \\
&\leq & |\langle \xi_{2N},\underbrace{\wp,\dots,\wp}_{2m-2}, \hpc^{k_1},\dots,\hpc^{k_n},\hpc^i,\hpc^j\rangle|\\
&&+ \sum_{\begin{subarray}{c}S\subset \{1,\dots,2m-2+n\}\\ |S|=1 \end{subarray}}\sum_{\alpha=0}^{2N+1}\sum_{\beta=0}^{2N+1}|\langle \wp,\wp,\xi_S,\xi_\alpha\rangle g^{\alpha \beta}\langle \xi_\beta, \xi_{\{1,\dots,2m-2+n\}-S},\hpc^i,\hpc^j\rangle|\\
&&+\sum_{\begin{subarray}{c}S\subsetneqq \{1,\dots,2m-2+n\}\\ |S|\geq 2  \end{subarray}}\sum_{\alpha=0}^{2N+1}\sum_{\beta=0}^{2N+1}|\langle \wp,\wp,\xi_S,\xi_\alpha\rangle g^{\alpha \beta}\langle \xi_\beta, \xi_{\{1,\dots,2m-2+n\}-S},\hpc^i,\hpc^j\rangle|\\
&&+ \sum_{\begin{subarray}{c}S\subset \{1,\dots,2m-2+n\}\\ \emptyset \neq S\neq\{1,\dots,2m-2+n\} \end{subarray}}\sum_{\alpha=0}^{2N+1}\sum_{\beta=0}^{2N+1}
|\langle \wp,\xi_S, \hpc^i,\xi_\alpha\rangle g^{\alpha \beta}\langle \xi_{\beta}, \xi_{\{1,\dots,2m-2+n\}-S},\wp,\hpc^j\rangle|\\
&\leq & a_{2m-2+n}\cdot (2m-2+n)! 9^{m-2}(2C)^{2m-2+n+D(2mN+\sum_{a=1}^n k_a+l,2m+n+1)}\\
&&+\frac{2m-2+n}{2}\cdot a_{2m-3+n}\cdot(2m-3+n)! 9^{m-2}(2C)^{2m-2+n+D(2mN+\sum_{a=1}^n k_a+l,2m+n+1)} \\
&&+\frac{1}{2}\sum_{k=2}^{2m-3+n}\bigg[\binom{2m-2+n}{k} a_{k}\cdot k!\cdot a_{2m-2+n-k}\cdot(2m-2+n-k)!\cdot 9^{m-1}\\
&&\cdot (2C)^{2m-2+n+D(2mN+\sum_{a=1}^n k_a+l,2m+n+1)}\bigg]\\
&&+\frac{1}{2}\sum_{k=1}^{2m-3+n}\binom{2m-2+n}{k}  a_{k}\cdot k!\cdot a_{2m-2+n-k}\cdot (2m-2+n-k)!\cdot 9^{m-1}\\
&&\cdot (2C)^{2m-2+n+D(2mN+\sum_{a=1}^n k_a+l,2m+n+1)}\\
&\leq &\frac{9}{2}\cdot a_{2m-3+n}\cdot(2m-2+n)! 9^{m-2}(2C)^{2m-2+n+D(2mN+\sum_{a=1}^n k_a+l,2m+n+1)} \\
&&+\frac{1}{2}\sum_{k=2}^{2m-3+n}\bigg[\binom{2m-2+n}{k} a_{k}\cdot k!\cdot a_{2m-2+n-k}\cdot(2m-2+n-k)!\cdot 9^{m-1}\\
&&\cdot (2C)^{2m-2+n+D(2mN+\sum_{a=1}^n k_a+l,2m+n+1)}\bigg]\\
&&+\frac{1}{2}\sum_{k=1}^{2m-3+n}\binom{2m-2+n}{k}  a_{k}\cdot k!\cdot a_{2m-2+n-k}\cdot (2m-2+n-k)!\cdot 9^{m-1}\\
&&\cdot (2C)^{2m-2+n+D(2mN+\sum_{a=1}^n k_a+l,2m+n+1)}\\
&= & (2m-2+n)!9^{m-1}(2C)^{2m-2+n+D(2mN+\sum_{a=1}^n k_a+l,2m+n+1)}\sum_{k=1}^{2m-3+n}a_{k}a_{2m-2+n-k}\\
&\leq & (2m-2+n)!9^{m-1}(2C)^{2m-2+n+D(2mN+\sum_{a=1}^n k_a+l,2m+n+1)} \cdot a_{2m-2+n},
\end{eqnarray*}
where in the third inequality we have used
\[
a_{2m-2+n}\leq 4 a_{2m-3+n}.
\]
Thus we have proved (\ref{eq-claim-2}). So we have 
\begin{equation}\label{eq-estimate-3}
|\langle \gamma_1,\dots,\gamma_n\rangle|\leq  \frac{1}{4(n-3)-2}\frac{(2n-6)!}{(n-3)!} (18C)^{n-3+D(\sum_{i=1}^n \mathrm{deg}_\bC(\gamma_i),n)}
\end{equation}
for $\gamma_i\in\mathbf{B}$, $1\leq i\leq n$, if there is at least one $\wp$ in $\gamma_1,\dots,\gamma_n$. Since
\[
\frac{(2n-6)!}{(n-3)!}
%\leq \frac{(2n-6)(2n-6)\cdot (2n-4)(2n-4)\dots 2\cdot 2}{(n-3)!}
\leq 4^{n-3} (n-3)!,
\]
combining (\ref{eq-estimate-3}) with Zinger's estimate (\ref{eq-estimate-0}), enlarging $C$ as the beginning of the proof, we obtain the conclusion.
\end{proof}

\begin{prop}\label{prop-evenconvergent}
	Let $\qdc$ be an even dimensional quadric hypersurface. Then the quantum cohomology of $\qdc$ is analytic in a neighborhood of $\mathbf0\in H^*(\qdc)$.
\end{prop}
\begin{proof}
This follows from theorem \ref{thm-estimate-primitive}, by a similar computation as in the proof of Proposition \ref{prop-oddconvergent}.
\end{proof}

\section{Proof of Gamma conjecture II for quadrics}\label{sec:proof-GammaII}

The main result of this section is the following theorem.
\begin{thm}\label{thm-GammaII-quadric}
A quadric hypersurface 
	$\qdc$ satisfies the Gamma conjecture II.
\end{thm}

This theorem follows from Proposition \ref{evencase} and \ref{oddcase}. 
%Note that a one-dimensional quadric is isomorphic to $\bP^1$, and a two-dimensional quadric is isomorphic to $\bP^1\times\bP^1$. 
The Gamma conjecture II for $\bP^1$ and $\bP^1\times\bP^1$ was already known \cite{CDG18,FaZh19,GGI16}. So we restrict ourselves to the cases of dimension$\geq3$.

\subsection{\texorpdfstring{$J$}{J}-function and \texorpdfstring{$P$}{P}-function}\label{Pfunction}

Recall that $\cZ^{coh}_\mathbf{0}:H^*(F)\to\cS_\mathbf{0}$ is a linear isomorphism. Following \cite[Definition 3.6.3]{GGI16}, Givental's $J$-function of $F$ is given by
\ban
J_F(t):=z^{\frac{\dim F}2}(\cZ^{coh}_\mathbf{0})^{-1}\mathbbm1=z^{\frac{\dim F}2}z^{-\rho}z^\mu L(\mathbf{0},z)^{-1}\mathbbm1,\textrm{ via }t=z^{-1}.
\nan
More explicitly, we have \cite[(3.6.7)]{GGI16}
\ban
J_F(t)=t^\rho\bigg[\mathbbm1+\suml_{i=0}^{s-1}T^i\suml_{\substack{\mathbf d\in\Eff(F)\\m\geq0}}\<\tau_m(T_i)\>_{\mathbf d}^Ft^{\int_{\mathbf d}c_1(F)}\bigg],
\nan
where $\{T^i\}$ is the basis in $H^*(F)$ satisfying $\int_FT^i\cup T_j=\delta_{ij}$.

For $F=\qdc_N$, a smooth quadric hypersurface in $\bP^{N+1}$, it follows from \cite[Theorem 9.1]{Giv96} that
\ban
J_{\qdc_N}(t)=t^{N\hpc}\suml_{d=0}^\infty\frac{\prodl_{m=1}^{2d}(2\hpc+m)}{\prodl_{m=1}^d(\hpc+m)^{N+2}}t^{Nd}.
\nan

\begin{remark}
Recall that the original Givental's big $J$-function \cite[(7)]{Giv04} is a formal function of $\mathbf t\in H^*(F)$ taking values in $H^*(F)$, given by
\ban
J_F(\mathbf t,\hbar;\qtp)=\hbar+\mathbf t+\suml_{i=0}^{s-1}T^i\suml_{\substack{\mathbf{d}\in\Eff(F)\\m,n\geq0}}{\frac{\qtp^\mathbf d}{n!\hbar^{m+1}}}\<\tau_m(T_i)\tau_0(\mathbf t)^n\>^F_{\mathbf d}.
\nan
It can be reduced to the form 
\ban
J_F(\mathbf t,\hbar;\qtp)=\hbar e^{\frac{\mathbf t}{\hbar}}\Big(\mathbbm1+\suml_{i=0}^{s-1}T^i\suml_{\substack{\mathbf{d}\in \Eff(F)\\m\geq0}}\frac{\qtp^\mathbf de^{\int_\mathbf d\tau}}{\hbar^{m+2}}\<\tau_m(T_i)\>^F_\mathbf d\Big),
\nan
for $\mathbf t\in H^0(F)\oplus H^2(F)$, after applying the fundamental class axiom and the divisor axiom.
Givental's $J$-function in this paper is obtained from the original one by setting
\ban
J_F(t)&=&J_F(\mathbf t=c_1(F)\log t,\hbar=1;\qtp=\mathbf{1})\\
&=&e^{c_1(X)\log t}\bigg(\mathbbm1+\suml_{i=0}^{s-1}T^i\suml_{\substack{\mathbf d\in\Eff(F)\\m\geq0}}\<\tau_m(T_i)\>_{\mathbf d}^Ft^{\int_{\mathbf{d}}c_1(F)}\bigg).
\nan
\end{remark}

In this note, instead of the well-known $J_F$, we will need a new function $P_F$ to understand $\cZ^{coh}_\mathbf 0$ (see Lemma \ref{f0}). The new function $P_F$ is defined as
\ba\label{defofP}
P_F(t):=z^{\frac{\dim F}2}(\cZ^{coh}_\mathbf{0})^*\mathbbm1=z^{\frac{\dim F}2}z^{\rho}z^{\mu}L(\mathbf{0},z)^*\mathbbm1,\textrm{ via }t=z^{-1}.
\na
Here for any $M\in\End(H^*(F))$, we let $M^*\in\End(H^*(F))$ be its dual with respect to the Poincar\'e pairing, i.e. $\<M\alpha,\beta\>^F=\<\alpha,M^*\beta\>^F$. For example, we have $(z^\rho)^*=z^\rho$, and $(z^\mu)^*=z^{-\mu}$. Since 
\ban
\<L(\mathbf{0},z)^*\alpha,\beta\>^F=\<\alpha,L(\mathbf{0},z)\beta\>^F=\<\alpha,\beta\>^F+\suml_{\substack{\mathbf d\in\Eff(F)\\m\geq0}}(\frac{-1}{z})^{m+1}\<\tau_m(\beta)\tau_0(\alpha)\>^F_{\mathbf d},
\nan
it follows that
\ban
L(\mathbf{0},z)^{*}\mathbbm1=\suml_{i=0}^{s-1}\<L(\mathbf{0},z)^{*}\mathbbm1,T_i\>^FT^i=\mathbbm1+\suml_{i=0}^{s-1}T^i\suml_{\substack{\mathbf d\in\Eff(F)\\m\geq0}}(\frac{-1}{z})^{m+1}\<\tau_m(T_i)\tau_0(\mathbbm1)\>^F_{\mathbf d}.
\nan
So by the definitions of $\rho$ and $\mu$, we have
\ban
P _F(t)&=&z^{\frac{\dim F}2}z^{\rho}z^\mu L(\mathbf{0},z)^{*}\mathbbm1\\
&=&t^{-\rho}\bigg[\mathbbm1+\suml_{i=0}^{s-1}T^i\suml_{\substack{\mathbf d\in\Eff(F)\\m\geq0}}(-t)^{m+1}\<\tau_m(T_i)\tau_0(\mathbbm1)\>^F_{\mathbf d}t^{-\deg T^i}\bigg]\\
&=&t^{-\rho}\bigg[\mathbbm1+\suml_{i=0}^{s-1}(-1)^{\deg T^i}T^i\suml_{\substack{\mathbf d\in\Eff(F)\\m\geq0}}\<\tau_m(T_i)\tau_0(\mathbbm1)\>^F_{\mathbf d}(-t)^{\int_dc_1(F)}\bigg].
\nan
Here we have used the dimension axiom in the last equality.

One can observe that $P_F(t)$ and $J_F(t)$ are closely related. In fact, we can determine $P_F(t)$ from $J_F(t)$. For $F=\qdc_N$, we have
\ba
P_{\qdc_N}(t)&=&t^{-N\hpc}\suml_{d=0}^\infty\frac{\prodl_{m=1}^{2d}(-2\hpc+m)}{\prodl_{m=1}^d(-\hpc+m)^{N+2}}(-t)^{Nd}=t^{-N\hpc}\suml_{d=0}^\infty\frac{\prodl_{m=1}^{2d}(2\hpc-m)}{\prodl_{m=1}^d(\hpc-m)^{N+2}}t^{Nd}\nonumber\\
&=&\suml_{d=0}^\infty\Big[\frac{\Gamma(\hpc-d)}{\Gamma(\hpc)}\Big]^{N+2}\frac{\Gamma(2\hpc)}{\Gamma(2\hpc-2d)}t^{N(d-\hpc)}.\label{PQ}
\na

\subsection{Even-dimensional case}

In this subsection, let $\qdc=\qdc_{2N}$ be a smooth quadric hypersurface in $\bP^{2N+1}$ with $N\geq2$. The main result of this subsection is Proposition \ref{evencase}. 

\subsubsection{Singular cohomology}
Recall that $\hpc$ stands for the hyperplane class. 
%Denote by $\hpc$ the pullback of the hyperplane class of $\bP^{2N+1}$ on $H^2(Q)$. Then 
We have $c_1(\qdc)=2N\hpc$, and a basis of $H^*(\qdc)$ is 
\ban
\hpc^0,\dots,\hpc^{2N},\wp=e_{2k}-e'_{2k},
\nan
such that $\wp\in H^{2N}(\qdc)$ is a primitive class with $\wp\cup\wp=\hpc^{2N}$. The Poincar\'e pairing on $H^*(\qdc)$ is
\[
\int_\qdc\hpc^i\cup \hpc^j=2\delta_{i+j,2N},\
\int_\qdc\hpc^i\cup\wp=0,\
\int_\qdc\wp\cup\wp=2.
\]
The Hodge grading operator $\mu\in\End(H^*(\qdc))$ is given by
\ban
\mu(\hpc^i)=(i-N)\hpc^i,\quad\mu(\wp)=0.
\nan

\subsubsection{Specialized quantum product}
Let $\circ:=\bullet_{\mathbf t=\mathbf0}$ be the quantum product on $H^*(\qdc)$ specialized at $\mathbf t=\mathbf0$. Then the product $\circ$ is determined by $3$-point invariants (see Lemma \ref{nonzerothreepointinvariants}). For example, we have
\ban
c_1(\qdc)\circ \hpc^i&=&\left\{\begin{array}{cc}
2N\hpc^{i+1},&0\leq i\leq 2N-2,\\
2N\hpc^{2N}+4N\hpc^0,&i=2N-1,\\
4N\hpc,&i=2N,
\end{array}\right.\\
c_1(\qdc)\circ\wp&=&0.
\nan
In particular, the characteristic polynomial of $(c_1(\qdc)\circ)\in\End(H^*(\qdc))$ is
\ban
x^{2N+2}-4(2N)^{2N}x^2.
\nan
Then $T:=2N\cdot4^{\frac1{2N}}$ is the spectral radius of $(c_1(\qdc)\circ)$.
 
\if{
\begin{lemma}\label{lem-Beauville}
We have $\hpc^{2N}\circ \hpc^{2N}=4\hpc^0$ and 
\[
\hpc^i\circ\wp=\left\{\begin{array}{cc}
\wp,&i=0,\\
0,&1\leq i\leq 2N-1,\\
-2\wp,&i=2N.
\end{array}\right.
\]
\end{lemma}
\begin{proof}
For the first equality, we have
\ban
\hpc^{2N}\circ \hpc^{2N}&=&(\hpc\circ \hpc^{2N-1}-2\hpc^0)\circ \hpc^{2N}\\
&=&\hpc^{2N-1}\circ(\hpc\circ \hpc^{2N})-2\hpc^{2N}\\
&=&\hpc^{2N-1}\circ2\hpc-2\hpc^{2N}\\
&=&2(\hpc^{2N}+2\hpc^0)-2\hpc^{2N}\\
&=&4\hpc^0.
\nan
For the second equality, for $1\leq i\leq 2N-1$, we have
\[\hpc^i\circ\wp=\hpc^{i-1}\circ \hpc\circ\wp=0,\]
and for $i=2N$, we have
\ban
\hpc^{2N}\circ\wp&=&(\hpc\circ \hpc^{2N-1}-2\hpc^0)\circ\wp=\hpc^{2N-1}\circ(\hpc\circ\wp)-2\wp=-2\wp.
\nan
\end{proof}
}\fi

Let $\zeta:=\exp\big(\frac{2\pi\mathbf{i}}{2N}\big)$, and set
\ban
v_\pm&=&2\mathbf{i}\big(\frac{\hpc^0}4-\frac{\hpc^{2N}}8\big)\pm\frac\wp2,\\
v_k&=&\frac{(-1)^k}{\sqrt N}\bigg(\frac{\hpc^0}2+\suml_{p=1}^{2N}\big(\frac{T\zeta^{-k}}{2N}\big)^{-p}\hpc^p\bigg),\quad k=0,1,\dots,2N-1.
\nan

\begin{cor}\label{semisimpleofQ}
The Frobenius algebra $(H^*(\qdc),\circ)$ is semisimple, and $v_+,v_-,v_0,\dots,v_{2N-1}$ form a normalized idempotent basis.
\end{cor}
\begin{proof}
By direct calculation, we can verify that
\begin{eqnarray*}
&c_1(\qdc)\circ v_\pm=0,\ c_1(\qdc)\circ v_k=T\zeta^{-k}v_k\ \mbox{for}\ 0\leq k\leq 2N-1,\\
&v_\pm\circ v_\pm=2\mathbf{i} v_\pm,\ v_\pm\circ v_\mp=0.
\end{eqnarray*}
So 
%we see that
\ban
0=(c_1(\qdc)\circ v_\pm)\circ v_k=v_\pm\circ(c_1(\qdc)\circ v_k)=T\zeta^{-k}v_\pm\circ v_k\Rightarrow v_\pm\circ v_k=0,\nan
and for $k\neq k'$,
\ban
T\zeta^{-k}v_k\circ v_{k'}=(c_1(\qdc)\circ v_k)\circ v_{k'}=v_k\circ(c_1(\qdc)\circ v_{k'})=T\zeta^{-k'}v_k\circ v_{k'}\Rightarrow v_k\circ v_{k'}=0.
\nan
We write 
\ban
\mathbbm1=a_+v_++a_-v_-+a_0v_0+\cdots+a_{2N-1}v_{2N-1}.
\nan
for some constants $a_{\pm},a_k(0\leq k\leq2N-1)$. Then 
\ban
v_k=v_k\circ\mathbbm1=v_k\circ a_kv_k\Rightarrow a_k\neq0\textrm{ and }v_k\circ v_k=\frac1{a_k}v_k.
\nan
So $\frac{v_+}{2\mathbf{i}},\frac{v_-}{2\mathbf{i}},a_0v_0,\dots,a_{2N-1}v_{2N-1}$ form an idempotent basis of $(H^*(\qdc),\circ)$. Moreover, we can check that
\[\<v_\pm,v_\pm\>^\qdc=\<v_k,v_k\>^\qdc=1,\]
which implies that $v_+,v_-,v_0,\dots,v_{2N-1}$ form a normalized idempotent basis.
\end{proof}

\subsubsection{Determination of flat sections}
Recall that 
\ban
(\nabla_{\mathbf 0})_{z\partial_z}=z\partial_z-\frac1z(c_1(\qdc)\circ)+\mu.
\nan
So for a $H^*(\qdc)$-valued function $f(z)=\suml_{i=0}^{2N}f_i(z)z^{i-N}\hpc^{2N-i}+f_\wp(z)\wp$, we have
\ba
(\nabla_{\mathbf 0})_{z\partial_z}f&=&z^{-N}\bigg\{\suml_{i=0}^{2N-2}z^i\Big[z\partial_zf_i-2Nf_{i+1}\Big]\hpc^{2N-i}\nonumber\\
&&\quad+z^{2N-1}\Big[z\partial_zf_{2N-1}-2N(f_{2N}+2f_0\cdot z^{-2N})\Big]\hpc^{2N-(2N-1)}\nonumber\\
&&\quad+z^{2N}\Big[z\partial_zf_{2N}-2N\cdot 2f_1\cdot z^{-2N}\Big]\hpc^{2N-(2N)}\bigg\}+z\partial_zf_\wp\wp.\label{nabla}
\na

\begin{lemma}\label{fiintermsoff0}
Let $f(z)=\suml_{i=0}^{2N}f_i(z)z^{i-N}\hpc^{2N-i}+f_\wp(z)\wp\in\cS_\mathbf 0$. Then $f_\wp(z)$ is a constant function, and $f(z)-f_\wp(z)\wp$  is determined by $f_0(z)$ via
\ban
f_i&=&(2N)^{-i}(z\partial_z)^if_0,\quad0\leq i\leq2N-1,\\
f_{2N}&=&(2N)^{-2N}(z\partial_z)^{2N}f_0-2f_0\cdot z^{-2N},
\nan
and 
\ba\label{equationinz}
\bigg\{(z\partial_z)^{2N+1}-(2N)^{2N}z^{-2N}\Big[4z\partial_z-2\cdot2N\Big]\bigg\}f_0=0.
\na
\end{lemma}
\begin{proof}
Since $(\nabla_{\mathbf 0})_{z\partial_z}f=0$, it follows from \eqref{nabla} that $\partial_zf_\wp=0$ and
\ban
\frac1{2N}z\partial_zf_i&=&f_{i+1},\quad0\leq i\leq2N-2,\\
\frac1{2N}z\partial_zf_{2N-1}&=&f_{2N}+2f_0\cdot z^{-2N},\\
\frac1{2N}z\partial_zf_{2N}&=&2f_1\cdot z^{-2N},
\nan
which gives the required equations.
\end{proof}

\if{
For $t:=z^{-1}$ and $D:=t\partial_t$, we have $D=-z\partial_z$ and equation \eqref{equationinz} is transformed to
\ba\label{equationint}
\{D^{2N+1}-(2N)^{2N}t^{2N}[4D+2\cdot2N]\}f_0=0.
\na
}\fi

\begin{lemma}\label{f0}
For $\gamma\in H^*(\qdc)$, we have
\ban
\cZ^{coh}_{\mathbf 0}(\gamma)=\begin{pmatrix}\hpc^0 & \hpc^1 & \dots & \hpc^{2N}& \wp\end{pmatrix}\cdot\begin{pmatrix}
z^{N}\Big(\frac{z\partial_z}{2N}\Big)^{2N}f_0-z^{-N}\cdot 2f_0\\
z^{N-1}\Big(\frac{z\partial_z}{2N}\Big)^{2N-1}f_0\\
\vdots\\
z^{-(N-1)}\Big(\frac{z\partial_z}{2N}\Big)f_0\\
z^{-N}f_0\\
f_\wp
\end{pmatrix},
\nan
with 
\ban
f_0(z)=\frac12\int_\qdc\gamma\cup P_\qdc(z^{-1})\textrm{ and }f_\wp=\frac12\int_\qdc\gamma\cup\wp.
\nan
Recall that $P_\qdc$ is defined in \eqref{defofP}.
\end{lemma}
\begin{proof}
From Lemma \ref{fiintermsoff0}, we only need to determine $f_0$ and $f_\wp$. For $f_0$, we have
\ban
z^{-N}f_0(z)=\frac12\<\cZ^{coh}_{\mathbf 0}(\gamma),\mathbbm1\>^\qdc=\frac12\<\gamma,(\cZ^{coh}_{\mathbf 0})^*\mathbbm1\>^\qdc=\frac12\<\gamma,z^{-N}P_\qdc(z^{-1})\>^\qdc,
\nan
where we have used $\int_\qdc \hpc^{2N+1}=2$ in the first equality. For $f_\wp$, we have
\ban
f_\wp&=&\frac12\<\cZ^{coh}_{\mathbf 0}(\gamma),\wp\>^\qdc\\
&=&\frac12\<\gamma,z^\rho z^\mu L(\mathbf{0},z)^*\wp\>^\qdc\\
&=&\frac12\<\gamma,z^\rho z^\mu\bigg(\wp+\suml_{i=0}^{2N}\suml_{\substack{d>0\\m\geq0}}(\frac{-1}z)^{m+1}\<\tau_m(\hpc^i)\tau_0(\wp)\>_d^\qdc\frac{\hpc^{2N-i}}{2}\\
&&+\suml_{\substack{d>0\\m\geq0}}(\frac{-1}z)^{m+1}\<\tau_m(\wp)\tau_0(\wp)\>_d^\qdc\frac{\wp}{2}\bigg)\>^\qdc\\
&=&\frac12\<\gamma,z^\rho z^\mu\bigg(\wp+\suml_{\substack{d>0\\m\geq0}}(\frac{-1}z)^{m+1}\<\tau_m(\wp)\tau_0(\wp)\>_d^\qdc\frac{\wp}{2}\bigg)\>^\qdc\\
&=&\frac12\<\gamma,\wp+\suml_{\substack{d>0\\m\geq0}}(\frac{-1}z)^{m+1}\<\tau_m(\wp)\tau_0(\wp)\>_d^\qdc\frac{\wp}{2}\>^\qdc,
\nan
where we have used $\int_\qdc\wp\cup\wp=2$ in the first equality, and the fourth equality comes from \cite[Lemma 1]{LePa04}. Note that $f_\wp$ is independent of $z$, which implies that
\ban
f_\wp=\frac12\int_\qdc\gamma\cup\wp.
\nan
\end{proof}

\subsubsection{Explicit computations}

Let $S_+$ and $S_-$ be the spinor bundles on $\qdc$ with
\ban
\Ch(S_\pm)=\frac{2^N}{1+e^{2\pi\mathbf{i} \hpc}}\pm(2\pi)^N\frac\wp2\textrm{ (see Theorem \ref{computationofch}).}
\nan
It is known that 
\ban
\mathcal E=(S_+,S_-,\cO,\cO(1),\dots,\cO(2N-1)),\textrm{ and }\mathcal E'=(S_-,S_+,\cO,\dots,\cO(2N-1)),
\nan
are full exceptional collections on $\mathcal D^b(\qdc)$ \cite{Kap88}. We have
\ban
\hat\Gamma_\qdc\Ch(S_\pm)&=&\frac{\Gamma(1+\hpc)^{2N+2}}{\Gamma(1+2\hpc)}\cdot\frac{2^{N}}{1+e^{2\pi\mathbf{i} \hpc}}\pm(2\pi)^N\frac\wp2,\\
\hat\Gamma_\qdc\Ch(\cO(k))&=&\frac{\Gamma(1+\hpc)^{2N+2}}{\Gamma(1+2\hpc)}\cdot e^{2k\pi\mathbf{i} \hpc},\quad0\leq k\leq 2N-1.
\nan
To study $\cZ^{K}_{\mathbf 0}(S_\pm)$ and $\cZ^{K}_{\mathbf 0}(\cO(k))$, we define
\ban
g_\pm(z)&=&\int_\qdc\hat\Gamma_\qdc\Ch(S_\pm)\cup P_\qdc(z^{-1}),\\
g_k(z)&=&\int_\qdc\hat\Gamma_\qdc\Ch(\cO(k))\cup P_\qdc(z^{-1}),\quad0\leq k\leq 2N-1.
\nan
Then from Lemma \ref{f0}, we have
\ban
\<\cZ^{K}_{\mathbf 0}(S_\pm)(z),\mathbbm1\>^\qdc&=&\frac{g_\pm(z)}{(2\pi z)^N},\\
\<\cZ^{K}_{\mathbf 0}(\cO(k))(z),\mathbbm1\>^\qdc&=&\frac{g_k(z)}{(2\pi z)^N},\quad0\leq k\leq 2N-1.
\nan

In this subsection, we study the asymptotic behavior of $g_\pm$ and $\cZ^{K}_{\mathbf 0}(S_\pm)$ in Lemma \ref{asymptoticofg} and \ref{asymptoticofspm}, and that of $g_k$ and $\cZ^{K}_{\mathbf 0}(\cO(k))$ in Lemma \ref{asymptoticofgk} and \ref{asymptoticofsk}. We will prove Gamma conjecture II for even-dimensional quadrics in Proposition \ref{evencase}.

\begin{lemma}\label{asymptoticofg}
	As $z\to0$ in any closed subsector of $|\arg z-\frac\pi{2N}|<\frac\pi2+\frac\pi{2N}$, we have
	\ban
	\frac{g_\pm(z)}{(2\pi z)^N}\sim-\frac{\mathbf{i}}{2}+\suml_{m=1}^\infty a_mz^m,\quad a_m\in\bC.
	\nan
\end{lemma}
\begin{proof}
Note that
\ban
\int_\qdc\wp\cup P_{\qdc}(z^{-1})=0.
\nan
Then from \eqref{PQ}, we have
\ban
g_\pm(z)&=&\suml_{d=0}^\infty\int_\qdc\frac{\Gamma(1+\hpc)^{2N+2}}{\Gamma(1+2\hpc)}\cdot\frac{2^{N}}{1+e^{2\pi\mathbf{i} \hpc}}\Big[\frac{\Gamma(\hpc-d)}{\Gamma(\hpc)}\Big]^{2N+2}\frac{\Gamma(2\hpc)}{\Gamma(2\hpc-2d)}z^{2N(\hpc-d)}\\
&=&\suml_{d=0}^\infty\textrm{Coeff}_{s^{2N}}\bigg(2\frac{\Gamma(1+s)^{2N+2}}{\Gamma(1+2s)}\cdot\frac{2^{N}}{1+e^{2\pi\mathbf{i} s}}\Big[\frac{\Gamma(s-d)}{\Gamma(s)}\Big]^{2N+2}\frac{\Gamma(2s)}{\Gamma(2s-2d)}z^{2N(s-d)}\bigg),
\nan
where in the second equality, we have used $\int_\qdc\hpc^{2N}=2$. From
\begin{equation}\label{eq-Gammafunction-identity-1}
\frac{\Gamma(2s)}{\Gamma(1+2s)}=\frac1{2s},\quad\frac{\Gamma(1+s)}{\Gamma(s)}=s,
\end{equation}
we get
\ban
g_\pm(z)&=&2^{N}\suml_{d=0}^\infty\textrm{Coeff}_{s^{2N}}\bigg(\frac{s^{2N+1}}{1+e^{2\pi\mathbf{i} s}}\cdot\frac{\Gamma(s-d)^{2N+2}}{\Gamma(2s-2d)}z^{2N(s-d)}\bigg).
\nan
Moreover, from the equality $\Gamma(x)\Gamma(1-x)=\frac{\pi}{\sin\pi x}$, we have
\ban
\frac1{1+e^{2\pi\mathbf{i} s}}&=&\frac1{1+e^{2\pi\mathbf{i}(s-d)}}\\
&=&\frac{e^{-\pi\mathbf{i}(s-d)}}{e^{-\pi\mathbf{i}(s-d)}+e^{\pi\mathbf{i}(s-d)}}\\
&=&\frac{e^{-\pi\mathbf{i}(s-d)}}{2\cos\pi(s-d)}\\
&=&\frac{e^{-\pi\mathbf{i}(s-d)}}{2\sin\pi\big(\frac12-(s-d)\big)}\\
&=&\frac{e^{-\pi\mathbf{i}(s-d)}}{2\pi}\Gamma\big(\frac12-(s-d)\big)\Gamma\big(\frac12+(s-d)\big),
\nan
and from the equality $\Gamma(x)\Gamma(x+\frac12)=\frac{2\sqrt\pi}{4^{x}}\Gamma(2x)$, we have
\ban
\frac1{\Gamma(2s-2d)}=\frac{2\sqrt\pi\cdot 4^{d-s}}{\Gamma(s-d)\Gamma(s-d+\frac12)}.
\nan
So we obtain
\ban
g_\pm(z)&=&\frac{2^{N}}{\sqrt\pi}\suml_{d=0}^\infty\textrm{Coeff}_{s^{2N}}\bigg(s^{2N+1}\Gamma\big(\frac12-(s-d)\big)\Gamma(s-d)^{2N+1}(4e^{\pi\mathbf{i}}z^{-2N})^{d-s}\bigg)\\
&=&\frac{2^{N}}{\sqrt\pi}\suml_{d=0}^\infty\res_{s=0}\Gamma\big(\frac12-(s-d)\big)\Gamma(s-d)^{2N+1}(4e^{\pi\mathbf{i}}z^{-2N})^{d-s}\\
&=&\frac{2^{N}}{\sqrt\pi}\suml_{d=0}^\infty\res_{s=-d}\Gamma\big(\frac12-s\big)\Gamma(s)^{2N+1}(4e^{\pi\mathbf{i}}z^{-2N})^{-s}\\
&=&\frac{2^{N}}{\sqrt\pi}\suml_{d=0}^\infty\res_{s=d}\Gamma\big(\frac12+s\big)\Gamma(-s)^{2N+1}(4e^{\pi\mathbf{i}}z^{-2N})^{s}.
\nan

Take a path $L$ going from $-\mathbf{i}\infty$ to $+\mathbf{i}\infty$, such that $k$ lies to the right of $L$, and $-\frac{1}{2}-k$ lies to the left of $L$, for all nonnegative integers $k$. 
Then by the residue theorem we get
\[
g_\pm(z)=\frac{2^{N}}{\sqrt\pi}\cdot\frac1{2\pi\mathbf{i}}\int_L\Gamma(-s)^{2N+1}\Gamma(\frac12+s)(4e^{\pi\mathbf{i}}z^{-2N})^{s}ds.
\]
Here the convergence issue is addressed in \cite[p. 144\textendash145]{Luk69}. In fact as a consequence  $g_\pm(z)$ is a multiple of a Meijer's $G$-function (we refer readers to \cite[Section 5.2]{Luk69} for the definition and the notation): 
\begin{equation*}
	g_\pm(z)=\frac{2^{N}}{\sqrt\pi}\cdot G^{2N+1,1}_{1,2N+1}\left(4e^{\pi\mathbf{i}}z^{-2N}\bigg|\begin{matrix}\frac12\\\underbrace{0,\dots,0}_{2N+1}\end{matrix}\right).
\end{equation*}
Now the asymptotic expansion of $g_\pm(z)$ follows from \cite[Section 5.7, Theorem 1]{Luk69}.
\end{proof}

\begin{lemma}\label{asymptoticofspm}
	As $z\to0$ in any closed subsector of $|\arg z-\frac\pi{2N}|<\frac\pi2+\frac\pi{2N}$, we have $\cZ^{K}_{\mathbf 0}(S_\pm)(z)\to v_\pm$.
\end{lemma}
\begin{proof}
	From Lemma \ref{f0}, we have
	\ban
	&&\cZ^{K}_{\mathbf 0}(S_\pm)(z)\\
	&=&\bigg[-\frac{g_\pm(z)}{(2\pi z)^N}\hpc^0+\frac12\frac{g_\pm(z)}{(2\pi z)^N}\hpc^{2N}\pm\frac\wp2\bigg]+\suml_{p=1}^{2N}\frac{z^p}{2}\big(\frac1{2\pi z}\big)^{N}\big(\frac{z\partial_z}{2N}\big)^pg_\pm(z)\hpc^{2N-p}\\
	&=&(\mathrm I)+(\mathrm{II}).
	\nan
	In the limit process described in the lemma, one can use Lemma \ref{asymptoticofg} to check that $(\mathrm I)\to v_\pm$. Now it suffices to show that $(\mathrm{II})\to0$, and to this end, we only need to show that for $p\geq1$,
	\ba\label{vanishinginspm}
	\frac{1}{(2\pi z)^N}\big(\frac{z\partial_z}{2N}\big)^pg_\pm(z)=O(1).
	\na
	
	We use Leibnitz rule to find that
	\ban
	&&\frac{1}{(2\pi z)^N}\big(\frac{z\partial_z}{2N}\big)^pg_\pm(z)\\
	&=&\frac{1}{(2\pi z)^N}\big(\frac{z\partial_z}{2N}\big)^p\bigg((2\pi z)^N\cdot\frac{g_\pm(z)}{(2\pi z)^N}\bigg)\\
	&=&\frac{1}{(2\pi z)^N}\suml_{p'=0}^p\binom{p}{p'}\big(\frac{z\partial_z}{2N}\big)^{p-p'}\bigg((2\pi z)^N\bigg)\cdot\big(\frac{z\partial_z}{2N}\big)^{p'}\bigg(\frac{g_\pm(z)}{(2\pi z)^N}\bigg)\\
	&=&\suml_{p'=0}^p\binom{p}{p'}\frac{1}{2^{p-p'}}\big(\frac{z\partial_z}{2N}\big)^{p'}\bigg(\frac{g_\pm(z)}{(2\pi z)^N}\bigg).
	\nan
	From \cite[Section 1.3, (3)]{Luk69}, we can use Lemma \ref{asymptoticofg} to see that for $ p'\geq0$,
	\ban
	\big(\frac{z\partial_z}{2N}\big)^{p'}\bigg(\frac{g_\pm(z)}{(2\pi z)^N}\bigg)\sim-\frac{\mathbf{i}}{2}\delta_{p',0}+\suml_{m=1}^\infty\big(\frac{m}{2N}\big)^{p'}a_mz^m.
	\nan
    This finishes the proof of \eqref{vanishinginspm}.
\end{proof}
\if{
\begin{prop}\label{spmrespect}
	$\cZ^{K}_{\mathbf 0}(S_\pm)$ respects $(0,v_\pm)$ with phase in $(0,\frac\pi{2N})$.
\end{prop}
\begin{proof}
	\textcolor{red}{????????????}
\end{proof}
}\fi
By the proof of Corollary \ref{semisimpleofQ}, $T\zeta^{-k}=2N\cdot 4^{\frac1{2N}}e^{-\frac{k\pi\mathbf{i}}{N}}$ is a simple eigenvalue of $(c_1(\qdc)\circ)$.
\begin{lemma}\label{asymptoticofgk}
	For $0\leq k\leq 2N-1$, as $z\to0$ in any closed subsector of $|\arg z+\frac{k\pi}{N}|<\pi+\frac\pi{2N}$, we have 
	\ban
	e^{\frac{T\zeta^{-k}}z}\frac{g_k(z)}{(2\pi z)^N}\sim\frac{(-1)^k}{2\sqrt{N}}+\suml_{m=1}^\infty b_m^{(k)}z^m,\quad b_m^{(k)}\in\bC.
	\nan
\end{lemma}
\begin{proof}
From \eqref{PQ}, we have
\ban
g_k(z)&=&\suml_{d=0}^\infty\int_\qdc e^{2k\pi\mathbf{i} \hpc}\frac{\Gamma(1+\hpc)^{2N+2}}{\Gamma(1+2\hpc)}\Big[\frac{\Gamma(\hpc-d)}{\Gamma(\hpc)}\Big]^{2N+2}\frac{\Gamma(2\hpc)}{\Gamma(2\hpc-2d)}z^{2N(\hpc-d)}\\
&=&\suml_{d=0}^\infty\textrm{Coeff}_{s^{2N}}\bigg(2e^{2k\pi\mathbf{i} s}\frac{\Gamma(1+s)^{2N+2}}{\Gamma(1+2s)}\Big[\frac{\Gamma(s-d)}{\Gamma(s)}\Big]^{2N+2}\frac{\Gamma(2s)}{\Gamma(2s-2d)}z^{2N(s-d)}\bigg),\\
\nan
where in the second equality, we have used $\int_\qdc\hpc^{2N}=2$. From (\ref{eq-Gammafunction-identity-1}) and
\begin{equation*}
\quad\frac1{\Gamma(2s-2d)}=\frac{2\sqrt\pi\cdot 4^{d-s}}{\Gamma(s-d)\Gamma(s-d+\frac12)},
\end{equation*}
we get
\ban
g_k(z)&=&2\sqrt\pi\suml_{d=0}^\infty\textrm{Coeff}_{s^{2N}}\bigg(e^{2k\pi\mathbf{i} s}s^{2N+1}\frac{\Gamma(s-d)^{2N+1}}{\Gamma(\frac12+s-d)}4^{d-s}z^{2N(s-d)}\bigg)\\
&=&2\sqrt\pi\suml_{d=0}^\infty\textrm{Coeff}_{s^{2N}}\bigg(s^{2N+1}\frac{\Gamma(s-d)^{2N+1}}{\Gamma(\frac12+s-d)}(\frac{e^{2k\pi\mathbf{i}}}4z^{2N})^{s-d}\bigg)\\
&=&2\sqrt\pi\suml_{d=0}^\infty\res_{s=0}\frac{\Gamma(s-d)^{2N+1}}{\Gamma(\frac12+s-d)}(\frac{e^{2k\pi\mathbf{i}}}4z^{2N})^{s-d}\\
&=&2\sqrt\pi\suml_{d=0}^\infty\res_{s=-d}\frac{\Gamma(s)^{2N+1}}{\Gamma(\frac12+s)}(\frac{e^{2k\pi\mathbf{i}}}4z^{2N})^{s}\\
&=&2\sqrt\pi\suml_{d=0}^\infty\res_{s=d}\frac{\Gamma(-s)^{2N+1}}{\Gamma(\frac12-s)}(\frac{e^{2k\pi\mathbf{i}}}4z^{2N})^{-s}.
\nan
Take a path $L$ going from $-\mathbf{i}\infty$ to $+\mathbf{i}\infty$, such that $k$ lies to the right of $L$ for all nonnegative integers $k$. 
Then  as the proof of lemma \ref{asymptoticofg}  we obtain
\begin{eqnarray*}
g_k(z)&=&2\sqrt\pi\cdot\frac{1}{2\pi\mathbf{i}}\int_{c-\mathbf{i}\infty}^{c+\mathbf{i}\infty}\frac{\Gamma(-s)^{2N+1}}{\Gamma(\frac12-s)}(4e^{-2k\pi\mathbf{i}}z^{-2N})^{s}ds\\
&=&2\sqrt\pi\cdot G^{2N+1,0}_{1,2N+1}\left(4e^{-2k\pi\mathbf{i}}z^{-2N}\bigg|\begin{matrix}\frac12\\\underbrace{0,\dots,0}_{2N+1}\end{matrix}\right).
\end{eqnarray*}
Then the asymptotic expansion of $g_k(z)$  follows from \cite[Section 5.7, Theorem 5]{Luk69}.
\end{proof}

\begin{lemma}\label{asymptoticofsk}
For $0\leq k\leq 2N-1$, as $z\to0$ in any closed subsector of $|\arg z+\frac{k\pi}{N}|<\pi+\frac\pi{2N}$, we have 
\ban
e^{\frac{T\zeta^{-k}}z}\cZ^{K}_{\mathbf 0}(\cO(k))(z)\to v_k.
\nan
\end{lemma}
\begin{proof} 
From Lemma \ref{f0}, we have
\ban
e^{\frac{T\zeta^{-k}}z}\cZ^{K}_{\mathbf 0}(\cO(k))(z)&=&\bigg[-e^{\frac{T\zeta^{-k}}z}\frac{g_k(z)}{(2\pi z)^N}\hpc^0+\frac12e^{\frac{T\zeta^{-k}}z}\frac{g_k(z)}{(2\pi z)^N}\hpc^{2N}\bigg]\\
&&\quad+\suml_{p=1}^{2N}\frac{z^p}{2}e^{\frac{T\zeta^{-k}}z}\big(\frac1{2\pi z}\big)^{N}\big(\frac{z\partial_z}{2N}\big)^pg_k(z)\hpc^{2N-p}.
\nan
Note that
\ban
&&\frac{z^p}{2}e^{\frac{T\zeta^{-k}}z}\big(\frac1{2\pi z}\big)^{N}\big(\frac{z\partial_z}{2N}\big)^pg_k(z)\\
&=&\frac{z^p}{2}e^{\frac{T\zeta^{-k}}z}\big(\frac1{2\pi z}\big)^{N}\big(\frac{z\partial_z}{2N}\big)^p\bigg(e^{-\frac{T\zeta^{-k}}z}\cdot(2\pi z)^N\cdot e^{\frac{T\zeta^{-k}}z}\frac{g_k(z)}{(2\pi z)^N}\bigg)\\
&=&\suml_{p_1+p_2+p_3=p}\binom{p}{p_1,p_2,p_3}\frac{z^p}{2}e^{\frac{T\zeta^{-k}}z}\big(\frac{z\partial_z}{2N}\big)^{p_1}\bigg(e^{-\frac{T\zeta^{-k}}z}\bigg)\\
&&\qquad\cdot\big(\frac1{2\pi z}\big)^{N}\big(\frac{z\partial_z}{2N}\big)^{p_2}\bigg((2\pi z)^N\bigg)\cdot\big(\frac{z\partial_z}{2N}\big)^{p_3}\bigg(e^{\frac{T\zeta^{-k}}z}\frac{g_k(z)}{(2\pi z)^N}\bigg).
\nan
We observe that
\ban
\frac{z^p}{2}e^{\frac{T\zeta^{-k}}z}\big(\frac{z\partial_z}{2N}\big)^{p_1}\bigg(e^{-\frac{T\zeta^{-k}}z}\bigg)&=&\frac12\Big(\frac{T\zeta^{-k}}{2N}\Big)^{p_1}z^{p-p_1}+o(z^{p-p_1}),\\
\big(\frac1{2\pi z}\big)^{N}\big(\frac{z\partial_z}{2N}\big)^{p_2}\bigg((2\pi z)^N\bigg)&=&\frac1{2^{p_2}},
\nan
and from \cite[Section 1.3, (3)]{Luk69}, we can use Lemma \ref{asymptoticofgk} to obtain
\ban
\big(\frac{z\partial_z}{2N}\big)^{p_3}\bigg(e^{\frac{T\zeta^{-k}}z}\frac{g_k(z)}{(2\pi z)^N}\bigg)&\sim&\frac{(-1)^k}{2\sqrt{N}}\delta_{p_3,0}+\suml_{m=1}^\infty\big(\frac{m}{2N}\big)^{p_3}b_m^{(k)}z^m.
\nan
So we conclude that 
\ban
\frac{z^p}{2}e^{\frac{T\zeta^{-k}}z}\big(\frac1{2\pi z}\big)^{N}\big(\frac{z\partial_z}{2N}\big)^pg_k(z)=\frac{(-1)^k}{4\sqrt N}\Big(\frac{T\zeta^{-k}}{2N}\Big)^{p}+o(1),\quad p\geq1.
\nan
As a consequence, we have
\ban
&&e^{\frac{T\zeta^{-k}}z}\cZ^{K}_{\mathbf 0}(\cO(k))(z)\\
&=&\bigg[-\frac{(-1)^k}{2\sqrt N}\hpc^0+\frac{(-1)^k}{4\sqrt N}\hpc^{2N}\bigg]+\suml_{p=1}^{2N}\frac{(-1)^k}{4\sqrt N}\Big(\frac{T\zeta^{-k}}{2N}\Big)^{p}\hpc^{2N-p}+o(1)\\
&\to&v_k.
\nan
\end{proof}

\begin{prop}\label{evencase}
    Even-dimensional quadrics satisfy Gamma conjecture II.
\end{prop}
\begin{proof}
The two-dimensional quadric is isomorphic to $\mathbb{P}^1\times \mathbb{P}^1$, thus the Gamma conjecture II is  known from \cite{FaZh19}. For $N\geq2$, by Proposition \ref{prop-evenconvergent}, the quantum cohomology of $\qdc=\qdc_{2N}$ is analytic and semisimple around $\mathbf 0\in H^*(\qdc)$, and $\mathcal D^b(\qdc)$ admits a full exceptional collection
	\ban
	(E_1,\dots,E_{2N+2})=(S_+,S_-,\cO,\dots,\cO(2N-1)).
	\nan
	 Assume that $B$ is an open neighborhood of $\mathbf 0$ properly-chosen with respect to $\{(u_i,\Psi_i)\}_{1\leq i\leq 2N+2}$, such that
	 \ban
	 (u_1(\mathbf 0),\Psi_1(\mathbf0))=(0,v_+),\quad (u_2(\mathbf 0),\Psi_2(\mathbf0))=(0,v_-),\\
	 (u_i(\mathbf 0),\Psi_i(\mathbf 0))=(T\zeta^{-(i-3)},v_{i-3}),3\leq i\leq2N+2.
	 \nan
	 Fix $\phi_\pm\in(0,\frac\pi{2N})$ such that $\phi_+>\phi_-$. Let
	 \ban
	 \phi_1=\phi_+,\quad\phi_2=\phi_-,\quad\phi_i=-\frac{(i-3)\pi}{N},3\leq i\leq2N+2.
	 \nan
	 From Lemma \ref{asymptoticofspm} and \ref{asymptoticofsk}, $\mathcal Z_B^K(E_i)$ respects $(u_i,\Psi_i)$ with phase $\phi_i$ around $\mathbf 0$. So we can choose $\mathbf t_0$ near $\mathbf 0$ with $u_i(\mathbf t_0)\neq u_j(\mathbf t_0)$ for $i\neq j$, such that $\mathcal Z_B^K(E_i)$ respects $(u_i,\Psi_i)$ with phase $\phi_i$ around $\mathbf t_0$. From Theorem \ref{sufficientcondition}, we see that $\qdc$ satisfies Gamma conjecture II.
\end{proof}

\subsection{Odd-dimensional case}

\begin{prop}\label{oddcase}
	$\qdc_{2N+1}$ satisfies Gamma conjecture II.
\end{prop}
\begin{proof}
	The one-dimensional case was known \cite{CDG18, GGI16}. For $N\geq2$, The proof of Gamma conjecture II for $\qdc_{2N+1}$ is similar to the $\qdc_{2N}$ case in the above subsection. In fact the $\qdc_{2N+1}$ case is easier, since there are no primitive classes. We leave the details to interested readers.
\end{proof}

%So the proof of  Theorem \ref{thm-GammaII-quadric} is completed.

\section*{Acknowledgements}
Both authors are grateful to Jianxun Hu for his encouragement and support, and would like to thank Sergey Galkin for his comments. Ke would also like to thank Xingbang Hao, Changzheng Li for helpful discussions, and especially  Di Yang
for sharing his knowledge on  Frobenius manifolds. The authors are partially supported by NSFC grant 11831017, and Ke is also supported by NSFC grants 11890662, 11771461, 11521101 and 11601534, and Hu by NSFC grant 11701579.

\textsc{School of Mathematics, Sun Yat-sen University, Guangzhou 510275, P.R. China}

 \emph{E-mail address:}  huxw06@gmail.com\\

\textsc{School of Mathematics, Sun Yat-sen University, Guangzhou 510275, P.R. China}

\emph{E-mail address:} kehuazh@mail.sysu.edu.cn

\end{document}